\documentclass[12pt,a4paper]{report}

\usepackage{latexsym}                   %Fonts
\usepackage{epsfig}                     %Figures
\usepackage{graphicx}
\usepackage{amssymb, amsmath, amsthm}   %Theorems
\usepackage{enumerate}                  %Enumeration Style
\usepackage{mathrsfs}                   %Math Fonts
\usepackage{german}
\usepackage[plain]{algorithm}
\usepackage{algorithmic}
\usepackage{makeidx}                    %index
\usepackage{array}
\usepackage{longtable}
\usepackage[nottoc]{tocbibind}          %bib, index in table of contents
\usepackage{geometry}										%fuer Layout
\usepackage{verbatim}										%Textdatei wird dargestellt (z.B. Algorithmus)

\newtheorem{thm}{Theorem}[section]
\newtheorem{lem}[thm]{Lemma}

\newtheorem{cor}[thm]{Corollary}

\theoremstyle{definition}
\newtheorem{defi}[thm]{Definition}
\theoremstyle{definition}
\newtheorem{rem}[thm]{Remark}
\theoremstyle{definition}
\newtheorem{example}[thm]{Example}

\newcommand{\bE}{{\mathbb E}}
\newcommand{\bJ}{{\mathbb J}}
\newcommand{\bN}{{\mathbb N}}
\newcommand{\bP}{{\mathbb P}}
\newcommand{\bR}{{\mathbb R}}

\newcommand{\cP}{{\mathcal P}}

\newcommand{\cS}{{\mathcal S}}
\newcommand{\cT}{{\mathcal T}}
\newcommand{\sC}{{\mathscr C}}
\newcommand{\sL}{{\mathscr L}}

\newcommand\ring[1]{\mathaccent23{#1}}
\def\Vr{\ring{V}}

\makeindex

\begin{document}
%%%%%%%%%%%%%%%%%%%%%%%%%%%%%%%%%%%%%%%%%%%%%%%%%%%%%%%

\thispagestyle{empty}

\begin{center}

{\Huge Technische Universit\"{a}t M\"{u}nchen}
\\
\vspace*{1.5cm}
{\huge \sc{ Zentrum Mathematik}}
\\
\vspace*{3cm}
{\Huge {\bf Stochastic Models}} \linebreak \\[1.5mm]
{\Huge {\bf for Speciation Events}}  \linebreak \\[1.5mm]
{\Huge {\bf in Phylogenetic Trees }}
\\
\vspace*{4cm}
{\Large Diplomarbeit}\linebreak \\
{\Large von}\linebreak \\
{\Large Tanja Gernhard}\\
\vspace*{3cm}

\begin{minipage}{11cm}
\begin{center}
{\large
\begin{tabular}{ll}
Aufgabensteller: & Prof. Dr. Rupert Lasser\\
Betreuer:  & Prof. Dr. Mike Steel\\
& \\
Abgabetermin: & 7. April 2006\\
\end{tabular} }
\end{center}
\end{minipage}

\vspace{2cm}
\end{center}

\newpage
\thispagestyle{empty}
\rule{0pt}{14cm}

\newpage

\thispagestyle{empty}
\rule{0pt}{16cm}

Hiermit erkl\"{a}re ich, dass ich die Diplomarbeit selbst\"{a}ndig
angefertigt und nur die angegebenen Quellen und Hilfsmittel
verwendet habe.\\ \vspace{1.0cm}

M\"{u}nchen, den 7. April 2006 \\ \vspace{0.5cm}

\begin{flushright}
...............................................................\\
\textit{Tanja Gernhard}
\end{flushright}

\newpage
\thispagestyle{empty}
\rule{0pt}{14cm}

\pagenumbering{roman}
\setcounter{page}{1}
%%%%%%%%%%%%%%%%%%%%%%%%%%%%%%%%%%%%%%%%%%%%%%%%%%%%%%%%%%%%%%%%%%%%%%%%%%%%%%%
%   Abstracts                                                                 %
%%%%%%%%%%%%%%%%%%%%%%%%%%%%%%%%%%%%%%%%%%%%%%%%%%%%%%%%%%%%%%%%%%%%%%%%%%%%%%%

%\chapter*{Abstract\markboth{Abstract}{Abstract}}
%\addcontentsline{toc}{chapter}{Abstract}
%% Mark 'Abstract' both even and odd markers
%Your abstract goes here.
%%%%%%%%%%%%%%%%%%%%%%%%%%%%%%%%%%%%%%%%%%%%%%%%%%%%%%%%%%%%%%%%%%%%%%%%%%%%%%%
%   Acknowledgements                                                          %
%%%%%%%%%%%%%%%%%%%%%%%%%%%%%%%%%%%%%%%%%%%%%%%%%%%%%%%%%%%%%%%%%%%%%%%%%%%%%%%
\chapter*{Acknowledgements\markboth{Acknowledgements}{Acknowledgements}}
\addcontentsline{toc}{chapter}{Acknowledgements}
%A year ago, I had no idea what my thesis should be about. When I addressed Mike Steel out of the blue about doing a little project with him, he invited me right away to come. That was the start of my involvement in phylogenetics which would not have been possible without Mike's generous help and support. Mike always had an open ear and great advice for any problems I had and strongly supported the idea of extending our little project to a whole thesis. With Mike, I finally found my field of research. Not only would I like to thank Mike for all his scientific support, but also for encouraged me to see the best spots in New Zealand via hiking, kayaking etc. and always had great suggestions! Thanks so much!
%
%
%I'd also like to thank Craig Moritz, Andrew Hugall, Arne Mooers and Rutger Vos who posed the questions which led to my thesis.
%
%Further, I'd like to thank the Friedrich-Ebert-Stiftung for the support throughout
%my time at university and the Allen Wilson Center for having me as a summer student while I was in New Zealand.
%
%Last but not least, I would like to thank my family and Markus for supporting me over the whole time, having good advice whenever I had to made a key decision and who were always encouraging me.
%
%
%
%%
%%
%%
%%
%%When I addressed Mike Steel out of the blue to spend some time in New Zealand, I didn't really
%%

First and foremost, I would like to thank my supervisor Mike Steel for making it possible for me to come to New Zealand, for the great support throughout my stay, for suggesting great problems to work on and
for very helpful discussions and advice.
%, and for making it possible I could give my first talk on a conference and possibly have my first paper published. 
Through my stay in New Zealand and my work with Mike, I finally found my area in research. %For the weekends, Mike always had great suggestions for us visitors, so I saw most spectacular spots in New Zealand. Thanks so much!
 
My thesis abroad and the great experience I had during that time would not have been possible without the support of my German supervisor Rupert Lasser. He encouraged me in any of my plans and let me have all the freedom I needed in choosing a topic for my thesis.

The three days of Daniel Ford's stay in Canterbury were probably the three most productive days of my thesis, while we implemented and optimized my algorithms. Daniel introduced me to Python which was a very convenient language for my research.

Talking to Erick Matsen during coffee breaks helped me to see things I was working on in a broader scientific perspective.
Mareike Fischer had very helpful comments for last improvements of my thesis.

I would also like to thank Craig Moritz, Andrew Hugall, Arne Mooers and Rutger Vos who posed the questions which led to my thesis.

The people and the friendly environment in the Biomath Department at Canterbury University made my stay most enjoyable. Special thanks go to Charles Semple who helped me very much when I first arrived so that I felt comfortable in New Zealand right away.

Further, thanks to the Friedrich-Ebert-Stiftung for the support throughout
my time at university and the Allan Wilson Center for hosting me as a summer student while I was in New Zealand.

Last but not least, I would like to thank my family and my boyfriend for supporting me in any possible way, for giving me good advice whenever I had to make a key decision, for always encouraging me and for providing me a home I always look forward going back to.

\newpage
\thispagestyle{empty}
\rule{0pt}{14cm}

%\addcontentsline{toc}

{\setlength{\baselineskip}{1.5\baselineskip} \tableofcontents }

\newpage
\pagenumbering{arabic} \setcounter{page}{1}

\chapter{Introduction} \label{introduction}
\section{Overview}

\begin{figure}[!b]
\begin{center}
%\resizebox{8cm}{!}{
\includegraphics[scale=0.5]{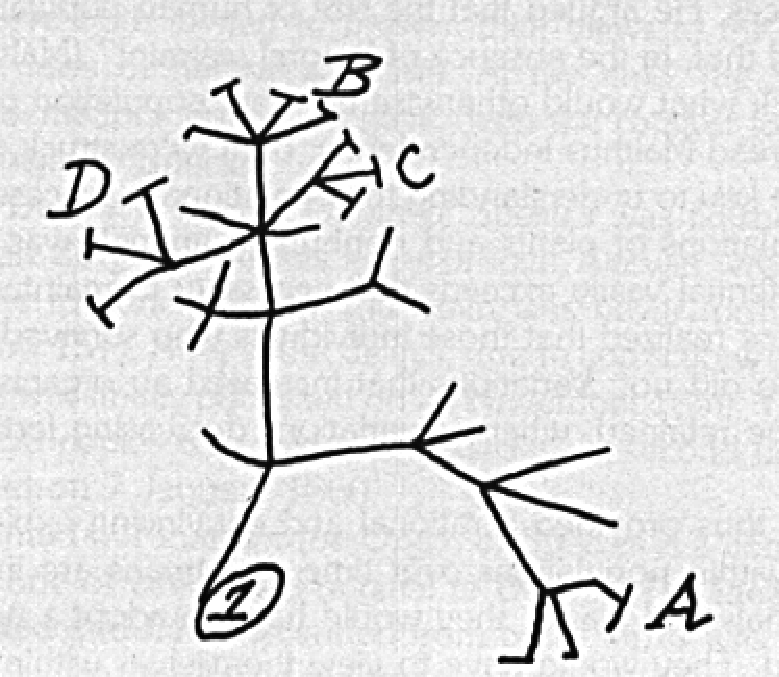}
%}
\caption{Darwin's first diagram of an evolutionary tree from his `First Notebook on Transmutation of Species' (1837).}
\label{FigDarwin}
\end{center}
\end{figure}

In 1837, Darwin published a first sketch of an evolutionary tree, see Fig. \ref{FigDarwin}. This new idea that all species evolved over time was under a lot of discussion and
not until the early 20th century was evolution generally accepted by the scientific community.
Since then, much research went into the field of evolution. With the help
of fossils, and by comparing the anatomy as well as the geographic occurrence of species, complex evolutionary trees have been created.

In an evolutionary tree, each leaf represents an existing species and all the interior vertices represent the ancestors. The edges of the tree show
the relationships between the species.

\begin{figure}
\begin{center}
\includegraphics[scale=0.5]{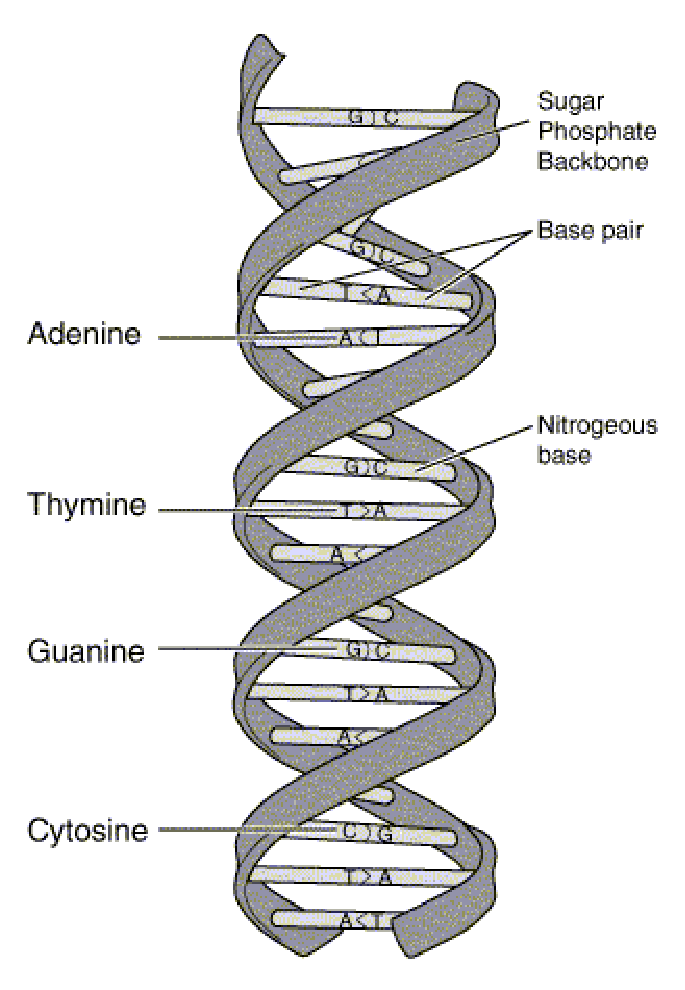}
\caption{The DNA - a double helix}
\label{DNADoubleHelix}
\end{center}
\end{figure}

The first step to modern evolutionary research was the discovery of the double helix structure of DNA (deoxyribonucleic acid) by Watson and Crick in 1953. The genetic code is a long chain of bases (Adenine, Cytosine, Guanine, Thymine) and triplets of these bases
encode the 20 amino acids.
A backbone of sugars and phosphates holds the bases together, see Fig. \ref{DNADoubleHelix}.
The amino acids in a cell form proteins according to the DNA code.
From a chemical point of view, life is nothing else than the functioning of proteins.
Since the DNA determines which proteins are built, a living organism can chemically be described by its DNA, the genetic information \cite{Wikipedia}.

Each cell of an organism has an identical copy of the DNA.
In eukaryotes, the DNA is found in a cell nucleus whereas in prokaryotes (archaea and bacteria), the DNA is not separated from the rest of the cell.

During reproduction, the DNA is transmitted to the offspring, so parents and children are similar in many ways (e.g. hair color,
blood group, disease susceptibility).

\begin{figure}
\begin{center}
\includegraphics[scale=0.5]{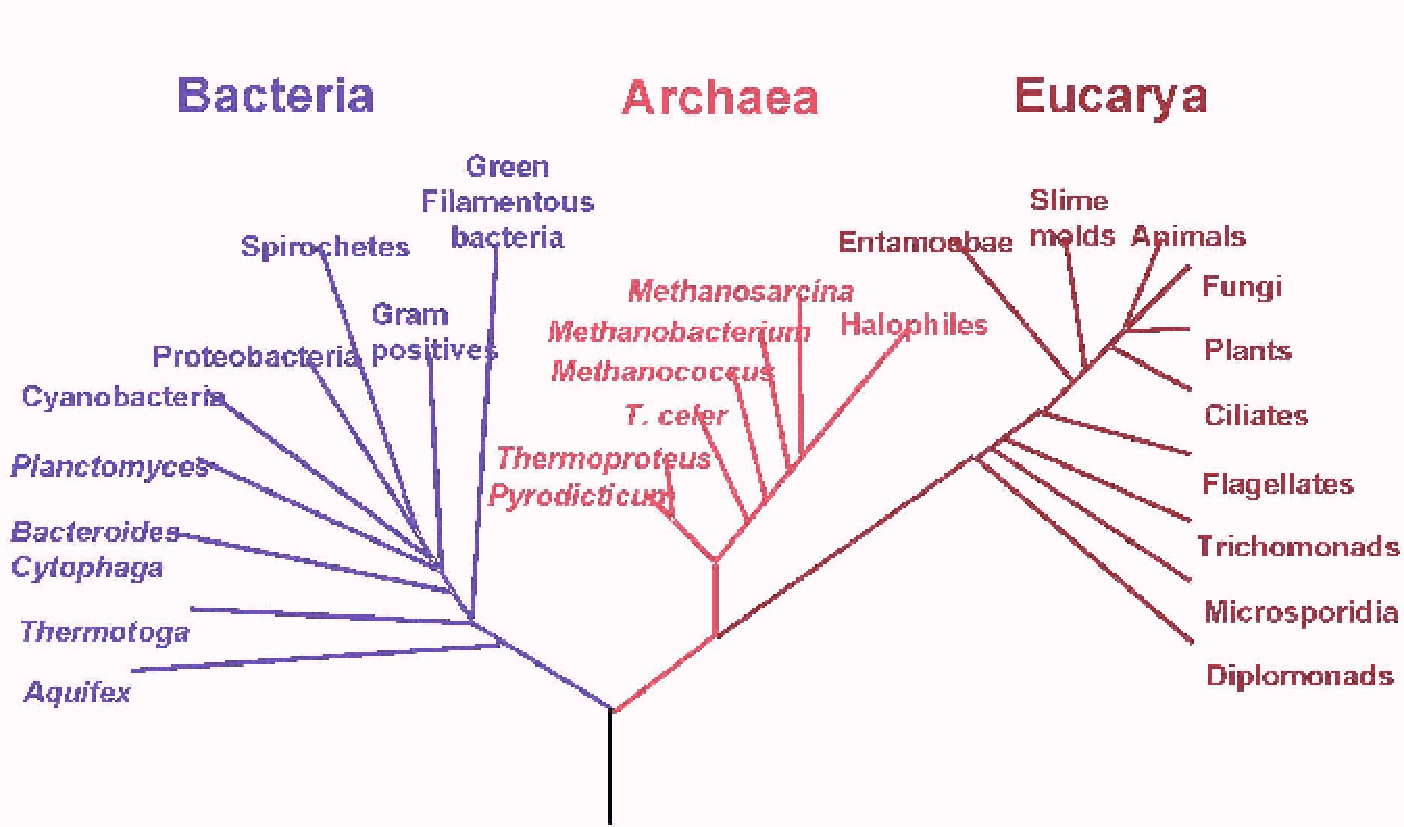}
\caption{Illustration of the tree of life by Carl R. Woese. There are three main branches, the bacteria, archaea and eucarya, source http://www.life.uiuc.edu/micro/faculty/faculty$_{-}$woese.htm.}
\label{FigTreeOfLife}
\end{center}
\end{figure}

It was not until 2003 that the complete human DNA code was described. Currently, the complete DNA sequence of several different species is known
(358 bacteria, 27 archae, 95 eukaryotes, see http://www.ncbi.nih.gov/).
By aligning the DNA of different species,
the similarities and differences of the DNA allow us to reconstruct lineages with more accuracy than before; for an example see Fig. \ref{FigTreeOfLife}.

It is noticeable that the same four DNA bases and the 20 amino acids are found in all organisms. This is strong evidence
for having one common ancestor to all the species.

Evolutionary trees are also called `phylogenetic trees'. If all the species in the tree have a common ancestor,
we call the tree a `rooted tree', the common ancestor is called the `root'.

I take a closer look at rooted phylogenetic trees. The shape of the tree is determined by how speciation occurred.
But since speciation is not understood well and is dependent on historical events which we might never be able to reconstruct,
a stochastic model for speciation is needed.
I investigate the Yule model and the uniform model, two very common models.\\
%We state a statistic to test the Yule model against the uniform model.

In my thesis, I develop the theory with a view to the following applications in biology.

Rutger Vos and Arne Mooers from the Simon Fraser University (Vancouver) recently constructed a supertree for the primates (i.e. lemurs, monkeys, apes and humans) as shown in Appendix \ref{Primates}.

In Section \ref{ExPrimatesYule}, we will see that the primate tree is much more likely to have evolved under the Yule than under the uniform model.

With the supertree method, the shape of the primate tree could be determined, but there was no information about the edge lengths, i.e. the time between speciation events.
In \cite{Vos2006}, edge lengths were estimated by simulations, assuming the (super)tree evolved under the Yule model. The authors concluded by asking for an analytical approach which I develop in Chapter \ref{ChaptRank}.

%By developing an analytical approach, we solved closely related problems like
%For such a tree, we will calculate the probability of one event being earlier than the second.
%This probability depends on the shape of the tree. Our calculation of the probability requires certain properties regarding the stochastic model
%for speciation
%which are in particular fulfilled by the Yule model.

Craig Moritz (UC Berkeley) and Andrew Hugall (University of Adelaide) worked with an evolutionary tree which had edge lengths assigned. The leaves were different types of snails. The snails either live in open forest or rain forest. Moritz and Hugall asked (pers. comm.) if the rate of speciation for open forest snails differs from the rate of speciation for rain forest snails. The rate of speciation is a measure of how fast a class of species produces splits in the evolutionary tree. Chapter \ref{SpeciRate} provides a linear algorithm for solving that problem.

%All those calculations depend a lot on the correctness of the stochastic model for speciation, in our case the Yule model.
%We will see that the trees under the Yule model are more balanced than the trees constructed from biological data. Experiences
%found that the Yule model is still a very good model for speciation though.

\section{Short guide to the thesis}

In Chapter \ref{StochModels}, two important stochastic models for binary phylogenetic trees are introduced - the uniform and the Yule model.
Those two models are discussed and the Kullbach-Liebler-distance between them is calculated.
The Kullbach-Liebler-distance turns out to be very useful in deciding whether a given tree evolved under the Yule or the uniform model.

Chapter \ref{Martingale} formulates a test statistic for that decision problem, the log-likelihood-ratio test. Instead of estimating
the power of the test by simulations, we provide an analytic bound for the power by introducing a martingale process on trees and applying the Azuma inequality.

The algorithms in Chapter \ref{ChaptRank} work in particular for trees under the Yule model. In order to verify that a tree evolved under Yule, the test provided in Chapter \ref{Martingale} can be applied before running the algorithms.

After having established all the necessary stochastic background, Chapter \ref{ChaptRank} provides a quadratic algorithm
for calculating the probability distribution of the rank
for a given interior vertex in a phylogenetic tree. The algorithm is called {\sc RankProb} and we assume that every rank function
on a given tree is equally likely.
That is in particular the case for the Yule model. The algorithm {\sc RankProb} is extended to non-binary trees as well, again we
assume that every rank function is equally likely. We call that algorithm {\sc RankProbGen}.
Calculating the probability
of having an interior vertex $u$ earlier in the tree than an interior vertex $v$ is calculated with the algorithm {\sc Compare} in quadratic time.
We coded up the algorithms {\sc RankProb} and {\sc Compare} in Python, see Appendix \ref{PythonCode}.
The chapter concludes with an analytical approach of estimating edge lengths in a given tree under the Yule model. This approach makes use of the algorithm {\sc RankProb}.

Chapter \ref{SpeciRate} looks at the rate of speciation. Given is a phylogenetic tree with the leaves being
divided into two classes $\alpha$ and $\beta$.
The edge lengths shall represent the time between two events.
We provide a linear algorithm for the expected time a species of class $\alpha$ exists until it speciates and two new species evolve. The average edge length is an estimate for the inverse of the rate of speciation.
An example for the classes $\alpha$ and $\beta$ could be rain forest snails and open forest snails.\\

%The trees in Appendix \ref{primates} were supplied by Rutger Vos and Arne Mooers. The coding of the algorithms in 

After introducing the stochastic models in Chapter \ref{StochModels}, the remaining results in that Chapter are new. The results in Chapter
\ref{Martingale}, \ref{ChaptRank} and \ref{SpeciRate} are new unless otherwise stated. Improvements on the algorithms in Chapter \ref{ChaptRank} and coding them up in Python was joint work with Daniel Ford.
Chapter \ref{ChaptRank} was the topic of my talk at the New Zealand Phylogenetics Conference in Kaikoura in February 2006 (http://www.math.canterbury.ac.nz/bio/kaikoura06/).\\% and parts of this chapter are published in .

The rest of this Chapter introduces the basic definitions from graph theory and phylogenetics needed for the thesis. Further,
some basic results for phylogenetic trees are stated.

\section{Graphs and Trees} \label{GraphsAndTrees}

\begin{defi} \index{graph} \label{DefGraph}
A $graph~G$ is an ordered pair $(V,E)$ consisting of a non-empty
set $V$ of $vertices$ \index{vertex} and a multiset $E$ of $edges$ \index{edge} each of which
is an element of $\{ \{x,y\}: x,y \in V \}$. The degree $\delta(v)$ \index{vertex!degree of} \label{DefDegreeV} of
a vertex $v \in V$ is the number of edges in $G$ that are incident
with $v$.
A $path$ \index{path} $p$ in $G$ from vertex $x \in V$ to vertex $y \in V$ is a sequence $p=(v_i)_{i=1, \ldots n}$, $v_i \in V$,
such that $x=v_1$, $y=v_n$, and $\{v_i, v_{i+1}\} \in E$ for $i=1, \ldots n-1$.
A graph $G$ is $connected$ \index{graph!connected} precisely if there exists a path from $x$ to $y$ for all $x, y \in V$.
A $cycle$ \index{cycle} in a graph is a path $p=(v_i)_{i=1, \ldots n}$ with $v_1 = v_n$.
The graph $G'=(V',E')$ is a $subgraph$ \index{subgraph}of $G$ if $V' \subseteq V$ and $E' \subseteq E$.
\end{defi}

\begin{figure}
\begin{center}
\input{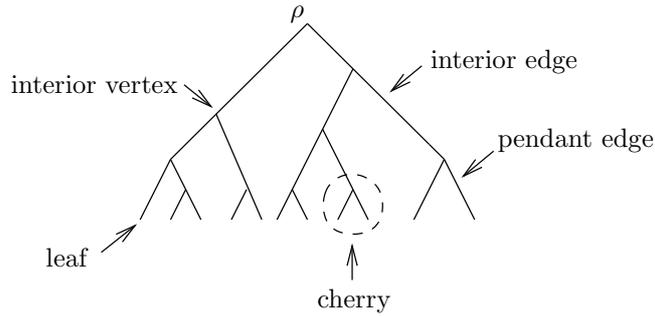}
\caption{A rooted binary tree}
\label{FigTree1}
\end{center}
\end{figure}

\begin{defi} \index{tree} \label{DefTree}
A $tree$ $T=(V,E)$ is a connected graph with no cycles. A connected subgraph of $T$ is a $subtree$ \index{subtree} of $T$.
A $rooted$ $tree$ \index{tree!rooted} is a tree that has exactly one distinguished vertex called the $root$ \index{root} which we denote by the letter $\rho$.
A vertex $v \in V$ with $\delta(v) \leq 1$ is called a $leaf$ \index{leaf}. The set of all leaves of $T$ is denoted by $L$.
A vertex which is not a leaf is called an $interior~vertex$. \index{vertex!interior}
Let $\Vr$ denote the set of all interior vertices of $T$. A $binary$ $tree$ \index{tree!binary} is a tree with $\delta(v)=3$ for all $v \in \Vr$.
A $rooted$ $binary$ $tree$ is a rooted tree with $\delta(v)=3$ for all $v \in \Vr \setminus \rho$ and $\delta(\rho)=2$.
Let $V' \subset V$. The subtree $T'=T|_{V'}$ \index{subtree} is the minimal (w.r.t. the number of vertices) connected subgraph of $T$ containing $V'$.
An edge which is incident with a leaf is called a $pendant$ edge. \index{edge!pendant} A non-pendant edge is called an $interior$ edge. \index{edge!interior}
Two distinct leaves of a tree form a $cherry$ \index{cherry} if they are adjacent to a common ancestor.
Let $v \in \Vr \setminus \rho$ with $\delta(v)=2$.
The vertex $v$ is $suppressed$ \index{vertex!suppressed} in $T$ if we delete $v$ with its two incident edges $e_1=(v_1,v), e_2=(v,v_2)$ and then add a new edge $e=(v_1,v_2)$.
For an example of a tree see Fig. \ref{FigTree1}.
\end{defi}

\begin{defi}
\index{phylogenetic $X$-tree} \index{phylogenetic tree} \label{DefPhyloTree}
Let $T=(V,E)$ be a rooted tree with leaf set $L \subset V$ and for all $v \in \Vr \setminus \rho$ is $\delta(v) \neq 2$.
Let $X$ be a non-empty finite set with $|X|=|L|$.
Let $\phi: X \rightarrow L$ be a bijection. Then $\cT=(T, \phi)$ is called a $phylogenetic$ $(X-)~tree$ with $labeling$ $function$ $\phi.$ \index{labeling function}
$X$ is called the {\it label set.} \index{label set}
A phylogenetic tree is also called a $labeled$ $tree$. \index{labeled tree} A $tree$ $shape$ \index{tree shape} is a phylogenetic tree without the labeling.
\end{defi}
\begin{rem}
In the following, for a phylogenetic tree $\cT$, we sometimes write $E_{\cT}$ instead of $E$,
$V_{\cT}$ instead of $V$,
$\Vr_{\cT}$ instead of $\Vr$
and $L_{\cT}$ instead of $L$. This notation clarifies to which tree the sets refer whenever we talk about several different trees.
\end{rem}
\begin{defi} \label{DefPartOrder} \index{partial order on a tree}
Let $T$ be a rooted tree. A partial order $\leq_{T}$ on V is obtained by setting $v_1 \leq_{T} v_2$ ($v_1, v_2 \in V$)
precisely if the path from the root $\rho$ to $v_2$ includes $v_1$. If $v_1 \leq_{T} v_2$,
we say $v_2$ is a $descendant$ \index{descendant} of $v_1$ and $v_1$ is an $ancestor$ \index{ancestor} of $v_2$.
If $v_1 \leq_{T} v_2$ and there is no $v_3 \in V$ with
$v_1 \leq_{T} v_3 \leq_{T} v_2$, we say $v_2$ is a $direct~descendant$ \index{descendant!direct} of $v_1$ and $v_1$ is a $direct~ancestor$ \index{ancestor!direct} of $v_2$.
The {\it number of direct descendants} of $v$ is $d(v)$. \label{DefNumbDes}
When we talk about a phylogenetic tree, we often write $\leq_{\cT}$ instead of $\leq_T$.
\end{defi}
\begin{figure}
\begin{center}
\input{treeintro02-2.pstex_t}
\caption[A rooted binary phylogenetic $X$-tree]{A rooted binary phylogenetic $X$-tree $\cT$ with
$X=\{a,b, \ldots ,k \}$ and the subtree $\cT'=\cT|_{\{f,h,i,j \}}$.} \label{FigTree2-2}
\end{center}
\end{figure}
\begin{defi} \label{DefSubtreeX'}
Let $\cT=(T,\phi)$ be a phylogenetic $X$-tree.
Let $X' \subset X$. The phylogenetic subtree $\cT'=\cT|_{X'}=(T',\phi')$ \index{subtree!phylogenetic}
is a phylogenetic tree where $T'$ is the tree $T|_{\phi(X')}$ with all degree-two vertices suppressed (except for the root).
The labeling function is $\phi' = \phi |_{X'}$.
The root of $\cT'$ is the vertex $\rho'$ which is minimal in the tree $T'$ under the partial order $\leq_{\cT}$ (see Fig. \ref{FigTree2-2}).
%Let $X = \{1, \ldots, n\}$ and $X' = \{1, \ldots, i\}, i \leq n$. We write $\cT|_{ \{1, \ldots, i\} }$ instead of $\cT|X'$.
Let $\cT'$ be a subtree of $\cT$. Denote the subtree $\cT | _{L_{\cT} \setminus L_{\cT'}}$ by $\cT \setminus \cT'.$

Let $v \in \Vr$ and let $X_v$ be the label set of all the leaves in $\cT$ which are descendants of $v$.
The subtree $\cT_v$ is $induced$ $by$ $v$ if $\cT_v=\cT|_{X_v}$. \label{DefSubtreev} \index{subtree!induced by $v$}
A binary phylogenetic tree is $balanced$ \index{balanced tree} if the two subtrees induced by the two direct descendants of the root have the same shape.
Otherwise, the tree is $unbalanced$. \index{unbalanced tree}
\end{defi}

\begin{figure}
\begin{center}
\input{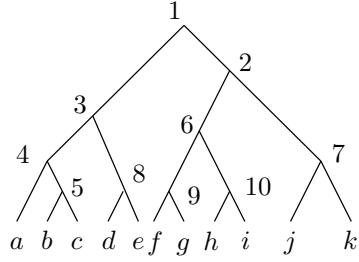}
\caption[A rooted binary phylogenetic ranked $X$-tree]{A rooted binary phylogenetic ranked $X$-tree with
$X=\{a,b, \ldots ,k \}$} \label{FigTree2}
\end{center}
\end{figure}

\begin{defi} \label{DefRank}
Let $\cT$ be a rooted phylogenetic tree. Let the function $r$ be a
bijection from the set of interior vertices $\Vr$ of $\cT$ into
$\{1,2, \ldots ,|\Vr|\}$ that satisfies the following property:
$$if~v_1\leq_{\cT} v_2,~then~r(v_1) \leq r(v_2)$$
$(\cT,r)$ is called a {\it phylogenetic ranked tree} (see Fig. \ref{FigTree2}). \index{phylogenetic tree!ranked} \index{ranked phylogenetic tree}
The function $r$ is called a {\it rank function} for $\cT$. \index{rank function}
A vertex $v$ with $r(v)=i$ is said to be in the $i-th~position$ of $\cT$ or $v$ has rank $i$.
We write $r_{\cT}$ instead of $r$ when it is not clear from the context to which tree the rank function $r$ refers.
Note that $r$ induces a linear order on the set $\Vr$. We define the set $r(\cT)$ as
$$r(\cT)= \{r: \ r \rm{ \ is \ a \ rank \ function \ on} \ \cT \}.$$
\end{defi}

The following Lemma has been shown in \cite{Steel2003} using poset theory. We will give an elementary proof using induction.
\begin{lem} \label{LemNumbRank}
Let $\cT$ be a rooted phylogenetic tree. For each $v \in \Vr$, let $\lambda_v$ denote the number of elements of $\Vr$ that are descendants of $v$.
Then the number of rank functions for $\cT$ is \index{rank function!number of}
\begin{equation} |r(\cT)| = \frac{|\Vr|!}{\displaystyle \prod_{v \in \Vr}\lambda_v} \label{EqnNumbRank} \end{equation}
Note that a vertex $v$ is a descendant of itself by definition, so $\lambda_v$ also counts the vertex $v$.
\end{lem}
\begin{proof}
This proof is done by induction over the number $n$ of interior vertices of a tree.
For $n=1$, there is only one rank function, the only interior vertex has rank $1$, which equals to $\frac{|\Vr|!}{\prod_{v \in \Vr}\lambda_v}=\frac{1!}{1}=1$.
Suppose that (\ref{EqnNumbRank}) is true for all trees with $n<k$ interior vertices.
Let $\cT$ be a tree with $k$ interior vertices. The degree of root $\rho$ is $\delta(\rho)=m$ where $m<k$.
$\cT$ has $m$ vertex-disjoint rooted subtrees $\cT_1,\cT_2, \ldots ,\cT_m$ induced by the direct descendants of $\rho$, and with $|\Vr_{\cT_i}|<k$.
Each subtree $\cT_i$ has $\frac{|\Vr_{\cT_i}|!}{\prod_{v \in \Vr_{\cT_i}}\lambda_v}$ different rank functions by the induction assumption.
Counting all the rank functions on $\cT$ is equivalent to counting the rank functions on each subtree $\cT_i$
and then combining the positions of the vertices of all the $\cT_i$ to get a linear order on $\Vr_{\cT}$, by preserving the order of the vertices of each $\cT_i$.
For a given rank function on each $\cT_i$, we can order all the interior vertices in $\frac{\left(\sum_i |\Vr_{\cT_i}| \right)!}{\prod_i \left(|\Vr_{\cT_i}|! \right)}$
different ways where the order within each $\cT_i$ is preserved. Multiplying by all the possible rank functions for each $\cT_i$ yields to
\begin{eqnarray}
|r(\cT)| &=& \frac{\displaystyle \left( \sum_{i=1}^m |\Vr_{\cT_i}| \right)!}{\displaystyle \prod_{i=1}^m \left(|\Vr_{\cT_i}|! \right)}
\left( \prod_{i=1}^m |r(\cT_i)| \right)  \nonumber \\
&=&
\frac{\displaystyle \left( \sum_{i=1}^m |\Vr_{\cT_i}| \right)!}{\displaystyle \prod_{i=1}^m \left(|\Vr_{\cT_i}|! \right)}
\left( \prod_{i=1}^m \frac{|\Vr_{\cT_i}|!}{\prod_{v \in \Vr_{\cT_i}}\lambda_v} \right)  \nonumber \\
& = & \left(\sum_{i=1}^m |\Vr_{\cT_i}| \right)!
\prod_{i=1}^m \frac{1}{\prod_{v \in \Vr_{\cT_i}}\lambda_v}  \nonumber \\
 & = & \frac{(|\Vr_{\cT}|-1)!}{\displaystyle \prod_{v \in \Vr_{\cT \setminus \rho}}\lambda_v} \notag \\
 & = & \frac{|\Vr_{\cT}|!}{\displaystyle \prod_{v \in \Vr_{\cT}}\lambda_v}. \nonumber
\end{eqnarray}
This establishes the induction step, and thereby the theorem.
\end{proof}
\begin{rem} \label{RemTreeSet}
In the following, all trees shall be rooted.
The set of all binary rooted phylogenetic trees with label set $X$ is denoted by $RB(X)$.
The set of all ranked binary rooted phylogenetic trees with label set $X$ is denoted by $rRB(X)$.
\end{rem}
\begin{rem} \label{RemBinTree}
A rooted binary phylogenetic tree with $n$ leaves has $|\Vr| = n-1$ interior vertices and $|E|=2(n-1)$ edges,
which is shown by induction in \cite{Steel2003}.
\end{rem}

\chapter{Stochastic Models on Trees} \label{StochModels}

Given a phylogenetic $X$-tree, we are interested in the probability of that tree from the set $RB(X)$ or $rRB(X)$, depending on whether the given
tree is ranked or not.
When defining a probability distribution on trees,
the probability of a labeled tree should be invariant under a different labeling.
This property is called \emph{exchangeability.} \index{exchangeabilitiy}

There are several stochastic models for binary phylogenetic $X$-trees, the
most common are the uniform and Yule model which we will introduce and compare.

In the following, for simplifying notation, any $X$ with $|X|=n$ shall be $X= \{1, 2, \ldots ,n\}$ and we write $RB(n)$, $rRB(n)$ instead of $RB(X)$, $rRB(X)$. \label{TreeSetn}
%Further Comb model, A generalization of those three models is the alpha model.

\section{The uniform model} 
\begin{figure}
\begin{center}
\input{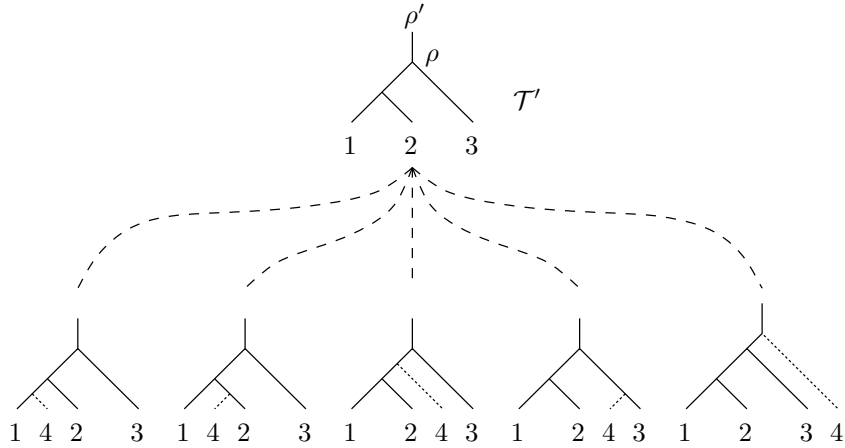}
\caption[Tree evolving under the uniform model]
{Tree evolving under the uniform model. Let $X=\{1,2,3,4\}$.
Given the tree $\cT'$ with label set $\{1,2,3\}$, which has probability $1/3$ under the uniform model, there are five possible edges to attach the leaf with label $4$.
Each of the five trees with label set $\{1,2,3,4\}$ has probability $1/5$ given $\cT'$.
So the overall probability of each tree with four leaves is $1/15$ under the uniform model.}
\label{FigTreest01}
\end{center}
\end{figure}

\index{uniform model}
Under the uniform model, a random element of $RB(n)$ is generated in the following way ({\it cf.} Figure \ref{FigTreest01}):
\begin{itemize}
\item Label the two leaves of a cherry with $1$ and $2$.
\item Add to the cherry a third edge connecting the root $\rho$ of the cherry and a new vertex $\rho'$ which is earlier than $\rho$. This extended cherry is denoted by $\cT$.
\item In each step, modify $\cT$ in the following way, until $\cT$ has $n$ leaves:
\begin{itemize}
\item Let the number of leaves of $\cT$ be $k$.
Choose an edge of $\cT$ randomly and with uniform probability and subdivide this edge to create a new vertex.
\item Add an edge from the new vertex to a new leaf.
\item Label the new leaf by $k+1$.
\end{itemize}
\item Remove from the tree $\cT$ the vertex $\rho'$ and its incident edge to get the binary rooted tree $\cT$.
\end{itemize}
In this way, each rooted binary phylogenetic $X$-tree has equal
probability (see \cite{Pinelis2003}). Obviously, the probability of a tree is invariant
under a different leaf labeling.
Note that it is not necessary to choose the elements of $X$ in the given order $1, 2, \ldots ,n$. We could choose the leaf labels
in any order. This will not be the case for the Yule model.
\begin{lem} \label{Lem-n!!}
For each $n \geq 2$,
$$(2n-3)!! = \frac{n! c_{n-1}} {2^{n-1}}$$
with $(2n-3)!! = (2n-3) \cdot (2n-5) \ldots 5 \cdot 3 \cdot 1$ and $c_n$ being the \label{DefCatalan} $n$-th Catalan number, $c_n = \frac {1} {n+1} {2n \choose n}$. \index{Catalan number}
\end{lem}
\begin{proof}
\begin{eqnarray}
{(2n-3)!!}
&=& \frac {(2n-3)!} {2^{n-2} \left( \frac{2n-4} {2} \right)!}
= \frac {(2n-3)!} {2^{n-2} (n-2)!}  \notag \\
&=& \frac {(2n-2)!} {2^{n-1} (n-1)!}
= \frac {\frac {(n-1)! (2(n-1))!} {2(n-1)!}}  {2^{n-1}}
= \frac  {n! \frac {1} {n} {2(n-1) \choose n-1}} {2^{n-1}} \notag \\
&=& \frac {n!c_{n-1}}{2^{n-1}}. \notag
\end{eqnarray}
\end{proof}
The following result is already shown in
\cite{Steel2003} by considering unrooted trees and defining a
bijection from unrooted to rooted trees. This proof is direct.
\begin{thm} \label{ThmNumbPhyloTree}
The number of binary rooted phylogenetic trees is $$|RB(n)| = (2n-3)!!$$ \index{phylogenetic tree!number of}
\end{thm}
\begin{proof}
The proof is done by induction over $n$. For $n=2$, we have $|RB(2)|=1$ and $(2 \cdot
2 - 3)!!=1$. Assume $|RB(n)| = (2n-3)!!$ holds for all $n \leq
k$, where $k \geq 2$. A tree $\cT_k$ with $k$ leaves has $2(k-1)$ edges (see Remark
(\ref{RemBinTree})). Denote the root of $\cT_k$ by $\rho_k$. The
$(k+1)$-th leaf $x$ can be attached to $\cT_k$ to any of the
$2(k-1)$ edges or a new root $\rho$ with edges $e_1=(\rho,\rho_k)$
and $e_2=(\rho, x)$ is added. So we can construct $2(k-1)+1 = 2k-1$
different trees from $\cT_k$. By the induction assumption, we have
$|RB(k)| = (2k-3)!!$. Therefore, $|RB(k+1)| = (2k-3)!! \cdot
(2k-1) = (2(k+1)-3)!!$ which proves the theorem.
\end{proof}
\begin{cor} \label{CorProbUnif}
Under the uniform model, the probability $\bP[\cT]$ of a tree $\cT$ chosen from the set $RB(n)$ is
$$\bP[\cT] = \frac{1} {(2n-3)!!} = \frac{2^{n-1}} {n!c_{n-1}}.$$ \index{uniform model!probability of $\cT$}
\end{cor}
\begin{proof}
Since a phylogenetic tree $\cT$ is chosen from $RB(n)$ uniformly at random in the uniform model, we have
$$\bP[\cT] = \frac {1} {|RB(n)|}.$$
By Theorem (\ref{ThmNumbPhyloTree}) and Lemma (\ref{Lem-n!!}), we get $\bP[\cT] = \frac{1} {(2n-3)!!} = \frac{2^{n-1}} {n!c_{n-1}}$.
\end{proof}
\section{The Yule model} 
\begin{figure}
\begin{center}
\input{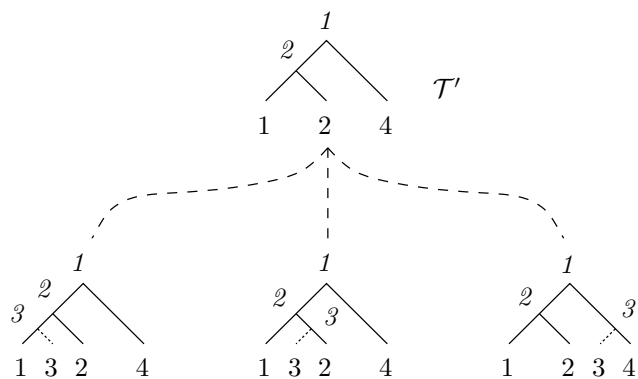}
\caption[Tree evolving under the Yule model]
{Ranked tree evolving under the Yule model. Let $X=\{1,2,3,4\}$.
Suppose the ranked tree $\cT'$ with label set $\{1,2,4\}$ evolved under the Yule model.
There are three possible pendant edges to attach the leaf with the remaining label $3$.
Each ranked tree with label set $\{1,2,3,4\}$ has probability $\frac{2^{4-1}}{4!(4-1)!}=1/18$ according to Theorem (\ref{ThmProbRankTree}).}
\label{FigTreest02}
\end{center}
\end{figure}

\index{Yule model}
Under the Yule model \cite{Yule1924,Harding1971}, a random element of $rRB(n)$ is generated in the following way ({\it cf.} Figure \ref{FigTreest02}):
\begin{itemize}
\item Two elements of $X$ are selected uniformly at random and the two leaves of a cherry are labeled by them.
This cherry is denoted by $\cT$ and its root has rank $1$.
\item In each step, modify $\cT$ in the following way, until $\cT$ has $n$ leaves:
\begin{itemize}
\item Let the number of leaves of $\cT$ be $k$.
Choose a pendant edge of $\cT$ uniformly at random and subdivide this edge to create a new interior vertex with rank $k$.
\item Add an edge from the new vertex to a new leaf.
\item Select an element of $X$ which is not in the label set of $\cT$ uniformly at random and label the new leaf by that element.
\end{itemize}
\end{itemize}
In other words, any pendant edge of a binary tree is equally likely to split and give birth to two new pendant edges.
The Yule model is therefore an explicit model of the process of speciation.
This makes it a very important model for the distribution on trees. Since the labels are added uniformly at random, the probability of a tree
is invariant under a different leaf labelling (i.e. dependent only on the `shape' of the tree).

Note that under the Yule model, at each moment in time, the probability of a speciation event is equal for all the current species.
For different points in time, these probabilities can be quite different though.

Under the Yule model, balanced trees are more likely than unbalanced trees whereas under the uniform model, every tree is equally likely.
Phylogenetic trees constructed for most sets of species tend to be more balanced than predicted by the uniform model,
but less balanced than predicted by the Yule model.
That can be explained in the following way.
In nature, we observe that a species, which has not given birth to new species for a long time, is not very likely to give birth in the future either.
The Yule model does not take this fact into account.
In \cite{Steel2001}, there is an extension of the Yule model described which takes care of that biological observation.
One special case of the extended Yule model assumes, that unless a species has undergone a speciation event within the last $\epsilon$ time interval,
it will never do so.
It is shown in \cite{Steel2001} that for sufficient small $\epsilon$, this model induces the uniform distribution.
So the uniform model can also be interpreted as a process of speciation.

The Yule and the uniform model can be put in a more general framework. In \cite{Aldous1996}, the beta-splitting model is introduced,
where the Yule and the uniform model are special cases. In \cite{Ford2005},
the alpha model is introduced and again, the Yule and the uniform model are special cases.
In both papers, a one parameter family of probability models on binary phylogenetic trees is introduced which interpolates continuously between
the Yule and the uniform model.

These models are far more complicated than the uniform and Yule model though, and since especially the Yule model is still a
reasonably good model for speciation, we will now focus on properties of the Yule model.
Theorem (\ref{ThmProbRankTree}) and Corollary (\ref{CorNumbRankTree})
have been established in \cite{Edwards1970}. Here we provide an alternative proof.
\begin{thm} \label{ThmProbRankTree}
The probability under the Yule model of generating a ranked binary phylogenetic tree $(\cT,r) \in rRB(n)$ is
\index{Yule model!probability of $(\cT,r)$}
$$\bP [\cT,r] = \frac{2^{n-1}}{n!(n-1)!}.$$
That is a uniform distribution over $rRB(n)$.
\end{thm}
\begin{proof}
We calculate the probability $\bP [\cT,r]$ by looking at the generation of the tree $\cT$.
%By ranking the interior vertices of $\cT$ according to when they were created in the generation of $\cT$, we get a ranked tree $(\cT,r)$.
In the first step of the generation, we have
%${n \choose 2} = \frac {n(n-1)} {2} $
$n$
possibilities to choose
the label for the left leaf of the cherry and $n-1$ possibilities to choose the label for the right leaf of the cherry.
So the probability for a certain cherry, with distinguishing between left and right vertex,
is $\frac {1} {n(n-1)}$, since the selection of the labels is uniformly at random.
The root of the cherry has rank $1$.
When adding a new leaf to a tree $\cT_k$ with $k$ leaves, we have $k$ possibilities to choose a pendant vertex and $n-k$
possibilities to choose a label.
So the probability of attaching a new labeled leaf to a certain edge is
$\frac{1} {k (n-k)}$ since we choose the pendant edge and the label uniformly at random.
The new interior vertex has rank $k$.
Let the new leaf be $x$. The leaf $x$ shall be on the right side of the new cherry.
With the process above, we get two equal trees precisely if every step of the tree generation process is equal for both trees.
While distinguishing between left and right child of an interior vertex, we count each phylogenetic tree $2^{|\Vr|}=2^{n-1}$ times.
Therefore, we get the following probability for the ranked phylogenetic tree $(\cT,r)$
$$\bP[\cT,r] = 2^{n-1} \frac {1} {n(n-1)} \frac {1}{2 (n-2)} \frac {1}{3 (n-3)} \ldots \frac {1}{(n-1) 1} = \frac{2^{n-1}}{n!(n-1)!}$$
Since $\bP[\cT,r]$ is independent of $\cT$ and $r$, we have a uniform distribution.
\end{proof}
\begin{cor} \label{CorNumbRankTree}
The number of ranked phylogenetic trees is
$$|rRB(n)|=\frac{n!(n-1)!}{2^{n-1}}$$ \index{ranked phylogenetic tree!number of}
\end{cor}
\begin{proof}
Since $\bP[\cT,r] = \frac{2^{n-1}}{n!(n-1)!}$ is uniform under the Yule model and probabilities add up to $1$,
we have $\frac{n!(n-1)!}{2^{n-1}}$ different ranked phylogenetic trees.
\end{proof}
\begin{lem} \label{LemYule}
Let $A$ be a finite set and for each $a \in A$, let $B(a)$ be a finite set and let $\Omega=\{(a,b): a \in A, b \in B(a) \}$.
Let $C=(C_1,C_2)$ be the (two-dimensional) random variable which takes a value in $\Omega$ selected uniformly at random,
i.e. $\bP[C=(a,b)] = {1}/{|\Omega|}$ for all $(a,b) \in \Omega$.
Then the conditional probability distribution $\bP[C=(a,b)|C_1=a]$ is uniform on $B(a)$.
\end{lem}
\begin{proof}
We have
$$\bP[C=(a,b)|C_1=a] = \frac{\bP[C=(a,b)]}{\bP[C_1=a]} = \frac{1}{|\Omega| \bP[C_1=a]}$$
which is independent of $b$ and therefore is uniform on $B(a)$.
\end{proof}
\begin{thm} \label{ThmYuleRankGivenT}
Assume a given binary phylogenetic tree $\cT$ with $n$ leaves evolved under the Yule model.
Then the probability of a rank function $r$ on a given tree $\cT$ is
\index{Yule model!probability of $r$ given $\cT$}
$$\bP[r | \cT] = \frac{\prod_{v \in \Vr}\lambda_v} {(n-1)!}$$
i.e. $\bP[r | \cT]$ is uniform over all rankings $r$ of $\cT$.
\end{thm}
\begin{proof}
Consider the probability distribution induced by the Yule model on $A=RB(n)$.
Let $B(a)$ be the set of all rankings for a tree $a \in A$
and let $\Omega=\{(a,b): a \in A, b \in B(a) \}$.
Let $C=(C_1,C_2)$ be the (two-dimensional) random variable which takes a value in $\Omega$.
The random variable $C$ is uniform on the set $\Omega$ by Theorem (\ref{ThmProbRankTree}) and we can apply Lemma (\ref{LemYule}) to obtain
$$\bP[C=(\cT,r)|C_1=\cT] = \bP[r | \cT] = \frac{1}{|\Omega| \bP[C_1=\cT]}$$
which shows that $\bP[r | \cT]$ is uniform over all rankings $r$ of $\cT$.
Since for a tree $\cT$, we have $\frac{|\Vr|!}{ \prod_{v \in \Vr}\lambda_v}$ possible rankings by (\ref{LemNumbRank}), and $|\Vr|=n-1$ for binary trees,
we get
$$\bP[r | \cT] = \frac {1} {\frac{|\Vr|!}{ \prod_{v \in \Vr}\lambda_v}} = \frac{\prod_{v \in \Vr}\lambda_v} {(n-1)!}.$$
\end{proof}
The following Corollary was established in \cite{Brown1994} using induction.
\begin{cor} \label{CorProbYule}
The probability of a binary phylogenetic tree $\cT \in RB(n)$ under the Yule model is
\index{Yule model!probability of $\cT$}
$$\bP[\cT] = \frac{2^{n-1}} {\displaystyle n! \prod_{v \in \Vr} \lambda_v}$$
where $\lambda_v$ is as defined in Lemma (\ref{LemNumbRank}).
\end{cor}
\begin{proof}
With Theorem (\ref{ThmProbRankTree}) and Theorem (\ref{ThmYuleRankGivenT}) we get
$$\bP[\cT] = \frac{\bP[\cT,r]}{\bP[r|\cT]} =  \frac{2^{n-1}}{n!(n-1)!} \cdot \frac {(n-1)!} {\prod_{v \in \Vr}\lambda_v}
=  \frac{2^{n-1}} {n! \prod_{v \in \Vr} \lambda_v}.$$
%By Theorem (\ref{ThmProbRankTree}) and Lemma
%(\ref{LemNumbRank}), and noting that $|\Vr|=n-1$ for a rooted binary tree, we have
%\begin{eqnarray}
%\bP[\cT) &=& \sum_{r \in r(\cT)} \bP[\cT,r) = \sum_{r \in r(\cT)}  \frac{2^{n-1}}{n!(n-1)!} \notag \\
%&=& \frac{|\Vr|!}{\displaystyle \prod_{v \in \Vr} \lambda_v} \cdot \frac{2^{n-1}}{n!(n-1)!} \notag \\
%&=& \frac{2^{n-1}} {\displaystyle n! \prod_{v \in \Vr} \lambda_v}. \notag
%\end{eqnarray}
\end{proof}
\begin{example}
Recall again the ranked tree $(\cT,r)$ in Fig. \ref{FigTree2}. In that tree, $X= \{a,b, \ldots ,k \}$ and $n=|X|=11$.
Let $\bP_Y[\cT,r]$ be the probability that the ranked tree $(\cT,r)$ evolved under the Yule model. With Theorem (\ref{ThmProbRankTree}), we get
$$\bP_Y[\cT,r] = \frac{2^{n-1}}{n!(n-1)!} = \frac{2^{10}}{11! \times 10!} \approx 0.71 \times 10^{-11}$$
With Corollary (\ref{CorProbYule}), we get
$$\bP_Y[\cT] = \frac{2^{n-1}} {n! \prod_{v \in \Vr} \lambda_v}
= \frac{2^{10}} {11! \times 1^5 \times 2 \times 3 \times 4 \times 5
\times 10} \approx 0.21 \times 10^{-7} $$ With Theorem
(\ref{ThmYuleRankGivenT}), we get
$$\bP_Y[r | \cT] = \frac{\prod_{v \in \Vr}\lambda_v} {(n-1)!} = \frac {1^5 \times 2 \times 3 \times 4 \times 5 \times 10} {10!} \approx 0.33 \times 10^{-3}$$
Let $\bP_U[\cT]$ be the probability that $\cT$ evolved under the uniform model. Then,
$$\bP_U[\cT] = 1/(2n-3)!! \approx 0.15 \times 10^{-8}$$
Since $\frac{\bP_Y[\cT]}{\bP_U[\cT]} \approx \frac{0.21}{0.15}\times 10^1 = 14 > 1$, i.e. $\bP_Y[\cT]>\bP_U[\cT]$, the tree $\cT$ (without a ranking) is more likely to have evolved under the Yule model.
\end{example}

\begin{rem}
In Chapter \ref{ChaptRank}, we want to calculate for a given phylogenetic tree $\cT$ the probability
$\bP[r(v)=i, r \in r(\cT) | \cT ]$ for a $v \in \Vr$ under the Yule model where $r(\cT)$ as defined in (\ref{DefRank}).
By Theorem (\ref{ThmYuleRankGivenT}), the rankings for $\cT$ all have the same probability, and therefore
$$\bP[r(v)=i, r \in r(\cT) | \cT ] = \frac{|\{r \in r(\cT) : r(v)=i \}| }{|r(\cT)|}.$$
For the value $|r(\cT)|$, a formula is stated in Lemma \ref{LemNumbRank}.
The value $|\{r \in r(\cT) : r(v)=i \}|$ will be calculated with the algorithm {\sc RankCount}.
%For that given tree $\cT$ we need to check if we can assume that $\cT$ evolved under the Yule model.
\end{rem}
\begin{rem} \label{RemCoales}
Another stochastic model on trees is the coalescent model. \index{coalescent model} The coalescent model starts with $n$ species and goes back in time. At each event,
two species are selected uniformly at random and the two species are joint together, the joint being a new species, the ancestor.
So after $n-1$ joining events, we are left with one species, the root of the tree.

With $i$ remaining species, we have ${i \choose 2}$ possibilities to choose two species for the joint.
The probability for a specific ranked tree is therefore
$$\bP[\cT,r] = \frac {1} {{n \choose 2}{n-1 \choose 2} \ldots {2 \choose 2}} = \frac{2^{n-1}}{n!(n-1)!}$$
which is equivalent to the Yule model.

Thus, the Yule model and the coalescent model are equivalent as long as edge lengths are not considered.
\end{rem}

\subsection{Did the primate tree evolve under Yule?} \label{ExPrimatesYule}
Consider the primate tree $\cT_p$ in Appendix \ref{Primates}. $\cT_p$ has $n=218$ leaves.
We want to calculate the value $\frac{\bP_Y[\cT_p]}{\bP_U[\cT_p]}$ in order to decide whether to favor the Yule model over the uniform model.
Note that $\bP_U[\cT]=\frac{2^{n-1}} {n!c_{n-1}}$ and $\bP_Y[\cT]=\frac{2^{n-1}} { n! \prod_{v \in \Vr} \lambda_v}$.

\begin{figure}
\begin{center}
\input{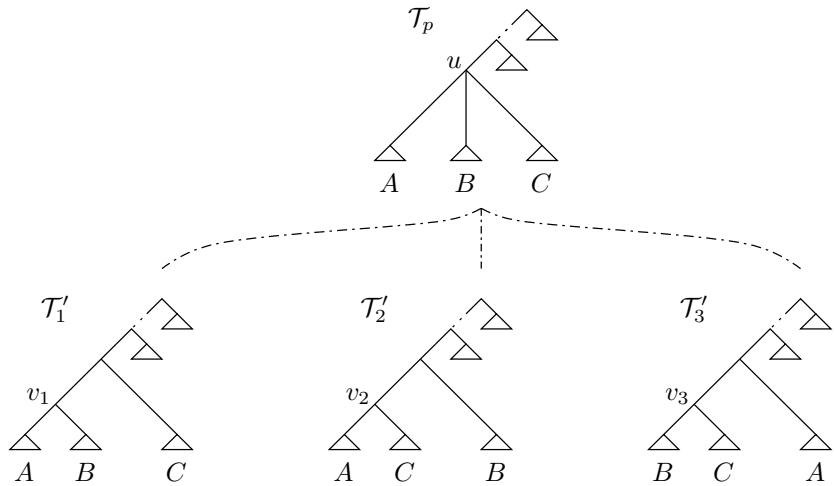}
\caption
{Vertex in $\cT_p$ with three direct descendants. There are three possible binary resolutions.% $\cT_{p_1}$, $\cT_{p_2}$ and $\cT_{p_3}$.
}
\label{FigTreeTestPrim1}
\end{center}
\end{figure}

In $\cT_p$, there are six vertices (vertex labels $48,63,148, 153,157$ and $200$) with more than two direct descendants because the exact resolution is unclear. Five of those vertices have three direct descendants.

%The interior vertex with label 148 has four leaves as direct descendants. There are two different shapes $t_1$ and $t_2$ for a binary tree with four leaves, see Fig. \ref{FigTreeTestPrim2}. Two leaves in the binary resolution form a cherry. Each labeling of the cherry is equally likely because of the exchangeability property.
%Therefore, we can assume that vertex 148 has three direct descendants, two leaves and one cherry with some labeling.

For each vertex with three direct descendants, there are three possible binary resolutions, see Fig. \ref{FigTreeTestPrim1}. \index{binary resolution}

%In general, let $\cT$ be a phylogenetic tree where each vertex has two or three direct descendants.
Let $u$ be a vertex of $\cT_p$ with three direct descandants.
Let $v$ be the additional vertex for a binary resolution of vertex $u$.
For the three different binary resolutions of vertex $u$, we also write $v_1, v_2, v_3$ instead of $v$, see Fig. \ref{FigTreeTestPrim1}.

Let $\cT'$ be a binary resolution of $\cT_p$. Let $\cT'_i$, $i=1,2,3$, be a binary resolution of $\cT_p$ where vertex $u$ is resolved as displayed in Fig. \ref{FigTreeTestPrim1}. Let $\lambda_{v(\cT')}$ be the number of descendants of $v$ in resolution $\cT'$. We want to estimate $\tilde{\lambda_{v}}$.
\begin{eqnarray}
\tilde{\lambda_{v}} &=& \frac{\displaystyle \sum_{\cT'} \lambda_{v(\cT')} \bP[\cT'] }{\displaystyle \sum_{\cT'}  \bP[\cT']} \notag \\
&=& \frac{\displaystyle \sum_{i=1}^3 \sum_{\cT'_i} \lambda_{v_{i}} \bP[\cT'_i] }{\displaystyle \sum_{i=1}^3 \sum_{\cT'_i} \bP[\cT'_i]} \notag \\
&=& \frac{\displaystyle \sum_{i=1}^3 \sum_{\cT'_i} \lambda_{v_{i}} \frac{2^n}{n!\displaystyle \prod_{w \in \Vr_{\cT'_i}} \lambda_w} }{\displaystyle \sum_{i=1}^3 \sum_{\cT'_i} \frac{2^n}{n!\displaystyle \prod_{w \in \Vr_{\cT'_i}} \lambda_w}} \notag \\
&=& \frac{\displaystyle \sum_{i=1}^3 \sum_{\cT'_i}  \frac{2^n}{n!\displaystyle \prod_{w \in \{ \Vr_{\cT'_i} \setminus v_i \} } \lambda_w} }{\displaystyle \sum_{i=1}^3 \frac{1}{\lambda_{v_i}}  \sum_{\cT'_i} \frac{2^n}{n!\displaystyle \prod_{w \in \{ \Vr_{\cT'_i} \setminus v_i \}} \lambda_w}  } \notag
\end{eqnarray}
Note that the inner sum is constant for all $i$, so we get
\begin{eqnarray}
\tilde{\lambda_{v}} &=& \frac{\displaystyle \sum_{\cT'_1} \frac{2^n}{n!\displaystyle \prod_{w \in \{ \Vr_{\cT'_1} \setminus v_1 \} } \lambda_w}  \sum_{i=1}^3  1 }{\displaystyle\sum_{\cT'_1} \frac{2^n}{n!\displaystyle \prod_{w \in \{ \Vr_{\cT'_1} \setminus v_1 \}} \lambda_w}  \sum_{i=1}^3  \frac{1}{\lambda_{v_i}} } \notag \\
&=& \frac{3} {\displaystyle \sum_{i=1}^3 \frac{1}{\lambda_{v_{i}}}} \notag
\end{eqnarray}
With this formula, we estimate the values $\tilde{\lambda_{v}}$ for the new vertex $v$ in the binary resolution of vertex $48,63,153,157$ and $200$.

\begin{figure}
\begin{center}
\input{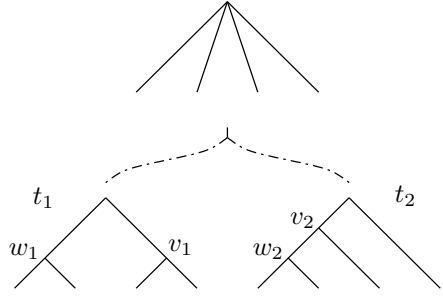}
\caption
{Vertex in $\cT_p$ with four leaf-descendants.% The probability of a labeled tree $\cT_1$ with shape $t_1$ under the Yule model is $\bP[\cT_1]=1/9$. There are three possible labellings of $t_1$, so $\bP[t_1]=3 \cdot 1/9 = 1/3$. Further, $\bP[t_2]=1-\bP[t_1]=2/3$.
}
\label{FigTreeTestPrim2}
\end{center}
\end{figure}

The interior vertex with label 148 has four leaves as direct descendants.
There are two different shapes $t_1$ and $t_2$ for a binary tree with four leaves, see Fig. \ref{FigTreeTestPrim2}.
In $t_1$, the new interior vertices $v_1$ and $w_1$ have the value $\lambda_{v_1} = 1$ and $\lambda_{w_1} = 1$. In $t_2$, the new vertex ${v_2}$ has $\lambda_{v_2} = 2$, the new vertex $w_2$ has $\lambda_{w_2} = 1$.
We set $\tilde{\lambda_{w}}=1$ in $\cT_p$ since $\lambda_{w_1} = 1$ and $\lambda_{w_2} = 1$.
We want to estimate $\tilde{\lambda_{v}}$, the value $\tilde{\lambda_{v}}$ shall be the weighted sum of the $\lambda_{v_i}$,
$$\tilde{\lambda_{v}} = \frac{\bP_Y[t_1]\lambda_{v_1}+\bP_Y[t_2]\lambda_{v_2}}{\bP_Y[t_1]+\bP_Y[t_2]}=1/3 \cdot 1 + 2/3 \cdot 2 = 5/3.$$
With those estimated values for $\tilde{\lambda_{v}}$, we now estimate $\frac{\bP_Y[\cT]}{\bP_U[\cT]}$.
Let $\cT_i, i=1,\ldots,m$, be the binary resolutions of $\cT$. We get
$$\frac{\bP_Y[\cT]}{\bP_U[\cT]} = \frac{\sum_i \bP_Y[\cT_i]}{\sum_i \bP_U[\cT_i]} \approx \frac {c_{n-1}} {\prod_{v \in \Vr_{\cT}} \lambda_v \cdot \prod \tilde{\lambda_{v}}} \approx 0.25 \times 10^{14}$$
which favors the Yule model over the uniform model.
Note that without the estimates for $\tilde{\lambda_v}$, we would have to calculate $\bP_Y[\cT_i]$ and $\bP_U[\cT_i]$ for the
$3^5 \times 15$ linear resolutions of $\cT$.

In Section \ref{EstEdgeLength}, we will assume that the primate tree $\cT_p$ evolved under the Yule model.

\section{Yule model vs. uniform model}
As we have seen in Corollary (\ref{CorProbUnif}), the probability of generating a given tree $\cT$ with $n$ leaves under the uniform model is
\label{ProbYU}
$$\bP_U[\cT] = \frac{2^{n-1}} {n!c_{n-1}}.$$
By Corollary (\ref{CorProbYule}), the probability of generating a given tree $\cT$ under the Yule model is
$$\bP_Y[\cT] = \frac{2^{n-1}} { n! \prod_{v \in \Vr} \lambda_v}.$$
The fraction of the two probabilities, the `Bayes factor' \cite{Everitt1998}, is
$$\frac {\bP_Y[\cT]} {\bP_U[\cT]} = \frac {c_{n-1}} {\prod_{v \in \Vr} \lambda_v}.$$
Given a tree $\cT$, we want to know if it evolved under the Yule or the uniform model.
The fraction $\frac {\bP_Y[\cT]} {\bP_U[\cT]}$ being bigger than $1$ suggests favoring the Yule model, the fraction being smaller than $1$
suggests favoring the uniform model.
So $\ln \left( \frac {\bP_Y[\cT]} {\bP_U[\cT]} \right) $ being bigger than $0$ suggests favoring the Yule model, the logarithm
being smaller than $0$ suggests favoring the uniform model.
% (see also `likelihood ratio statistic' in \cite{Everitt1998}).
In the following, we want to calculate the expected value $\bE_Y [\ln \left( \frac {\bP_Y[\cT]} {\bP_U[\cT]} \right) ]$, given the tree $\cT$ evolved
under the Yule model. We will see that $\bE_Y \left[ \ln \left( \frac {\bP_Y[\cT]} {\bP_U[\cT]} \right) \right]$
is the `Kullbach-Liebler' distance (defined below) between
$\bP_Y$ and $\bP_U$, and show that it goes to infinity with increasing $n$.
Further, $\bE_U \left[ \ln \left( \frac {\bP_U[\cT]} {\bP_Y[\cT]} \right) \right]$
goes to infinity with increasing $n$.
Therefore, for $n$ large enough, the value $\ln \left( \frac {\bP_Y[\cT]} {\bP_U[\cT]} \right)$
is relevant to the question of testing whether a tree evolved under the Yule or uniform model.
In Section \ref{StatTest}, we will actually test the Yule model against the uniform model.
\subsection{The Kullbach-Liebler distance} \index{Kullbach-Liebler distance}
\begin{defi} \index{entropy} \label{DefEntropy}
Let $X$ be a discrete random variable which takes values in the finite set $\Omega=\{w_1,w_2, \ldots ,w_n \}$ with associated probabilities
$\{ p(\omega_1), p(\omega_2) , \ldots , p(\omega_n) \}$. We call this probability distribution $p$.
The \emph{information content} of an event $\omega \in \Omega$ is
$$I(\omega) = - \ln p(\omega)$$ \index{information content}
The $entropy$ $\bJ_p$ of the probability distribution $p$ is defined as
$$\bJ_p = \bE[I(X)] = -\sum_{\omega \in \Omega} p(\omega) \ln p(\omega)$$
\end{defi}
In \cite{McKenzie2000}, Chapter 7, the entropy $\bJ_Y$ for the Yule distribution over $RB(n)$
and the entropy $\bJ_U$ for the uniform distribution over $RB(n)$
are calculated.
Recall that for two functions $f(n)$ and $g(n)$, we write $f(n) \sim g(n)$ precisely if $\lim_{n \rightarrow \infty} \frac {f(n)} {g(n)} = 1$.\\

For $\bJ_Y$, one has (from \cite{McKenzie2000}) \index{entropy!Yule} \index{Yule model!entropy}
\begin{equation} \label{EqnEntrYule}
\bJ_Y=n \sum_{k=2}^{n-1} \frac{g(k)}{k+1}
\end{equation}
where $g(k)=\frac{1-k}{k} \ln\frac{k-1}{2}+\ln\frac{k}{2} + \ln(k+1) - \frac{1}{k}\ln k!$. Asymptotically, one has
\begin{equation} \label{EqnEntrYuleAs}
\bJ_Y - n \ln(n) + c_1n \sim -\frac{1}{2} \ln(n)
\end{equation}
where $c_1 = \ln(2)\ln(\frac{200}{49e})+\ln(9)\ln (\frac{7}{10}) + 2 {\rm Li}_2(\frac{7}{4}) - 2 {\rm Li}_2(\frac{5}{2})-1 \approx 0.493$
and ${\rm Li}_2(x)=\int_1^x \frac {\ln t}{1-t} dt$.\\

For $\bJ_U$, one has (again from \cite{McKenzie2000}) \index{entropy!uniform} \index{uniform model!entropy}
\begin{equation} \label{EqnEntrUnif}
\bJ_U = \ln|RB(n)| = \ln(2n-3)!!
\end{equation}
and asymptotically
\begin{equation} \label{EqnEntrUnifAs}
\bJ_U - n \ln (n)  + c_2n \sim - \ln(n)
\end{equation}
where $c_2=1-\ln(2) \approx 0.307$.
\begin{defi} \label{DefKL}
Let $p$ and $q$ be probability distributions over a finite set
$\Omega$. The {\it Kullbach-Liebler distance} between $p$ and $q$ is
defined as
$$d_{KL}(p,q) = \sum_{\omega \in \Omega} p(\omega) \ln \frac {p(\omega)} {q(\omega)}.$$
\end{defi}
\begin{rem}
The Kullbach-Liebler distance is positive definite, i.e. $d_{KL}(p,q) \geq 0$ with $d_{KL}(p,q) = 0$ iff $p=q$.
Notice that $d_{KL}(p,q) = \infty$ iff there exists a $u \in \Omega$ with $p(u) > 0,$ $q(u) = 0$.
For $p=\bP_Y$ and $q=\bP_U$, both $d_{KL}(p,q)$ and $d_{KL}(q,p)$ are finite,
since $\bP_Y[\cT] > 0$ and $\bP_U[\cT] > 0$ for all $\cT \in RB(n)$.
Note that the Kullbach-Liebler distance between $p$ and $q$ is not symmetric, i.e. we have $d_{KL}(p,q) \neq d_{KL}(q,p)$ in general.
\end{rem}
\begin{rem} \label{RemKLExpV}
Note that the Kullbach-Liebler distance between the probability distributions $p$ and $q$ over the set $\Omega$ equals the following expected value
$$d_{KL}(p,q) = \sum_{\omega \in \Omega} p(\omega) \ln \frac {p(\omega)} {q(\omega)}=\bE_p [\ln \frac {p} {q}  ].$$
\end{rem}
\begin{lem} \label{LemKullLieb}
Let $\Omega$ be a finite set.
Let $p$ be any probability distribution over $\Omega$, and let $q$ be the uniform distribution over $\Omega$.
Then
$$d_{KL} (p,q) = \bJ_q - \bJ_p.$$
\end{lem}
\begin{proof}
By assumption, $q(\omega) = 1 / { | \Omega|}$ for all $\omega \in \Omega$.
From the definition of $d_{KL}(p,q)$, it follows that
\begin{eqnarray}
d_{KL}(p,q) &=& \sum_{\omega \in \Omega} p(\omega) \ln \frac {p(\omega)} {q(\omega)} \notag \\
&=& \sum_{\omega \in \Omega} p(\omega) \ln p(\omega) - \sum_{\omega \in \Omega} p(\omega) \ln q(\omega) \notag \\
&=& -\bJ_p - \sum_{\omega \in \Omega} p(\omega) \ln \frac{1}{|\Omega|} \notag \\
&=& -\bJ_p - \left( \ln \frac{1}{|\Omega|} \right) \sum_{\omega \in \Omega} p(\omega) \notag \\
&=& -\bJ_p - \ln \frac{1}{|\Omega|} \notag \\
&=& -\bJ_p - \sum_{\omega \in \Omega} \frac{1}{|\Omega|}  \ln \frac{1}{|\Omega|} \notag \\
&=& \bJ_q - \bJ_p .\notag
\end{eqnarray}
\end{proof}
\subsection{Kullbach-Liebler distance between $\bP_Y$ and $\bP_U$}
In the following, we calculate the Kullbach-Liebler distance between the Yule distribution
$\bP_Y$ and the uniform distribution $\bP_U$ over $RB(n)$.
\begin{thm} \label{ThmKLYU}  \index{Kullbach-Liebler distance!Yule-uniform}
Let $\bP_Y$ be the Yule distribution and $\bP_U$ be the uniform distribution over $RB(n)$.
The Kullbach-Liebler-distance between those two distributions is
$$d_{KL} (\bP_Y,\bP_U) = \ln(2n-3)!! - n \sum_{k=2}^{n-1} \frac{g(k)}{k+1}$$
where $g(k)$ is again defined as $g(k)=\frac{1-k}{k} \ln\frac{k-1}{2}+\ln\frac{k}{2} + \ln(k+1) - \frac{1}{k}\ln k!$.
Asymptotically, we have
$$d_{KL} (\bP_Y,\bP_U) -c_Yn \sim  - 1/2 \ln(n)$$
with $c_Y \approx 0.186$.
\end{thm}
\begin{proof}
From Lemma (\ref{LemKullLieb}), we have $d_{KL} (\bP_Y,\bP_U) = \bJ_U - \bJ_Y$. With Equations (\ref{EqnEntrYule}) and (\ref{EqnEntrUnif}), we get
$d_{KL} (\bP_Y,\bP_U) = \ln(2n-3)!! - n \sum_{k=2}^{n-1} \frac{g(k)}{k+1}$.
For the asymptotic behavior, we get with Equation (\ref{EqnEntrYuleAs}) and (\ref{EqnEntrUnifAs})
\begin{eqnarray}
\bJ_U   - n \ln (n)  + c_2n - (\bJ_Y - n \ln(n) + c_1n) &\sim& - \ln(n) + 1/2 \ln(n) \notag \\
\bJ_U - \bJ_Y  - (c_1 - c_2)n &\sim& - 1/2 \ln(n) \notag \\
\bJ_U - \bJ_Y  - c_Yn &\sim& - 1/2 \ln(n) \notag
\end{eqnarray}
where $c_Y= c_1-c_2 \approx 0.186$.
\end{proof}
\begin{cor}
For the expected value $\bE_Y [\ln \frac {\bP_Y} {\bP_U}]$, we get
$$\bE_Y [\ln \frac {\bP_Y} {\bP_U}] -c_Yn \sim - 1/2 \ln(n) $$
So $\bE_Y [\ln \frac {\bP_Y} {\bP_U}] \rightarrow \infty$ for $n \rightarrow \infty$.
\end{cor}
\begin{proof}
With Theorem (\ref{ThmKLYU}), we get
\begin{eqnarray}
\bE_Y [\ln \frac {\bP_Y} {\bP_U}] -c_Yn  &=& \sum_{\cT \in RB(n)} \bP_Y[\cT] \ln \frac {\bP_Y[\cT]} {\bP_U[\cT]} -c_Yn  \notag \\
&=& d_{KL} (\bP_Y,\bP_U) -c_Yn  \notag \\
&\sim&- 1/2 \ln(n) \notag
\end{eqnarray}
That implies $d_{KL} (\bP_Y,\bP_U) \sim c_Yn$ and since $c_Y > 0$, we have $\bE_Y [\ln \frac {\bP_Y} {\bP_U}] \rightarrow \infty$ for $n \rightarrow \infty$.
\end{proof}
\subsection{Kullbach-Liebler distance between $\bP_U$ and $\bP_Y$}
\begin{figure}
\begin{center}
\input{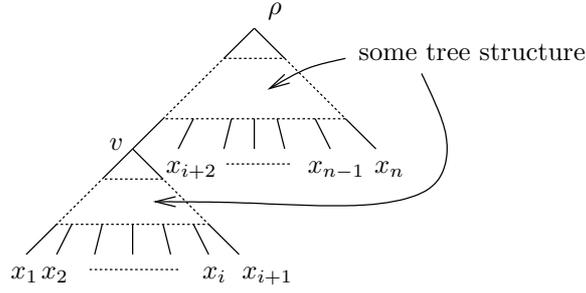}
\caption[Counting the pairs $(\cT,v)$ in Lemma (\ref{LemDoubleSum})]
{Counting the pairs $(\cT,v)$ in Lemma (\ref{LemDoubleSum}). The variables $(x_1, \ldots , x_{i+1})$ take any distinct values from $X'$,
the variables $(x_{i+2}, \ldots , x_{n-1}, x_n)$ take any distinct values from $X''$.}
\label{FigTreest03}
\end{center}
\end{figure}

In the following, we calculate the Kullbach-Liebler distance between the uniform distribution
$\bP_U$ and the Yule distribution $\bP_Y$ over $RB(n)$.
\begin{lem} \label{LemCentBin}
The central binomial coefficient ${2m \choose m}$ can be written as
$${2m \choose m} = 2^{2m} \prod_{j=1}^m \frac{2j-1}{2j}.$$
\end{lem}
\begin{proof}
\begin{eqnarray}
{2m \choose m} &=& \frac {(2m)!}{m!m!} = \frac {2^{2m} \cdot 2m \cdot (2m-1) \cdot (2m-2) \ldots 3 \cdot 2 \cdot 1}
 {2m \cdot 2m \cdot 2(m-1) \cdot 2(m-1) \ldots 4 \cdot 4 \cdot 2 \cdot 2}\notag \\
&=& 2^{2m} \prod_{j=0}^{m-1} \frac {2m-2j-1}{2(m-j)}\notag \\
&=& 2^{2m} \prod_{j=1}^{m} \frac {2j-1}{2j}.\notag
\end{eqnarray}
\end{proof}
\begin{lem} \label{LemDoubleSum}
For the set $RB(n)$, we have
$$\sum_{\cT \in RB(n)} \sum_{v \in \Vr_{\cT}} \ln \lambda_v = \sum_{i=1}^{n-1} \ln i {n \choose i+1} |RB(i+1)| |RB(n-i)|$$
where $\lambda_v$ is defined as in Lemma (\ref{LemNumbRank}).
\end{lem}
\begin{proof}
We have $\lambda_v \in \{1, 2, \ldots ,(n-1) \}$ since a binary tree $\cT$ with $n$ leaves has $n-1$ interior vertices. We rewrite the double sum as
$$\sum_{\cT \in RB(n)} \sum_{v \in \Vr_{\cT}} \ln \lambda_v = \sum_{i=1}^{n-1} \ln i \cdot | \{(\cT,v) : \cT \in RB(n), v \in \Vr_{\cT}, \lambda_v = i \} |$$
To calculate $| \{(\cT,v) : \cT \in RB(n), v \in \Vr_{\cT}, \lambda_v = i \} |$,
we have to count all the pairs $(\cT,v)$ with $v \in \Vr_{\cT}$ having exactly $i$ interior nodes as descendants.
For a binary tree, this is equivalent to $v$ having $i+1$ leaves as descendants ({\it cf.} Figure \ref{FigTreest03}).
So for an interior vertex $v$, we choose a subset $X'$ of $X$ consisting of $i+1$ elements, which shall label the leaf descendants of $v$.
We have ${n \choose i+1}$ possibilities to choose those $i+1$ elements.
There are $|RB(i+1)|$ possibilities to build up a binary tree with leaf set $X'$ and root $v$.
Let $X'' = (X \setminus X') \cup v$, so $|X''|=n-i$. For the set $X''$, there are $|RB(n-i)|$ possible binary trees.
Combining all those possibilities yields
$$| \{\cT,v : \cT \in RB(n), v \in \Vr_{\cT}, \lambda_v = i \} | = {n \choose i+1} |RB(i+1)| |RB(n-i)|$$
which proves the Lemma.
\end{proof}
\begin{thm} \label{ThmKLUY}
For the distance $d_{KL}(\bP_U,\bP_Y)$, it holds that
%$$d_{KL}(\bP_U,\bP_Y) = n \sum_{i=2}^{n-1} \left[ \frac {\ln i}{i+1} \prod_{j=1}^{n-i-1} \frac {(j+i)(2j-i)}{(2(j+i)-1)j} \right] -\ln c_{n-1}$$
$$d_{KL}(\bP_U,\bP_Y) = n S_n -\ln c_{n-1}$$
where $S_n = \sum_{i=2}^{n-1} \left[ \frac {\ln i}{i+1} \prod_{j=1}^{n-i-1} \frac {1-\frac{1}{2j}}{1-\frac{1}{2(j+i)}} \right]$ and
$c_n$ are the Catalan numbers as defined in Lemma (\ref{Lem-n!!}). \index{Kullbach-Liebler distance!uniform-Yule}
\end{thm}
\begin{proof}
By definition of the Kullbach-Liebler distance and with Corollary (\ref{CorProbUnif}) and (\ref{CorProbYule}) and setting $N=|RB(n)|$, we have,
\begin{eqnarray}
d_{KL}(\bP_U,\bP_Y) &=& \sum_{\cT \in RB(n)} \bP_U[\cT] \ln \frac {\bP_U[\cT]} {\bP_Y[\cT]} \notag \\
&=& \sum_{\cT \in RB(n)} \frac{2^{n-1}} {n!c_{n-1}} \ln \left[ \frac { \frac{2^{n-1}} {n!c_{n-1}}}
{\frac{2^{n-1}} { n! \prod_{v \in \Vr_{\cT}} \lambda_v}} \right] \notag \\
&=& \sum_{\cT \in RB(n)} \frac{1} {N} \ln \left[ \frac {\prod_{v \in \Vr_{\cT}} \lambda_v} {c_{n-1}} \right] \notag \\
&=& \frac {1}{N} \left[ \sum_{\cT \in RB(n)} \sum_{v \in \Vr_{\cT}} \ln \lambda_v \right] - \ln c_{n-1} \notag \\
&=& \frac {1}{N} s  - \ln c_{n-1} \label{EqnKLdist1}
\end{eqnarray}
where $s=\sum_{\cT \in RB(n)} \sum_{v \in \Vr_{\cT}} \ln \lambda_v$.
With Lemma (\ref{LemDoubleSum}) and Lemma(\ref{Lem-n!!}), we get
\begin{eqnarray}
s &=& \sum_{\cT \in RB(n)} \sum_{v \in \Vr_{\cT}} \ln \lambda_v  \notag \\
&=& \sum_{i=2}^{n-1} \ln i {n \choose i+1} |RB(i+1)| |RB(n-i)| \notag \\
&=& \sum_{i=2}^{n-1} \ln i {n \choose i+1} \frac{ c_{i} (i+1)!} {2^{i}}  \cdot \frac{c_{n-i-1} (n-i)!}{2^{n-i-1}} \notag \\
&=& \frac{n!}{2^{n-1}} \sum_{i=2}^{n-1} \ln i \frac {(i+1)!(n-i)!} {(i+1)! (n-i-1)!}  c_{i}  c_{n-i-1} \notag \\
&=& \frac{N}{c_{n-1}} \sum_{i=2}^{n-1} \ln i \cdot (n-i) \cdot c_{i} c_{n-i-1} \notag \\
&=& \frac{N n}{{2(n-1) \choose n-1}} \sum_{i=2}^{n-1} \frac {\ln i} {i+1}  {2i \choose i} {2(n-i-1) \choose n-i-1} \notag
\end{eqnarray}
With Lemma (\ref{LemCentBin}) we get
\begin{eqnarray}
s &=& \frac{N n}{2^{2(n-1)} \prod_{j=1}^{n-1} \frac{2j-1}{2j}} \sum_{i=2}^{n-1} \left[ \frac {\ln i} {i+1} 2^{2i} \prod_{j=1}^i \frac{2j-1}{2j}  2^{2(n-i-1)} \prod_{j=1}^{n-i-1} \frac{2j-1}{2j} \right] \notag \\
&=& N n   \sum_{i=2}^{n-1} \left[ \frac {\ln i} {i+1} \prod_{j=1}^{n-1} \frac{2j}{2j-1} \prod_{j=1}^i \frac{2j-1}{2j}  \prod_{j=1}^{n-i-1} \frac{2j-1}{2j} \right] \notag \\
&=& N n   \sum_{i=2}^{n-1} \left[ \frac {\ln i} {i+1} \prod_{j=i+1}^{n-1} \frac{2j}{2j-1}  \prod_{j=1}^{n-i-1} \frac{2j-1}{2j} \right] \notag \\
&=& N n   \sum_{i=2}^{n-1} \left[ \frac {\ln i} {i+1} \prod_{j=1}^{n-i-1} \frac{2(j+i)}{2(j+i)-1}  \prod_{j=1}^{n-i-1} \frac{2j-1}{2j} \right] \notag \\
&=& N n   \sum_{i=2}^{n-1} \left[ \frac {\ln i} {i+1} \prod_{j=1}^{n-i-1} \frac {(j+i)(2j-1)}{(2(j+i)-1)j} \right] \notag \\
&=& N n   \sum_{i=2}^{n-1} \left[ \frac {\ln i} {i+1} \prod_{j=1}^{n-i-1} \frac {2j-1}{2j-\frac {2j} {2(j+i)}} \right] \notag \\
&=& N n   \sum_{i=2}^{n-1} \left[ \frac {\ln i}{i+1} \prod_{j=1}^{n-i-1} \frac {1-\frac{1}{2j}}{1-\frac{1}{2(j+i)}} \right] \notag
%$$d_{KL}(\bP_U,pY) = \frac {1}{|RB(n)|} \sum_{\cT \in RB(n)} \sum_{v in \Vr} \ln \lambda_v - \ln c_{n-1}$$
\end{eqnarray}
Combining this result with Equation (\ref{EqnKLdist1}) establishes the theorem.
\end{proof}
\begin{lem} \label{LemCataAsy}
The asymptotic behavior of the $n$-th Catalan number $c_n$ is
$$c_n \sim n \ln 4$$
\end{lem}
\begin{proof}
With the Stirling formula, $\ln n! \sim n \ln n - n$ (see \cite{Bronstein2000}),we get
\begin{eqnarray}
\ln c_n &=& -\ln(n+1) + \ln{2n \choose n} \notag \\
&=& -\ln(n+1) + \ln(2n)! - 2\ln n! \notag \\
&\sim& -\ln(n+1) + 2n\ln 2n - 2n - 2 n \ln n + 2n \notag \\
&=& -\ln(n+1) + 2n\ln 2 \notag \\
&\sim& n\ln 4 \notag
\end{eqnarray}
\end{proof}
\begin{thm} \label{ThmKLUYAsy}
The Kullbach-Liebler distance between $\bP_U$ and $\bP_Y$ is asymptotically
$$d_{KL}(\bP_U,\bP_Y) \sim c_U n$$
where $c_U$ is a positive constant.
\end{thm}
\begin{proof}
From Theorem (\ref{ThmKLUY}), we have
$$d_{KL}(\bP_U,\bP_Y) = n \cS_n - \ln c_{n-1}$$
with $\cS_n = \sum_{i=2}^{n-1} \left[ \frac {\ln i}{i+1} \prod_{j=1}^{n-i-1} \frac {1-\frac{1}{2j}}{1-\frac{1}{2(j+i)}} \right]$ and $c_n$ being the $n$-th Catalan number.
By Lemma (\ref{LemCataAsy}), it holds $c_{n-1} \sim n \ln 4$.
In Section \ref{ApSn}, we show that $$\ln 4 < 1.44 < S_n < S'+N$$ for all $n \geq 200$ with $S'$ and $N$ being some fixed constants. This yields to
$$d_{KL}(\bP_U,\bP_Y) = n S_n -\ln c_{n-1} \sim n S_n - n \ln 4 \sim c_U n$$
with $c_U$ being a positive constant.
\end{proof}

\begin{cor}
We obtain
$$\bE_U [\ln \frac {\bP_U}{\bP_Y}] \rightarrow \infty \qquad {\rm for \ } n \rightarrow \infty $$
since $\bE_U [\ln \frac {\bP_U}{\bP_Y}] = d_{KL}(\bP_U,\bP_Y)$ by Remark (\ref{RemKLExpV}).
\end{cor}

\subsection{Calculating $S_n$} \label{ApSn}
In Theorem (\ref{ThmKLUY}), we obtain the following formula for the Kullbach-Liebler distance between $\bP_U$ and $\bP_Y$:
$$d_{KL}(\bP_U,\bP_Y) = n S_n -\ln c_{n-1}$$
with $S_n = \sum_{i=2}^{n-1} \left[ \frac {\ln i}{i+1} \cdot a_{n,i} \right]$
and $a_{n,i}=\prod_{j=1}^{n-i-1} \frac {1-\frac{1}{2j}}{1-\frac{1}{2(j+i)}}$.
In the following, we will calculate an upper and a lower bound for $S_n$.
Note that $\{a_{n,i} , n \in \bN\}$ is monotone decreasing for fixed $i$ and $a_{n,i}>0$. So $\lim_{n \rightarrow \infty} a_{n,i}$ exists.
%$$S := \lim_{n \rightarrow \infty} S_n$$
$$a_i := \lim_{n \rightarrow \infty} a_{n,i} = \prod_{j=1}^{\infty} \frac {1-\frac{1}{2j}}{1-\frac{1}{2(j+i)}} = \prod_{j=1}^i \left( 1-\frac{1}{2j} \right)>0$$
$$S_n' := \sum_{i=2}^{n-1} \left[ \frac {\ln i}{i+1} \cdot a_{i} \right]$$
With the property
$$\ln(1-x) = -x - \sum_{i=2}^{\infty} \frac {x^i}{i} \leq -x$$
for $0 \leq x<1$ (see \cite{Zwillinger1996}) and the property
$$\sum_{j=1}^i \frac{1}{j} \geq \int_{1}^{i} \frac {1}{x} dx = \ln(i)$$
we get the following:
\begin{eqnarray}
\ln a_i &=& \sum_{j=1}^i \ln ( 1-\frac{1}{2j} ) \notag \\
&\leq& -\frac{1}{2} \sum_{j=1}^i \frac{1}{j} \notag \\
&\leq& -\frac{1}{2} \ln (i) \notag
\end{eqnarray}
So we have
$$a_i \leq \frac{1}{\sqrt{i}}$$
In the following, we show that $S_n'$ converges.
\begin{eqnarray}
S_n' &=& \sum_{i=2}^{n-1} \left[ \frac {\ln i}{i+1} \cdot a_{i} \right] \notag \\
&\leq& \sum_{i=2}^{n-1} \frac{\ln i}{i^{3/2}} \notag
\end{eqnarray}
Since $\sum_{i=2}^{\infty} \frac{\ln i}{i^{3/2}}$ converges, it follows that $\{S_n', n \in \bN \}$ is bounded.
The sequence $\{ S_n', n \in \bN \}$ is monotone increasing since $\frac {\ln i}{i+1} \cdot a_{i} >0$ for all $i \in \bN, i \geq 2$.
So $\lim_{n \rightarrow \infty} S_n'$ exists and we define $$\lim_{n \rightarrow \infty} S_n' := S'.$$
Now we calculate an upper and a lower bound for $S_n$.
Since $a_{i,n} \rightarrow a_i$, there exists an $N \in \bN$ s.t. $a_{i,n} < (1 + 1/S') a_i $ for all $n>N$.
$$S_n = \sum_{i=2}^{n-1} \left[ \frac {\ln i}{i+1} \cdot a_{n,i} \right] < (N-1) + \sum_{i=N+1}^{n-1} \left[ \frac {\ln i}{i+1} \cdot (1 + 1/S') a_i \right] < (N-1)+(1+1/S')S_{n}' $$
Since $S_n'$ is monotone increasing, we get
$$S_n < (N-1)+(1+1/S')S_{n}' < (N-1)+(1+1/S')S'$$
which yields to
$$S_n < S' + N.$$
Since $a_{i,n} > a_i$, we have
$$S_n = \sum_{i=2}^{n-1} \left[ \frac {\ln i}{i+1} \cdot a_{n,i} \right]  > \sum_{i=2}^{n-1} \left[ \frac {\ln i}{i+1} \cdot a_{i} \right] = S_n'$$
So we get $S_n > S_n'$ for all $n$. With Maple, I calculated $S_{200}' \approx 1.44 > \ln 4$.
Overall, we have $$\ln 4 < 1.44 < S_n < S'+N$$
for all $n \geq 200$.

\chapter{Trees and Martingales} \label{Martingale}
In this chapter, we have a closer look at the process of the tree generation. We will see that
the tree generation is a certain stochastic process, a martingale.
Under the uniform model, the martingale fulfills the conditions for the Azuma inequality.

We make use of this property at the end of the chapter. We test the Yule model against the uniform model with the log-likelihood-ratio test.
With the Azuma inequality, we find an analytical bound for the power of the test.
Since the algorithms in Chapter \ref{ChaptRank} work in particular for trees under the Yule model, it will be useful to have a test for deciding whether a tree evolved under Yule.

First, we provide some basic definitions and properties on conditional probability and martingales.
\section{Conditional probability and martingales}
%This section introduces martingales and proves some properties.
%The results in this sections are from \cite{Ross1996} which will be needed in the following sections.
\begin{defi} \index{conditional expectation}
Let $X$ (resp. $Y$) be a discrete random variable which takes
values $\{x_i, i \in \bN\}$ (resp. $\{ y_i, i \in \bN\}$). The {
\it conditional expectation} $$Z=\bE[X|Y]=\sum_{j} x_j
\bP[X=x_j|Y]$$ is a random variable. Z takes values
$$z_i = \sum_{j} x_j \bP[X=x_j|Y=y_i]$$ on the set $\{ Y=y_i \}$ with probability $\bP[Z=z_i]=\bP[Y=y_i]$.
\end{defi}
The two equations in the next Lemma are stated in \cite{Ross1996} with a brief verification. We will give a full proof.
\begin{lem}  \label{LemCondExp}
Let $X$ (resp. $Y$, $U$) be a discrete random variable which takes
values $\{x_i, i \in \bN\}$ (resp. $\{ y_i, i \in \bN\}$, $\{ u_i,
i \in \bN\}$). Further, assume $\bE[|X|] < \infty $. Then, we get the following two equalities:
\begin{eqnarray}
\bE[X] &=& \bE[\bE[X|Y]] \label{EqnCondExp1} \\
\bE[X|U] &=& \bE[\bE[X|Y,U]|U] \label{EqnCondExp2}
\end{eqnarray}
\end{lem}
\begin{proof}
Let $Z=\bE[X|Y]$. We obtain Equation (\ref{EqnCondExp1}) from
\begin{eqnarray}
\bE[\bE[X|Y]] &=& \sum_{i} z_i \bP[Z=z_i] \notag \\
&=& \sum_i \sum_j x_j \bP[X=x_j|Y=y_i] \bP[Y=y_i] \notag \\
&=& \sum_i \sum_j x_j \bP[X=x_j,Y=y_i] \qquad \qquad (\ast)\notag \\
&=& \sum_j \sum_i x_j \bP[X=x_j,Y=y_i] \notag \\
&=& \sum_j x_j \bP[X=x_j] \notag \\
&=& \bE[X] \notag
\end{eqnarray}
The summation order in ($\ast$) can be changed since $\bE[|X|]<\infty$.

It is left to verify (\ref{EqnCondExp2}). Let $W=\bE[X|Y,U]$. The random variable $W$ takes a value
$$w_{j_1,j_2} = \sum_k x_k \bP[X=x_k | Y=y_{j_1}, U=u_{j_2}]$$
with probability $\bP[Y=y_{j_1}, U=u_{j_2}]$ where $j_1 \in \bN$
and $j_2 \in \bN$. Let $Z=\bE[W|U]$. The random variable $Z$ takes
a value
$$z_i = \bE[W|U=u_i]$$
with probability $\bP[U=u_i]$ where $i \in \bN$. We
transform $z_i$ to
\begin{eqnarray}
z_i &=& \bE[W|U=u_i] \notag \\
%&=& \bE[\bE[X|Y, U=u_i]|U=u_i] \notag \\
&=& \sum_{j_1,j_2} w_{j_1,j_2} \bP[W=w_{j_1,j_2} | U=u_i] \notag \\
&=& \sum_{j_1,j_2} \sum_k x_k \bP[X=x_k | Y=y_{j_1}, U=u_{j_2}] \bP[Y=y_{j_1}, U=u_{j_2} | U=u_i] \notag \\
&=& \sum_{j_1} \sum_k x_k \bP[X=x_k | Y=y_{j_1}, U=u_i] \bP[Y=y_{j_1}| U=u_i] \notag \\
&=& \sum_{j_1} \sum_k x_k \bP[X=x_k, Y=y_{j_1}, U=u_i] / \bP[U=u_i] \qquad \qquad (\ast \ast) \notag \\
&=& \sum_k \sum_{j_1} x_k \bP[X=x_k, Y=y_{j_1}, U=u_i] / \bP[U=u_i]  \notag \\
&=& \sum_k x_k \bP[X=x_k, U=u_i] / \bP[U=u_i] \notag \\
&=& \sum_k x_k \bP[X=x_k|U=u_i] \notag \\
&=& \bE[X|U=u_i] \notag
\end{eqnarray}
The summation order in $(\ast \ast)$ can be changed since $\bE[|X|]<\infty$.
So we obtain $$\bE[\bE[X|Y,U]|U=u_i] =
\bE[X|U=u_i]$$ for all $i \in \bN$, i.e.
$\bE[\bE[X|Y,U]|U] = \bE[X|U]$.
%$$\bE[X] = \bE[\bE[X|Y,U]]$$
%by Equation (\ref{EqnCondExp1}). Calculating the expected values conditional on the value of $U$, we get (\ref{EqnCondExp2}).
\end{proof}
\begin{defi} \index{martingale}
A stochastic process $\{Z_n, n \in \bN \}$ is called a {\it
martingale} if
$$\bE[|Z_n|] < \infty \qquad \forall n \in \bN$$
and
\begin{equation}
\bE[Z_{n+1}|Z_1, Z_2, \ldots ,Z_n]= Z_n. \label{EqnMart}
\end{equation}
\end{defi}
\begin{rem}
Taking expectations of (\ref{EqnMart}) with Equation (\ref{EqnCondExp1}) gives
$$\bE[Z_{n+1}] = \bE[Z_n].$$
\end{rem}
The results of Lemma (\ref{LemMartY}) and Theorem (\ref{ThmMartingale}) are already stated in \cite{Ross1996}.
Again, the following proofs are more detailed.
\begin{lem} \label{LemMartY}
Let $\{Z_n, n \in \bN \}$ be a discrete stochastic process with $\bE[|Z_n|]<\infty$. Let ${\textbf Y}$ be a
vector of discrete random variables. If
$$\bE[Z_{n+1}|Z_1, \ldots, Z_n, {\textbf Y}] = Z_n$$
then $\{Z_n\}$ is a martingale.
\end{lem}
\begin{proof}
It holds $\bE[Z_n|Z_1, \ldots ,Z_n] = Z_n$ since $\bE[Z_n|Z_1=z_1, \ldots ,Z_n=z_n]= z_n$.
With that property and with Equation (\ref{EqnCondExp2}), we get
\begin{eqnarray}
\bE[Z_{n+1}|Z_1, \ldots ,Z_n] &=& \bE[\bE[Z_{n+1}|Z_1, \ldots ,Z_n, {\textbf Y}]|Z_1, \ldots ,Z_n] \notag \\
&=& \bE[Z_n|Z_1, \ldots ,Z_n] \notag \\
&=& Z_n. \notag
\end{eqnarray}
\end{proof}
\begin{thm} \label{ThmMartingale}
Let $X, Y_1, Y_2, \ldots$ be discrete random variables such that $\bE[|X|] < \infty$ and let
$$Z_n = \bE[X|Y_1, \ldots Y_n]$$
for all $n \in \bN$. Then $\{Z_n , n \in \bN \}$ is a martingale.
\end{thm}
\begin{proof}
%Since $X$ takes values from a finite set, it holds $\bE[|X|]< \infty$.
With Equation (\ref{EqnCondExp1}), we get
$\bE[|Z_n|] = \bE[|\bE[X|Y_1, \ldots ,Y_n]|] \leq  \bE[\bE[|X| |Y_1, \ldots ,Y_n]] = \bE[|X|] < \infty$.
To check the second condition for a martingale, it is, by
Lemma (\ref{LemMartY}), sufficient to show that $\bE[Z_{n+1}|Z_1,
\ldots Z_n, Y_1, \ldots, Y_n]=Z_n$. We have
\begin{eqnarray}
\bE[Z_{n+1}|Z_1, \ldots Z_n, Y_1, \ldots, Y_n] &=& \bE[Z_{n+1}|Y_1, \ldots, Y_n] \notag \\
&=& \bE[\bE[X|Y_1, \ldots, Y_{n+1}]|Y_1, \ldots, Y_n] \notag \\
&=& \bE[X|Y_1, \ldots, Y_n] \qquad ({\rm from \ } (\ref{EqnCondExp2})) \notag \\
&=& Z_n \notag
\end{eqnarray}
which proves the theorem.
\end{proof}
\subsection{The Azuma inequality}
Let $\{Z_i, i \in \bN \}$ be a martingale. If the random varialbes $Z_i$ do not change too fast over time, Azuma's inequality gives us some
bounds on their probabilities.

The following theorem, the Azuma inequality, is stated in \cite{Ross1996} with a detailed proof.
\begin{thm}[Azuma's Inequality] \label{ThmAzuma} \index{Azuma's inequality}
Let $\{Z_i, i \in \bN \}$ be a martingale with $\bE[Z_i] = \mu$. Let $Z_0=\mu$ and suppose that for nonnegative constants $\alpha_j$, $\beta_j$, $j \geq 1$,
$$-\alpha_j \leq Z_{j} - Z_{j-1} \leq \beta_j.$$
Then for any $i \geq 0$, $a>0$:
\begin{eqnarray}
(i) & & \bP[Z_i -  \mu \geq a] \leq \exp \{- \frac{2a^2}{\sum_{j=1}^i (\alpha_j+\beta_j)^2} \} \notag \\
(ii) & & \bP[Z_i -  \mu \leq -a] \leq \exp \{- \frac{2a^2}{\sum_{j=1}^i (\alpha_j+\beta_j)^2} \} \notag
\end{eqnarray}
\end{thm}
The following corollary will be very useful for the next section.
\begin{cor} \label{CorAzuma}
Let $\{Z_i, i \in \bN \}$ be a martingale with $\bE[Z_i] = \mu$. Let $Z_0=\mu$ and suppose that for a nonnegative constant $\sC$, $j \geq 1$,
$$|Z_j - Z_{j-1}| \leq \sC$$
Then for any $i \in \bN$:
$$ \bP[Z_i \leq 0] \leq \exp \{- \frac{\mu^2}{2 i \sC^2} \} $$
\end{cor}
\begin{proof}
Let $\alpha_i = \beta_i = \sC$ for all $i \in \bN$ and $a=\mu$. Then inequality $(ii)$ in Theorem (\ref{ThmAzuma}) establishes the corollary.
\end{proof}

\section{A martingale process on trees under the uniform model} \label{SecMartUni} \index{martingale on trees}
In this section, we assume that a tree $\cT \in RB(n)$ evolved under the uniform model. Consider the following setting:
\begin{itemize}
\item Let $h_U: RB(n) \rightarrow \bR$ with $h_U(\cT) = \ln
\frac{\bP_U[\cT]}{\bP_Y[\cT]} = \ln \frac{\prod_{v \in
\Vr_{\cT}} \lambda_v}{c_{n-1}}$.
\item For $j \in \{1, \ldots
,n\}$, let $Y_j: RB(n) \rightarrow RB(j)$ with $Y_j(\cT) =
\cT|_{\{1 \ldots j\}}$.
\item For $j > n$, let $Y_j: RB(n) \rightarrow RB(n)$
with $Y_j(\cT) = \cT$.
\item Let $Z_i = \bE[h_U|Y_1, \ldots Y_i]$.
\end{itemize}
We have $\bE[|h_U(\cT)|] < \infty$ since $\cT$ is chosen from the finite set $RB(n)$ and $\max_{\cT \in RB(n)} |h_U(\cT)| <\infty$.
With Theorem (\ref{ThmMartingale}), we obtain that $\{Z_i, i \in \bN \}$ is
a martingale.
Note that $$Z_i = \bE[h_U|Y_1, \ldots Y_i] = \bE[h_U|Y_i].$$ For all $i \geq n$, we have
$$Z_i = \bE[h_U(\cT)|Y_i=\cT] = h_U(\cT).$$
The expected value $\mu_U$ of $Z_n$ is, with Remark (\ref{RemKLExpV}),
$$\mu_U = \bE[Z_n]  = \bE[h_U(\cT)]=d_{KL}(\bP_U,\bP_Y).$$
Theorem (\ref{ThmKLUYAsy}) shows
$$d_{KL}(\bP_U,\bP_Y) \sim c_U n$$
which means
$$\mu_U \sim c_U n.$$

In the following, we want to apply Azuma's inequality to the tree martingale $\{Z_i, i \in \bN \}$.
First, set $Z_0 := \bE[Z_n] = d_{KL}(\bP_U,\bP_Y)$.
To apply Azuma's inequality, we have to verify $|Z_i - Z_{i-1}| \leq \sC_U$ for all $i \in \bN$.
\begin{itemize}
\item For $i=1$, note that by definiton, we have $$Z_1 = \bE[h_U(\cT)|Y_1] =  \bE[h_U(\cT)] = d_{KL}(\bP_U,\bP_Y) = Z_0$$
so $|Z_1 - Z_0| = 0$.
\item For $i \geq n$, note that $Z_i = \bE[h_U(\cT)|\cT]=  h_U(\cT)$. So $|Z_i - Z_{i-1}| =0$ for all $i>n$.
\item Section (\ref{SecAzumaConst}) will establish $|Z_i - Z_{i-1}| \leq \ln n$ for $2 \leq i \leq n$.
\end{itemize}
With Corollary (\ref{CorAzuma}), we then have
\begin{eqnarray}
\bP[Z_n \leq 0] &\leq& \exp \{- \frac{\mu_U^2}{2 n (\ln n)^2} \} \notag \\
&\sim& \exp \{- \frac{c_U^2 n}{2 (\ln n)^2} \} \rightarrow 0 \qquad {\rm for} \ n \rightarrow \infty \notag
\end{eqnarray}
Note that $Z_n = h_U(\cT)= \ln \frac{\bP_U[\cT]}{\bP_Y[\cT]}$.
So for a tree $\cT$ generated under the uniform model, the probability that $\bP_U[\cT]$ is smaller than $\bP_Y[\cT]$ tends to $0$ quickly with $n$
as the number of leaves tends to $\infty$.
Therefore the Bayes factor $\frac {\bP_U[\cT]}{\bP_Y[\cT]}$ is a very good indicator as to whether a `big' tree evolved under the uniform model or not.

\subsection{Calculating a bound in the Azuma inequality} \label{SecAzumaConst}
Let $\{Z_i, i \in \bN \}$ be the tree martingale introduced above.
We can transform $Z_i$ into
\begin{eqnarray}
Z_i &=& \bE[h_U|Y_i] \notag \\
&=& \sum_{\cT \in RB(n)} h_U(\cT) \bP[\cT | Y_i] \notag \\
&=& \sum_{\cT \in RB(n)} \ln \frac{\prod_{v \in \Vr_{\cT}} \lambda_v}{c_{n-1}} \bP[\cT | Y_i] \notag \\
&=& \sum_{\cT \in RB(n)} \left( \sum_{v \in \Vr_{\cT}} \ln \lambda_v - \ln c_{n-1} \right) \bP[\cT | Y_i] \notag \\
&=& \left[ \sum_{\cT \in RB(n)} \left( \sum_{v \in \Vr_{\cT}} \ln
\lambda_v  \right) \bP[\cT | Y_i] \right] - \ln c_{n-1}\notag
\end{eqnarray}
The random variable $Z_i$ therefore takes values
$$ z_{i,t} = \left[ \sum_{\cT \in RB(n)} \left( \sum_{v \in \Vr_{\cT}} \ln \lambda_v  \right) \bP[\cT | Y_i = t] \right] - \ln c_{n-1} $$
for all $t \in RB(i)$.

Assuming that $\cT$ was generated under the uniform
model, i.e.
$$\bP[\cT | Y_i = t] = \frac{\bP[\cT]}{\bP[t]} = \frac{|RB(i)|}{|RB(n)|}$$
we get, for $t \in RB(i)$,
\begin{eqnarray}
z_{i,t} &=& \left[ \sum_{\substack{\cT \in RB(n) \\ \cT|_{\{1, \ldots, i\}}=t}} \left( \sum_{v \in \Vr_{\cT}} \ln \lambda_v  \right) \frac{|RB(i)|}{|RB(n)|} \right] - \ln c_{n-1}\notag \\
&=& \frac{|RB(i)|}{|RB(n)|} \left[ \sum_{\substack{\cT \in RB(n) \\
\cT|_{\{1, \ldots, i\}}=t}} \sum_{v \in \Vr_{\cT}} \ln
\lambda_v  \right] - \ln c_{n-1}.\notag
\end{eqnarray}
Let $\cT$ be a binary phylogenetic tree. For the subtree $\cT|_{\{1, \ldots, i\}}$, we will write $\cT(i)$.
The set of all binary phylogenetic trees with leave set $\{1, \ldots , i-1, i+1, \ldots n\}$ shall be $RB(n,i)$.
In the following, we will calculate an upper bound for $|Z_i-Z_{i-1}|$.
Note that $$|Z_i-Z_{i-1}| = \max_{t \in RB(i)} |z_{i,t} - z_{(i-1),t(i-1)}|.$$
The difference $|z_{i,t} - z_{(i-1),t(i-1)}|$ is
\begin{eqnarray}
\Delta_{i,t} &=& |z_{i,t} - z_{(i-1),t(i-1)}| \notag \\
&=& \left| \frac{|RB(i)|}{|RB(n)|}  \sum_{\substack{\cT \in RB(n) \\ \cT(i)=t}} \sum_{v \in \Vr_{\cT}} \ln \lambda_v   - \frac{|RB(i-1)|}{|RB(n)|}  \sum_{\substack{\cT \in RB(n) \\ \cT(i-1)=t(i-1)}} \sum_{v \in \Vr_{\cT}} \ln \lambda_v \right|   \notag \\
&=& \frac{|RB(i-1)|}{|RB(n)|} \left| (2i-3)  \sum_{\substack{\cT \in RB(n) \\ \cT(i)=t}} \sum_{v \in \Vr_{\cT}} \ln \lambda_v  - \sum_{\substack{t' \in RB(i) \\ t'(i-1)=t(i-1) }} \sum_{\substack{ \cT \in RB(n) \\ \cT(i)=t'}} \sum_{v \in \Vr_{\cT}} \ln \lambda_v   \right| \notag \\
&=& \frac{|RB(i-1)|}{|RB(n)|} \left| \sum_{\substack{t' \in RB(i) \\ t'(i-1)=t(i-1) }}  \sum_{\substack{\cT \in RB(n) \\ \cT(i)=t}} \sum_{v \in \Vr_{\cT}} \ln \lambda_v  - \sum_{\substack{t' \in RB(i) \\ t'(i-1)=t(i-1) }} \sum_{\substack{ \cT \in RB(n) \\ \cT(i)=t'}} \sum_{v \in \Vr_{\cT}} \ln \lambda_v   \right| \notag \\
&=& \frac{|RB(i-1)|}{|RB(n)|} \left| \sum_{\substack{t' \in RB(i) \\ t'(i-1)=t(i-1) }}  \left( \sum_{\substack{\cT \in RB(n) \\ \cT(i) =t}} \sum_{v \in \Vr_{\cT}} \ln \lambda_v  - \sum_{\substack{ \cT \in RB(n) \\ \cT(i)=t'}} \sum_{v \in \Vr_{\cT}} \ln \lambda_v \right)  \right| \notag \\
&=& \frac{|RB(i-1)|}{|RB(n)|} \left| \sum_{\substack{t' \in RB(i) \\ t'(i-1)=t(i-1) }}  \sum_{\substack{\cT'  \in RB(n,i) \\ \cT'(i-1)=t(i-1) }}  \left( \sum_{\substack{\cT \in RB(n) \\ \cT \setminus i = \cT' \\ \cT(i) =t}} \sum_{v \in \Vr_{\cT}} \ln \lambda_v  - \sum_{\substack{ \cT \in RB(n) \\ \cT \setminus i = \cT'  \\ \cT(i)=t'}} \sum_{v \in \Vr_{\cT}} \ln \lambda_v \right)  \right| \notag \\
&\leq& \frac{|RB(i-1)|}{|RB(n)|}  \sum_{\substack{t' \in RB(i) \\ t'(i-1)=t(i-1) }}  \sum_{\substack{\cT'  \in RB(n,i) \\ \cT'(i-1)=t(i-1) }}  \left| \sum_{\substack{\cT \in RB(n) \\ \cT \setminus i = \cT' \\ \cT(i) =t}} \sum_{v \in \Vr_{\cT}} \ln \lambda_v  - \sum_{\substack{ \cT \in RB(n) \\ \cT \setminus i = \cT'  \\ \cT(i)=t'}} \sum_{v \in \Vr_{\cT}} \ln \lambda_v   \right| \notag
\end{eqnarray}
Define $$s := \left| \sum_{\substack{\cT \in RB(n) \\ \cT \setminus i = \cT' \\ \cT(i) =t}} \sum_{v \in \Vr_{\cT}} \ln \lambda_v  - \sum_{\substack{ \cT \in RB(n) \\ \cT \setminus i = \cT'  \\ \cT(i)=t'}} \sum_{v \in \Vr_{\cT}} \ln \lambda_v   \right|.$$

\begin{figure}
\begin{center}
\input{treemart01.pstex_t}
\caption{Tree $\cT$ where leaf $i$ is moved}
\label{TreeMart01}
\end{center}
\end{figure}

\noindent
Consider the tree $\cT$ in Fig. \ref{TreeMart01}. Moving leaf $i$ to a new position will change $\lambda_v$ of a vertex $v$, if $v$ is on the path $P$ from $v_i$ to $v_i'$.
The change
of $\lambda_v$, when $v <_{\cT} v_i$, is $\lambda_v^{new}=\lambda_v - 1$. For the other vertices on that path, we have
$\lambda_v^{new}=\lambda_v + 1$.
So we get, with the property $\ln x - \ln y = \ln x/y$,
\begin{eqnarray}
s &=&    \left| \sum_{\substack{\cT \in RB(n) \\ \cT \setminus i = \cT' \\ \cT(i) =t}}
\left( \sum_{\substack{v \in \Vr_{\cT} \setminus v_i \\ v \in P \\ v <_{\cT} v_i  }} \left( \ln \frac {\lambda_v}  {\lambda_v - 1} \right) +
\sum_{\substack{v \in \Vr_{\cT} \setminus v_i \\ v \in P \\ v <_{\cT} v_i' }} \left( \ln \frac {\lambda_v} {\lambda_v + 1} \right) +
\ln \lambda_{v_i} - \ln \lambda_{v_i}' \right)  \right| \notag \\
&\leq& \sum_{\substack{\cT \in RB(n) \\ \cT \setminus i = \cT' \\ \cT(i) =t}}
\left| \sum_{\substack{v \in \Vr_{\cT} \setminus v_i \\ v \in P \\ v <_{\cT} v_i }} \left( \ln \frac{\lambda_v}{\lambda_v - 1} \right) +
\sum_{\substack{v \in \Vr_{\cT} \setminus v_i \\ v \in P \\ v <_{\cT} v_i' }} \left( \ln \frac{\lambda_v}{\lambda_v + 1} \right) +
\ln \lambda_{v_i} - \ln \lambda_{v_i}' \right| \notag \\
&=&     \sum_{\substack{\cT \in RB(n) \\ \cT \setminus i = \cT' \\ \cT(i) =t}}
\left| \sum_{\substack{v \in \Vr_{\cT} \setminus v_i \\ v \in P \\ v <_{\cT} v_i }} \left( \ln \frac {\lambda_v}{\lambda_v - 1} \right) +
\sum_{\substack{v \in \Vr_{\cT} \setminus v_i \\ v \in P \\ v <_{\cT} v_i' }} \left( \ln \frac {\lambda_v }{\lambda_v + 1} \right) +
s' \right| \notag
%&\leq& \sum_{\substack{\cT \in RB(n) \\ \cT \setminus i = \cT' \\ \cT(i) =t}}
%\sum_{k=1}^{n-1} \frac{k+1}{k} \notag
\end{eqnarray}
with
$$s' = \left\{
\begin{array}{ll}
    \sum_{i=\lambda_{v_i}'+1}^{\lambda_{v_i}} \ln \frac {i}{i-1} & \hbox{if $\lambda_{v_i}' \leq \lambda_{v_i}$} \\
    \sum_{i=\lambda_{v_i}+1}^{\lambda_{v_i}'} \ln \frac {i-1}{i} & \hbox{if $\lambda_{v_i} < \lambda_{v_i}'$} \\
\end{array}
\right.$$
Note that for any $v, w \in P$ with $v,w <_{\cT} v_i$ or $v,w <_{\cT} v_i'$, we have $\lambda_v \neq \lambda_w$. That yields to
\begin{eqnarray}
s &\leq&      \sum_{\substack{\cT \in RB(n) \\ \cT \setminus i = \cT' \\ \cT(i) =t}}
 \sum_{k=1}^{n-1} \ln \frac{k+1}{k}  \notag
\end{eqnarray}
Overall, we get, with using the property $\ln(1+x) < x$ for $x>0$,
\begin{eqnarray}
|z_{i,t} - z_{(i-1),t(i-1)}| &\leq& \frac{|RB(i-1)|}{|RB(n)|}  \sum_{\substack{t' \in RB(i) \\ t'(i-1)=t(i-1) }}  \sum_{\substack{\cT'  \in RB(n,i) \\ \cT'(i-1)=t(i-1) }}   \sum_{\substack{\cT \in RB(n) \\ \cT \setminus i = \cT' \\ \cT(i) =t}} \sum_{k=1}^{n-1} \ln \frac {k + 1}{k} \notag\\
&=& \frac{|RB(i)|}{|RB(n)|}  \sum_{\substack{\cT'  \in RB(n,i) \\ \cT'(i-1)=t(i-1) }}  \sum_{\substack{\cT \in RB(n) \\ \cT \setminus i = \cT' \\ \cT(i) =t}} \sum_{k=1}^{n-1} \ln \left( 1+ \frac {1}{k} \right)  \notag\\
&=& \frac{|RB(i)|}{|RB(n)|}  \sum_{\substack{\cT \in RB(n) \\ \cT(i)=t }} \sum_{k=1}^{n-1} \ln \left( 1+ \frac {1}{k} \right)  \notag \\
&=& \sum_{k=1}^{n-1} \ln \left( 1+ \frac {1}{k} \right) \notag \\
&<& \sum_{k=1}^{n-1} \frac {1}{k} \notag \\
&<& \int_1^n \frac {1}{x} dx \notag \\
&=& \ln n. \notag
\end{eqnarray}
Therefore,
$$
|Z_i - Z_{i-1}| = \max_{t \in RB(i)} |z_{i,t} - z_{(i-1),t(i-1)}| \leq \ln n.
$$

\section{A martingale process on trees under the Yule model}
In this section, we assume that a tree $\cT$ evolved under the Yule model. Consider the following setting:
\begin{itemize}
\item Let $h_Y(\cT) = - h_U(\cT) =  \ln \frac{\bP_Y[\cT]}{\bP_U[\cT]}$.
\item For $j \in \{1, \ldots
,n\}$, let $Y_j: RB(n) \rightarrow RB(j)$ with $Y_j(\cT) =
\cT|_{\{1 \ldots j\}}$.
\item For $j > n$, let $Y_j: RB(n) \rightarrow RB(n)$
with $Y_j(\cT) = \cT$.
\item Let $ \tilde{Z}_i = \bE[h_Y|Y_1, \ldots Y_i]$.
\end{itemize}
Since $h_Y = - h_U$, the process $\{ \tilde{Z}_i, i \in \bN \}$ is a martingale with the same argumentation as in Section \ref{SecMartUni}.
Further, from Section \ref{SecMartUni}, we get
$$\tilde{Z}_i= -\left[ \sum_{\cT \in RB(n)} \left( \sum_{v \in \Vr_{\cT}} \ln
\lambda_v  \right) \bP[\cT | Y_i] \right] + \ln c_{n-1}$$
and
$$ \tilde{z}_{i,t} = -\left[ \sum_{\cT \in RB(n)} \left( \sum_{v \in \Vr_{\cT}} \ln \lambda_v  \right) \bP[\cT | Y_i = t] \right]
 + \ln c_{n-1}$$
for all $t \in RB(i)$.

\section{Hypothesis testing: Did $\cT$ evolve under the Yule model?} \label{StatTest}
\index{hypothesis test}

In this section, the hypothesis that a given tree $\cT$ evolved under the Yule model is tested against the uniform model.

In \cite{McKenzie2000-2}, a test between the Yule and the uniform model is developed by counting cherries. It is shown that the number of cherries
in a tree is normally distributed with different expected values for the two models. The power of the test stated in \cite{McKenzie2000-2} is above $0.90$
for trees with more than $80$ leaves. The power is only stated as an asymptotic result though.

We will give an analytic result for the power of the log-likelihood-ratio test for the Yule model against the uniform model.

\bigskip
First, we recall the basics about hypothesis testing. In a hypothesis test, we test for a given dataset $x$ if a hypothesis $H_0$
is rejected in favor of a hypothesis $H_1$ or if $H_0$ is accepted.
The hypothesis test is characterized by a decision rule, it decides if $H_0$ is accepted.

\bigskip
\noindent
The Type I error of a hypothesis test is \index{Type I error}
$$\alpha = \bP[H_0 \ rejected \ |H_0 \ true].$$
The Type II error of a hypothesis test is \index{Type II error}
$$\beta = \bP[H_0 \ retained \ |H_1 \ true].$$

\noindent
The power of the test is $1-\beta$. \index{power of a test}

\bigskip
The next Lemma, the Neyman-Pearson Lemma (see \cite{Ross1996}), states that for a given Type I error, the likelihood-ratio
test is the test with the smallest Type II error.
\begin{lem}[Neyman-Pearson Lemma] \index{Neyman-Pearson Lemma}
When performing a hypothesis test between two point hypotheses $H_0$ and $H_1$,
then the likelihood-ratio test which rejects $H_0$ in favor of $H_1$ when \index{likelihood-ratio test}
$$\frac{ \bP[x|H_0 \ true]}{ \bP[x|H_1 \ true]} \leq k$$
with $k$ being some positive constant, is the most powerful test of size $\alpha$,
where $\alpha  = \bP[\frac{ \bP[x|H_0 \ true]}{ \bP[x|H_1 \ true]} \leq k |H_0 \ true]=  \bP[H_0 \ rejected |H_0 \ true]$ as defined above.
\end{lem}
Note that the log-likelihood-ratio test, \index{log-likelihood-ratio test} i.e. rejecting $H_0$ if $$\ln \frac{ \bP[x|H_0 \ true]}{ \bP[x|H_1 \ true]} \leq \ln k $$
is equivalent to the likelihood-ratio test.
We will test the Yule model against the uniform model with the log-likelihood-ratio test to get the smallest Type II error.

\bigskip
\noindent
Let $H_0$ and $H_1$ be the following hypotheses.
\bigskip
%\noindent

\qquad $H_0$: \qquad $\cT$ evolved under the Yule model

\qquad $H_1$: \qquad $\cT$ evolved under the uniform model

\bigskip
\noindent
The decision rule for this test shall be:
\begin{itemize}
\item $\tilde{Z}_n = \ln \frac{\bP_Y [\cT]}{\bP_U [\cT]} > 0 \Rightarrow $ accept $H_0$.
\item $\tilde{Z}_n = \ln \frac{\bP_Y [\cT]}{\bP_U [\cT]} \leq 0 \Rightarrow$ reject $H_0$.
\end{itemize}

\bigskip
The Type I and Type II error can be obtained with simulations, i.e. construct a lot of trees with $n$ leaves under the Yule model
and estimate $\alpha$ and $\beta$.

%With the results from the previous sections, we can provide analytical bounds for the Type I and Type II error.
With the results from the previous sections, we can provide an analytical bound for the Type II error.

%\bigskip
%\noindent
%A bound for the Type I error of this test is, with Corollary (\ref{CorAzuma}),
%\begin{eqnarray}
%\alpha = \bP[H_0 \ rejected \ |H_0 \ true] &=& \bP_Y [\ln \frac{\bP_Y [\cT]}{\bP_U [\cT]} \leq 0] \notag \\
%&\leq& \exp\{ - \frac {\mu_Y^2} {2n \sC_Y^2}\} \notag
%\end{eqnarray}

\bigskip
\noindent
A bound for the Type II error of this test is, with Corollary (\ref{CorAzuma}) and Theorem (\ref{ThmKLUY}),
\begin{eqnarray}
\beta = \bP[H_0 \ retained \ |H_1 \ true] &=& \bP_U [\ln \frac{\bP_Y [\cT]}{\bP_U [\cT]} > 0] \notag \\
&=& \bP_U [\ln \frac{\bP_U [\cT]}{\bP_Y [\cT]} < 0] \notag \\
&\leq& \exp\{ - \frac {\mu_U^2} {2n \sC_U^2}\} \notag \\
&\leq& \exp\{ - \frac {\mu_U^2} {2n (\ln n)^2}\} \notag \\
&=& \exp\{ - \frac {(n S_n -\ln c_{n-1})^2} {2n (\ln n)^2}\} \label{EqnPower}
\end{eqnarray}
with $S_n$ and $c_n$ as defined in Theorem (\ref{ThmKLUY}). Asymptotically, we get, with Theorem (\ref{ThmKLUYAsy}),
\begin{eqnarray}
\beta &\sim& \exp\{ - \frac {(c_U n)^2} {2n (\ln n)^2}\} \notag \\
&\leq& \exp\{ - \frac {((1.44-\ln 4) n)^2} {2n (\ln n)^2}\} \notag \\
&\approx&  \exp\{ - 0.00144 \frac {n} { (\ln n)^2}\} \notag
\end{eqnarray}

\noindent
So the power of the test, $1-\beta$, tends to $1$ as $n$ tends to $\infty$. 

%The bound for the power of the test is, however, very poor. This is due to two reasons.
%
%When calculating the bound for the Azuma inequality, we did a lot of rough estimations. The bound $\ln n$ can probably be improved.
With the current bound, the power of the test, calculated by Equation (\ref{EqnPower}), is bigger than $0.85$ only
for trees with more than $600$ leaves. It is probably possible to improve the bound for the Azuma inequality though.
If the current bound, $\ln n$, could be improved to $1/4 \ln n$, the power of the test would be bigger than $0.90$ for trees with more than $50$ leaves.
A bound of $1/2 \ln n$ would result in a power bigger than $0.90$ for trees with more than $170$ leaves.

\chapter{The Rank Function} \label{ChaptRank}
Consider the primate tree in Appendix \ref{Primates}.
Was speciation event with label 76 more likely to be an early event in the tree or a late event?
What is the probability that 76 was the 6th speciation event?
Was it more likely that speciation event 76 happened before speciation event 162 or 162 before 76?
This chapter will provide an answer to those questions, under the assumption that each rank function is equally likely, which is, in particular, the case under the Yule model.

The algorithms {\sc RankProb}, {\sc Compare} and an algorithm for obtaining the expected rank and variance for a vertex were implemented in Python. The code is attached in Appendix \ref{PythonCode}. This is joint work with Daniel Ford from Stanford University.

In Section \ref{EstEdgeLength}, we will show how to estimate edge lengths in a tree by calculating the probability distribution of the rank of a vertex. This question was posed by Arne Mooers and Rutger Vos, who constructed the primate supertree and wanted to estimate the edge lengths for it (see \cite{Vos2006}).
\section{Probability distribution of the rank of a vertex}
Let $\cT$ be a binary phylogenetic tree.
Specifying an order for the speciation events (i.e. the interior nodes) in $\cT$ is equivalent to introducing a rank function on $\cT$.
In this chapter, we are interested in the distribution of the possible ranks for a certain vertex, i.e. we want to know the
probability of $r(v)=i$ for a given~$v \in \Vr$. In other words, we want to calculate $\bP[r(v)=i| \cT]$,
with $r \in r(\cT)$, $r(\cT)$ is the set of possible rank functions on the tree $\cT$.
If every rank function on a given tree is equally likely, we have
\begin{equation} \label{EqnUnif}
\bP[r(v)=i| \cT] = \frac{| \{r: r(v)=i, r \in r(\cT) \}|}{| r(\cT)|}
\end{equation}
A formula for the denominator is given in Lemma (\ref{LemNumbRank}).
The enumerator will be calculated in polynomial time by algorithm {\sc RankCount}.\\

\noindent
Examples of stochastic models on phylogenetic trees where each rank function is equally likely:
\begin{itemize}
\item For the Yule model, we have seen in Theorem (\ref{ThmYuleRankGivenT}), that $\bP[r| \cT]$ is the uniform distribution.
\item As we have seen in Remark (\ref{RemCoales}), the coalescent model has the same probability distribution on rooted binary ranked trees as the Yule model.
So $\bP[r| \cT]$ is the uniform distribution.
\item In the uniform model
%(also known as PDA model), 
no rank function is induced when a tree is generated. We can assume though that for a given tree $\cT$,
each rank function is equally likely. Then, Equation (\ref{EqnUnif}) holds as well.
\end{itemize}
\begin{defi} \label{DefiAlpha}
Let $\cT$ be a rooted phylogenetic tree. Define
$$\alpha_{\cT,v}(i) := |\{r: r(v)=i,~r \in r(\cT) \}|$$
${\rm for}\ v\in \Vr,i \in 1, \ldots ,|\Vr|$. In other words,
$\alpha_{\cT,v}(i)$ denotes the number of rank functions $r$ for
$\cT$ in which $v$ comes in the $i$-th position.
\end{defi}
The following results will be needed in the next sections.
\begin{lem} \label{LemSequence}
Let
$$x^1 = \{x_1^1, x_2^1 \ldots x_{n_1}^1\}$$
$$x^2 = \{x_1^2, x_2^2 \ldots x_{n_2}^2\}$$
$$\vdots$$
$$x^d = \{x_1^d, x_2^d \ldots x_{n_d}^d\}$$
be $d$ disjoint sets with the linear order $x_1^i < x_2^i < \ldots < x_{n_i}^i$ for each $i \in \{1, \ldots , d\}$.
The number $\sL$ of possible linear orders on the set $ x^1 \cup x^2 \cup \ldots \cup x^d $,
with the linear order of each original set $x^i$ being preserved, is
$$\sL=\frac {\displaystyle \left( \sum_{i=1}^d n_i \right)!} {\displaystyle \prod_{i=1}^d n_i!} $$
\end{lem}
\begin{proof}
The number $\tilde{\sL}$ of linear orders of the $\sum_{i=1}^d n_i$
elements of $ x^1 \cup x^2 \cup \ldots \cup x^d $, allowing any
order on $x^i$, is $\tilde{\sL} = \left(\sum_{i=1}^d n_i \right)!$. The
number $\tilde{\sL}_i$ of linear orders of the $n_i$ elements of $x^i $ is $(n_i)!$. Since for $\sL$, we only allow the linear order
$x_1^i < x_2^i < \ldots < x_{n_i}^i$ on $x^i$, it holds
$$\sL = \frac{\tilde{\sL}} {\displaystyle \prod_{i=1}^d \tilde{\sL}_i} = \frac {\displaystyle \left( \sum_{i=1}^d n_i \right)!}
{\displaystyle \prod_{i=1}^d n_i!} $$
\end{proof}
\begin{cor} \label{CorSequenceBin}
For $d=2$ in Lemma (\ref{LemSequence}), we have
$$\sL = {n_1+n_2 \choose n_1}$$
possible linear orders on $x^1 \cup x^2$, preserving the linear order on $x^1$ and $x^2$.
\end{cor}
\begin{proof}
From Lemma (\ref{LemSequence}) follows
$$\sL = \frac {\displaystyle \left( \sum_{i=1}^2 n_i \right)!} {\displaystyle \prod_{i=1}^2 n_i!}
= \frac {(n_1+n_2)!} {(n_1)! (n_2)!} = {n_1+n_2 \choose n_1}$$
\end{proof}

\begin{rem} \label{nchoosek}
The values ${n \choose k}$ for all $n,k \leq N$ ($n,k,N \in \bN$) can be calculated in $O(N^2)$, {\it cf.} Pascal's Triangle. In Appendix \ref{PythonCode}, a dynamic programming version for calculating ${n \choose k}$ is implemented.
Thus, after $O(N^2)$ calculations, any value ${n \choose k}$ with $n,k \leq N$ can be obtained
in constant time in an algorithm.
\end{rem}
\subsection{Polynomial-time algorithms}
In the following, we give a polynomial algorithm to determine
$\alpha_{\cT,v}(i)$ for $v \in \Vr$ and $i=1, \ldots ,|\Vr|$
in a binary phylogenetic tree $\cT$.\\

\noindent
%\begin{algorithm} %[h!]
      %\begin{alg}
        {\bf Algorithm}: \textsc{RankCount}($\cT,v$)  \index{algorithm R\textsc{ankCount}}\\
        {\bf Input}: A rooted binary phylogenetic tree $\cT$ and an interior vertex $v$.\\
        {\bf Output}: The values of $\alpha_{\cT,v}(i)$ for ${i=1, \ldots ,|\Vr|}$.
        \begin{algorithmic}[1]
          \STATE Denote the vertices of the path from $v$ to root $\rho$ with \\ $(v=x_1, x_2, \ldots ,x_n=\rho)$.
          \STATE Denote the subtree of $\cT$, consisting of root $x_m$ and all its descendants, by $\cT_m$ for $m=1, \ldots ,n$.
          ({\it cf.} Figure \ref{FigRankCount1}).
          \FOR{$m=1, \ldots ,n$}
            \FOR{$i=1, \ldots , |\Vr_{\cT}|$}
                \STATE $\alpha_{\cT_m,v}(i):=0$
            \ENDFOR
          \ENDFOR
          \STATE $\alpha_{\cT_1,v}(1):=\frac{|\Vr_{\cT_1}|!}{\displaystyle \prod_{v \in \Vr_{\cT_1}}\lambda_v}$
          \FOR{$m = 2, \ldots ,n$}      %LINE 9 mentioned in Compare!!!!
            %\STATE $\cT_{m-1}':=\cT_m \setminus (\cT_{m-1}\cup x_m)$ \qquad ({\it cf.} Figure \ref{FigRankCount2})
            \STATE $\cT_{m-1}':=\cT_m |_{L_{\cT_m} \setminus L_{\cT_{m-1}}}$ \qquad ({\it cf.} Figure \ref{FigRankCount2})
            \STATE $R_{\cT_{m-1}'}:=\frac{|\Vr_{\cT_{m-1}'}|!}{\displaystyle \prod_{v \in \Vr_{\cT_{m-1}'}}\lambda_v}$
            \FOR{$i=m, \ldots ,|\Vr_{\cT_m}|$}
                \STATE $M := \min\{|\Vr_{\cT_{m-1}'}|, i-2\}$
                \STATE $\displaystyle \alpha_{\cT_m,v}(i):=$ \\ \quad $\sum_{j=0}^{M} \alpha_{\cT_{m-1},v} (i-j-1) R_{\cT_{m-1}'}
                {|\Vr_{\cT_{m-1}}|+|\Vr_{\cT_{m-1}'}|-(i-1) \choose |\Vr_{\cT_{m-1}'}|-j} {i-2 \choose j}$ \qquad $(\ast)$
            \ENDFOR
          \ENDFOR
         % \FOR{$i=1, \ldots , (i-2-M_{m-1})$}
%                \STATE $\alpha_{\cT_n,v}(i):=0$
%          \ENDFOR
          \STATE RETURN $\alpha_{\cT,v}:=\alpha_{\cT_n,v}$
        \end{algorithmic}
        \bigskip
      %\end{alg}
%\end{algorithm}
\begin{figure}
\begin{center}
\input{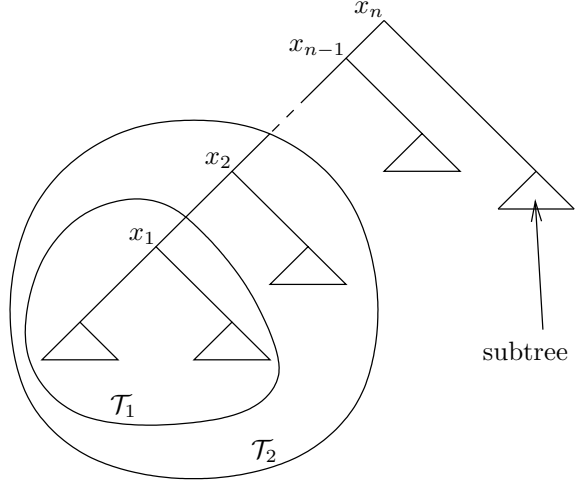}
\caption{Labeling the tree for {\sc RankCount}}
\label{FigRankCount1}
\end{center}
\end{figure}
\begin{figure}
\begin{center}
\input{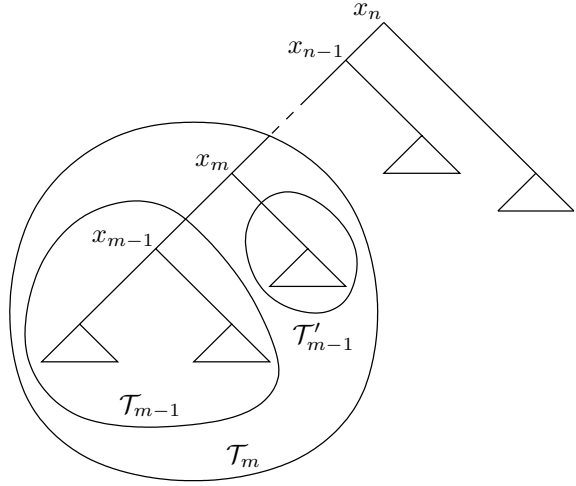}
\caption{Labeling the tree for recursion in {\sc RankCount}}
\label{FigRankCount2}
\end{center}
\end{figure}
\begin{thm} \label{ThmRankCount}
{\sc RankCount} returns the quantities $$\alpha_{\cT,v}(i) = |\{r: r(v)=i, r \in r(\cT) \}|$$
for each given $v\in \Vr$ and all $i \in 1, \ldots ,|\Vr|$.
\end{thm}
\begin{proof}
We have to show that all the $\alpha_{\cT_m,v}(i)$ produced by {\sc RankCount} equal the $\alpha_{\cT_m,v}(i)$ defined in (\ref{DefiAlpha}).
In the following, we denote the values $\alpha_{\cT_m,v}(i)$ produced by the algorithm with $\alpha_{\cT_m,v}^{Alg}(i)$ and $\alpha_{\cT_m,v}(i)$
shall denote the number of rank functions with $r(v)=i$ as defined in (\ref{DefiAlpha}). We will show $\alpha_{\cT_m,v}(i)=\alpha_{\cT_m,v}^{Alg}(i)$
for $m=1, \ldots ,n,~i=1, \ldots ,|\Vr_{\cT}|$. This is done by induction over $m$.\\
For $m=1$, $\alpha_{\cT_1,v}(1)=\alpha_{\cT_1,v}^{Alg}(1)$ since (\ref{LemNumbRank}) holds. Vertex $v$ is the root of $\cT_1$,
so $\alpha_{\cT_1,v}(i)=0$ for all $i>1$.\\
Let $m=k$ and $\alpha_{\cT_m,v}(i)=\alpha_{\cT_m,v}^{Alg}(i)$
holds for all $m<k$. $\alpha_{\cT_k,v}(i)=0$ clearly holds for all
$i > |\Vr_{\cT_k}|$ since $r_{\cT_k}: v \rightarrow \{1, \ldots
,|\Vr_{\cT_{k}}|\}$. So it is left to verify that the
term $(\ast)$ returns the right values
for $\alpha_{\cT_k,v}(i)$.
Assume that the vertex $v$ is in the $(i-j-1)$-th position in $\cT_{k-1}$ (with $i-j-1>0$) for some rank function $r_{\cT_{k-1}}$ and $v$ shall be in the
$i$-th position in $\cT_k$. We want to combine the linear order in the tree $\cT_{k-1}$ induced by $r_{\cT_{k-1}}$ with a linear order in $\cT_{k-1}'$ induced by $r_{\cT_{k-1}'}$
to get a linear order on $\cT_k$. The first $j$ vertices of $\cT_{k-1}'$ must be inserted between vertices of $\cT_{k-1}$ with lower rank than $v$ so that $v$
ends up to be in the $i$-th position of the tree $\cT_k$. We will count the number of possibilities to do so.
The tree $\cT_{k-1}'$ has $$R_{\cT_{k-1}'}=\frac{|\Vr_{\cT_{k-1}'}|!}{\displaystyle \prod_{v \in \Vr_{\cT_{k-1}'}}\lambda_v}$$ possible rank functions.
Combining a rank function $r_{\cT_{k-1}}$ with a rank function $r_{\cT_{k-1}'}$ for getting a rank function $r_{\cT_k}$ with $r_{\cT_k}(v)=i$
means inserting the first $j$ vertices of $\cT_{k-1}'$ anywhere between the first $(i-j-2)$ vertices of $\cT_{k-1}$. There are
$${(i-j-2)+j \choose j} = {i-2 \choose j}$$ possibilities according to Corollary \ref{CorSequenceBin}.
For combining the $|\Vr_{\cT_{k-1}}|-(i-j-1)$ vertices of rank larger than $v$ in $\cT_{k-1}$ with the remaining $|\Vr_{\cT_{k-1}'}|-j$ vertices in $\cT_{k-1}'$,
we have $${|\Vr_{\cT_{k-1}}|-(i-j-1)+|\Vr_{\cT_{k-1}'}|-j \choose |\Vr_{\cT_{k-1}'}|-j} = {|\Vr_{\cT_{k-1}}|+|\Vr_{\cT_{k-1}'}|-(i-1) \choose |\Vr_{\cT_{k-1}'}|-j}$$ possibilities.
This follows again from Corollary \ref{CorSequenceBin}.
The number of rank functions $r_{\cT_{k-1}}$ with $r_{\cT_{k-1}}(v)=i-j-1$ is $\alpha_{\cT_{k-1},v} (i-j-1)$ by the induction assumption.
Multiplying all those possibilities gives
$$\alpha_{\cT_{k-1},v} (i-j-1) R_{\cT_{k-1}'} {|\Vr_{\cT_{k-1}}|+|\Vr_{\cT_{k-1}'}|-(i-1) \choose
|\Vr_{\cT_{k-1}'}|-j} {i-2 \choose j}$$
$\alpha_{\cT_k,v}(i)$ is then the sum over all possible $j$ which is equal to the term $(\ast)$ for $\alpha_{\cT_k,v}^{Alg}(i)$. This establishes the theorem.
\end{proof}
\begin{thm} \label{ThmRankCountTime}
The runtime of {\sc RankCount} is $O(|\Vr|^2)$.
\end{thm}
\begin{proof}
Note that the number of rank functions $R_{\cT}=\frac{|\Vr_{\cT}|!}{\prod_{v \in \Vr_{\cT}}\lambda_v}$ on a tree $\cT$ with $\Vr$ interior vertices can be calculated in $O(|\Vr|)$, i.e. in linear time.

Further, note that the combinatorial factors ${n \choose k}$ for all $n,k \leq |\Vr|$ can be calculated in advance in quadratic time, see Remark (\ref{nchoosek}).
In the algorithm, those factors can then be obtained in constant time.

Contributions to the runtime from each line in {\sc RankCount} (the runtime is always w.r.t. $|\Vr|$):\\
{\it Line 1--2:} linear time\\
{\it Line 3--7:} quadratic time\\
{\it Line 8:} linear time\\
{\it Line 9--16:} quadratic time since:\\
\indent
{\it Line 11:} $R_{\cT_{m-1}'}$ can be calculated in $O(|\Vr|)$. This has to be done for
$m=1, \ldots ,n$, so overall the runtime for calculating all $R_{\cT_{m-1}'}$ is no more than $O(|\Vr|^2)$ since $n \leq |\Vr|$.\\
\indent
{\it Line 14:} We add up all calculations needed for obtaining $\alpha_{\cT_m,v}(i)$, $m=1, \ldots ,n$, $i=1, \ldots ,|\Vr_{\cT_m}|$:
$$\sum_{m=2}^n|\Vr_{\cT_m}||\Vr_{\cT_{m-1}'}| \leq \sum_{m=2}^n|\Vr||\Vr_{\cT_{m-1}'}| = |\Vr| \sum_{m=2}^n|\Vr_{\cT_{m-1}'}| \leq |\Vr|^2$$
The last inequality holds since the vertices of the $\cT_{m}'$, $m=1, \ldots ,n-1$, are distinct.
Therefore, line 14 contributes a quadratic runtime.\\
{\it Line 17:} constant time\\

So overall, the runtime is no more than $O(|\Vr|^2)$. Figure \ref{FigRunTimeBin} shows a tree for which the runtime of
{\sc RankCount} is actually quadratic.
\begin{figure}
\begin{center}
\resizebox{8cm}{!}{
\input{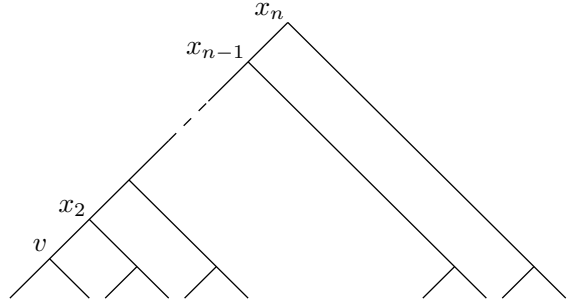}
}
\caption{Illustration for runtime of {\sc RankCount}}
\label{FigRunTimeBin}
\end{center}
\end{figure}
Counting all the calculations for term $(\ast)$ in the algorithm for the tree in \ref{FigRunTimeBin} yields to
\begin{eqnarray}
\sum_{m=2}^{n} \sum_{i=m}^{|\Vr_{\cT_m}|} |\Vr_{\cT_{m-1}'}|+1 &=& \sum_{m=2}^{n} \sum_{i=m}^{|\Vr_{\cT_m}|} 2 \notag \\
&=& \sum_{m=2}^{n} 2(|\Vr_{\cT_m}|-(m-1)) \notag \\
&=& \sum_{m=2}^{n} 2((2m-1)-(m-1)) \notag \\
&=& \sum_{m=2}^{n} 2m \notag \\
&=& n(n+1)-2 \notag
\end{eqnarray}
Since $n=(|\Vr|+1)/2$, we have a quadratic runtime.

%In the $m$-th step of the iteration, we have to calculate $\alpha_{\cT_m,v}(i)$ for $i=1, \ldots , |\Vr_{\cT_m}|=m$. So overall, we have to do
%$$\sum_{m=1}^n m = \frac{n(n+1)}{2}$$ calculations for the sum in $(\ast)$, which gives a runtime of $O(|\Vr|^2)$ since $n=|\Vr|$.
\end{proof}
\begin{cor} \label{CorProbRank}
The probability $\bP[r(v)=i|\cT]$ can be calculated in $O(|V|^2)$. We have
\begin{equation} \label{EqnProbRank} \bP[r(v)=i|\cT]=\frac{\alpha_{\cT,v}(i)}{ \sum_{i=1}^{|\Vr|} \alpha_{\cT,v}(i)}=
\frac{\alpha_{\cT,v}(i) \prod_{v \in \Vr}\lambda_v} {|\Vr|!}. \end{equation}
\end{cor}
\begin{proof}
The first equality in (\ref{EqnProbRank}) follows from basic probability theory. The second equality holds since
$\frac{|\Vr|!}{\prod_{v \in \Vr}\lambda_v}= \sum_i \alpha_{\cT,v}(i)$ by (\ref{LemNumbRank}).
The complexity of the runtime follows from (\ref{ThmRankCountTime}).
\end{proof}
\begin{rem} \label{RemExpVar}
We will write $\bP[r(v)=i]$ instead of $\bP[r(v)=i|\cT]$ in the following. With $\bP[r(v)=i]$ from Corollary (\ref{CorProbRank}), the
expected value $\mu_{r(v)}$ and the variance $\sigma_{r(v)}^2$ for
$r(v)$ can be calculated by
\begin{eqnarray}
\mu_{r(v)} &=& \sum_{i=1}^{|\Vr|} i\bP[r(v)=i] \notag \\
\sigma_{r(v)}^2 &=& \sum_{i=1}^{|\Vr|} i^2 \bP[r(v)=i] - \mu_{r(v)}^2 \notag
\end{eqnarray}
\end{rem}
\begin{example}
We will illustrate the algorithm \textsc{RankCount} for the tree in Figure~\ref{FigExpRankCount}. We get the following values:\\

\begin{figure}
\begin{center}
\input{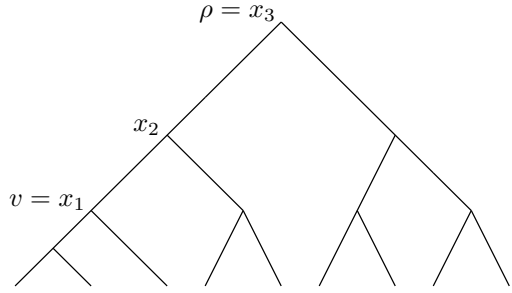}
\caption{Tree to illustrate the algorithm {\sc RankCount}} \label{FigExpRankCount}
\end{center}
\end{figure}

\noindent
$\alpha_{\cT_1,v}(1)=\frac{2!}{2 \cdot 1}=1$\\
\\
$\alpha_{\cT_2,v}(2)=\alpha_{\cT_{m-1},v}(1) 1 {2 + 1 - 1 \choose 1} {0 \choose 0} = 2$\\
$\alpha_{\cT_2,v}(3)=\alpha_{\cT_{m-1},v}(1) 1 {2 + 1 - 2 \choose 1} {1 \choose 0} = 1$\\
$\alpha_{\cT_2,v}(4)=0$\\
\\
$\alpha_{\cT_3,v}(3)=\alpha_{\cT_{m-1},v}(2) 2 {4 + 3 - 2 \choose 3} {1 \choose 0} + \alpha_{\cT_{m-1},v}(1) 2 {4 + 3 - 2 \choose 2} {1 \choose 1} = 40 + 0 = 40$\\
$\alpha_{\cT_3,v}(4)=\alpha_{\cT_{m-1},v}(3) 2 {4 + 3 - 3 \choose 3} {2 \choose 0} + \alpha_{\cT_{m-1},v}(2) 2 {4 + 3 - 3 \choose 2} {2 \choose 1} = 8 + 48 = 56$\\
$\alpha_{\cT_3,v}(5)=\alpha_{\cT_{m-1},v}(3) 2 {4 + 3 - 4 \choose 2} {3 \choose 1} + \alpha_{\cT_{m-1},v}(2) 2 {4 + 3 - 4 \choose 1} {3 \choose 2} = 18 + 36 = 54$\\
$\alpha_{\cT_3,v}(6)=\alpha_{\cT_{m-1},v}(3) 2 {4 + 3 - 5 \choose 1} {4 \choose 2} + \alpha_{\cT_{m-1},v}(2) 2 {4 + 3 - 5 \choose 0} {4 \choose 3} = 24 + 16 = 40$\\
$\alpha_{\cT_3,v}(7)=\alpha_{\cT_{m-1},v}(3) 2 {4 + 3 - 6 \choose 0} {5 \choose 3} = 20$\\
$\alpha_{\cT_3,v}(8)=0$\\
\\
With $\alpha_{\cT_3,v}=\alpha_{\cT,v}$, we get
\begin{eqnarray}
\bP[r(v)=1] &=& 0 \notag \\
\bP[r(v)=2] &=& 0 \notag \\
\bP[r(v)=3] &=& \frac{40}{40 + 56 + 54 + 40 + 20} = \frac{40}{210} = \frac{20}{105} \notag \\
\bP[r(v)=4] &=& \frac{28}{105} \notag \\
\bP[r(v)=5] &=& \frac{27}{105} \notag \\
\bP[r(v)=6] &=& \frac{20}{105} \notag \\
\bP[r(v)=7] &=& \frac{10}{105} \notag \\
\bP[r(v)=8] &=& 0 \notag
\end{eqnarray}
Therefore, the expected value $\mu_{r(v)}$ is
$$\mu_{r(v)} = \sum_{i=1}^8 i\bP[r(v)=i) = \frac{497}{105} \approx 4.73$$
and the variance $\sigma_{r(v)}^2$ is
$$\sigma_{r(v)}^2 = \sum_{i=1}^8 i^2 \bP[r(v)=i) - \mu_{r(v)}^2 = \frac{2513}{105} - \frac{497^2}{105^2} = \frac{344}{225} \approx 1.53$$
\end{example}
%\bigskip

\begin{rem}
Note that $\bP[r(v)=i] = \frac{\alpha_{\cT,v}(i)}{\sum_{j} \alpha_{\cT,v}(j)}$. Common factors in all $\alpha_{\cT,v}(i), i=1, \ldots ,|\Vr_{\cT_v}|$ will therefore cancel out.
\end{rem}

The next algorithm, {\sc RankProb}, is a modification of {\sc RankCount} such that common factors of $\alpha_{\cT,v}(i), i=1, \ldots ,|\Vr_{\cT_v}|$, will not be included. Therefore, the numbers we have to deal with in the algorithm stay smaller and the number of calculations is reduced.\\

\noindent
%\begin{algorithm} %[h!]
      %\begin{alg}
        {\bf Algorithm}: \textsc{RankProb}($\cT,v$)  \index{algorithm R\textsc{ankProb}}\\
        {\bf Input}: A rooted binary phylogenetic tree $\cT$ and an interior vertex $v$.\\
        {\bf Output}: The probabilities $\bP[r(v)=i]$ for ${i=1, \ldots ,|\Vr|}$.
        \begin{algorithmic}[1]
          \STATE Denote the vertices of the path from $v$ to root $\rho$ with \\ $(v=x_1, x_2, \ldots ,x_n=\rho)$.
          \STATE Denote the subtree of $\cT$, consisting of root $x_m$ and all its descendants, by $\cT_m$ for $m=1, \ldots ,n$.
          ({\it cf.} Figure \ref{FigRankCount1}).
          \FOR{$m=1, \ldots ,n$}
            \FOR{$i=1, \ldots , |\Vr_{\cT}|$}
                \STATE $\tilde{\alpha}_{\cT_m,v}(i):=0$
            \ENDFOR
          \ENDFOR
          \STATE $\tilde{\alpha}_{\cT_1,v}(1):=1$
          \FOR{$m = 2, \ldots ,n$}      
            \STATE $\cT_{m-1}':=\cT_m |_{L_{\cT_m} \setminus L_{\cT_{m-1}}}$ \qquad ({\it cf.} Figure \ref{FigRankCount2})
            \FOR{$i=m, \ldots ,|\Vr_{\cT_m}|$}
                \STATE $M := \min\{|\Vr_{\cT_{m-1}'}|, i-2\}$
                \STATE $\displaystyle \tilde{\alpha}_{\cT_m,v}(i):=\sum_{j=0}^{M} \tilde{\alpha}_{\cT_{m-1},v} (i-j-1)
                {|\Vr_{\cT_{m-1}}|+|\Vr_{\cT_{m-1}'}|-(i-1) \choose |\Vr_{\cT_{m-1}'}|-j} {i-2 \choose j}$
            \ENDFOR
          \ENDFOR
          \FOR{$i=1, \ldots , |\Vr_{\cT}|$}
            \STATE $\bP[r(v)=i] := \frac{\tilde{\alpha}_{\cT_n,v}(i)}{\sum_{j} \tilde{\alpha}_{\cT_n,v}(j)}$
          \ENDFOR
          \STATE RETURN $\bP[r(v)=i], i=1,\ldots,|\Vr|.$
        \end{algorithmic}
        \bigskip
      %\end{alg}
%\end{algorithm}

\begin{thm} \label{ThmRankProb}
{\sc RankProb} returns the quantities $$\bP[r(v)=i]$$
for each given $v\in \Vr$ and all $i \in 1, \ldots ,|\Vr|$. The runtime is $O(|\Vr|^2)$.
\end{thm}
\begin{proof}
Note that the structure of {\sc RankProb} is the same as the structure of {\sc RankCount}.
The only difference is that common factors to $\alpha_{\cT_m,v}(i)$ for all $i$ are not included.
Those common factors do not change the probabilities since they cancel out once calculating the probabilities.
Therefore, since {\sc RankCount} works correct, also {\sc RankProb} works correct.

It is left to verify the runtime. The only time consuming step in {\sc RankProb} is line 13. This line is of the same complexity as line 14 in {\sc RankCount}. Line 14 in {\sc RankCount} contributed a quadratic time. Therefore, the runtime of {\sc RankProb} is quadratic as well.
\end{proof}

\subsection{Non-binary trees and ranks}
Let $\cT$ be a non-binary phylogenetic tree. Assume that any possible rank function on $\cT$ is equally likely.
With that assumption, we have
$$\bP[r(v)=i] = \frac {\alpha_{\cT,v}(i)}{|r(\cT)|}.$$
To calculate these probabilities, the algorithm {\sc RankProb} can be generalized to non-binary trees. We call the generalized algorithm {\sc RankProbGen}.
\\

\noindent
%\begin{algorithm} %[h!]
      %\begin{alg}
        {\bf Algorithm} {\sc RankProbGen} ($\cT,v$) \index{algorithm R{\sc ankProbGen}}\\
        {\bf Input}: A rooted phylogenetic tree $\cT$ and an interior vertex $v$.\\
        {\bf Output}: The probabilities $\bP[r(v)=i]$ for ${i=1, \ldots ,|\Vr|}$.
        \begin{algorithmic}[1]
          \STATE Denote the vertices of the path from $v$ to root $\rho$ with \\ $(v=x_1, x_2, \ldots ,x_n=\rho)$.
          \STATE Denote the subtree of $\cT$, consisting of root $x_m$ and all its descendants, by $\cT_m$ for $m=1, \ldots ,n$.
          \FOR{$m = 1, \ldots ,n$}
            \FOR{$i=1, \ldots ,|\Vr_{\cT}|$}
               \STATE $\tilde{\alpha}_{\cT_m,v}(i)=0$
            \ENDFOR
          \ENDFOR
          \STATE $\tilde{\alpha}_{\cT_1,v}(1)=1$
          \FOR{$m = 2, \ldots ,n$}
                    \STATE Label the subtree $\cT_m \setminus \cT_{m-1}$ by $\cT_{m-1}'$ ({\it cf.} Figure \ref{FigGen01})
                \STATE $M= \min \{ |\Vr_{\cT_{m-1}'}|-1, i-2 \}$
                \FOR{$i=m, \ldots ,|\Vr_{\cT_m}|$}
                        \STATE $\displaystyle \tilde{\alpha}_{\cT_m,v}(i):=\sum_{j=0}^{M} \tilde{\alpha}_{\cT_{m-1},v} (i-j-1)
                        {|\Vr_{\cT_{m-1}}|+|\Vr_{\cT_{m-1}'}|-1-(i-1) \choose |\Vr_{\cT_{m-1}'}|-1-j} {i-2 \choose j}$
                \ENDFOR
          \ENDFOR
          \FOR{$i=1, \ldots , |\Vr_{\cT}|$}
            \STATE $\bP[r(v)=i] = \frac{\tilde{\alpha}_{\cT_n,v}(i)}{\sum_{j} \tilde{\alpha}_{\cT_n,v}(j)}$
          \ENDFOR
          \STATE RETURN $\bP[r(v)=i], i=1,\ldots,|\Vr|$.
        \end{algorithmic}
        \bigskip
      %\end{alg}
%\end{algorithm}

\begin{figure}
\begin{center}
\input{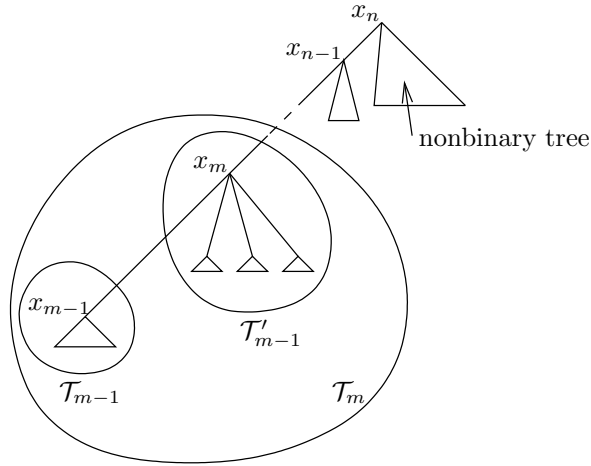}
\caption{Labelling the tree for algorithm {\sc RankProbGen}.}
\label{FigGen01}
\end{center}
\end{figure}

\begin{thm}
{\sc RankProbGen} returns the probabilities $$\bP[r(v)=i]$$
for each given $v\in \Vr$ and all $i \in 1, \ldots ,|\Vr|$.
The runtime is $O(|\Vr|^2)$.
\end{thm}
\begin{proof}
The algorithm is the same as {\sc RankProb}. The only difference is that in each step, we define $\cT_{m-1}' := \cT_m \setminus \cT_{m-1}$, i.e. the root of $\cT_m'$ is $x_m$.
For any rank function on $\cT_m'$, we now insert the first $j$ elements (excluding the root $x_m$) before the vertex $v$.
The number of ways to insert these vertices is counted analogously to the proof of Theorem (\ref{ThmRankCount}).
The number of possible rank functions on $\cT_m'$ does not have to be calculated, since these factors cancel out when calculating the probabilities.

Since we do the same iterations as in {\sc RankProb}, the algorithm {\sc RankProbGen} has quadratic runtime as well.
\end{proof}
\section{Comparing two interior vertices}
\begin{figure}
\begin{center}
\input{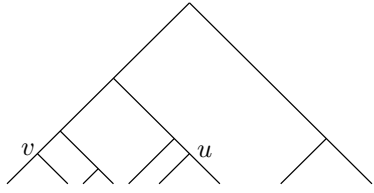}
\caption{What is the probability that vertex $u$ has smaller rank than vertex $v$?}
\label{FigCompare08}
\end{center}
\end{figure}

Assume again that every rank function on a binary phylogenetic tree $\cT$ is equally likely.
We want to compare two interior vertices $u$ and $v$ of $\cT$. Was $u$ more likely before $v$ or $v$ before $u$ ({\it cf.} Fig. \ref{FigCompare08})?
In other words, we want to know the probability
\label{ProbComp}
$$\bP_{u<v}:=\bP[r(u) < r(v) | \cT]$$
where $r(T)$ is the set of all possible rank functions on $\cT$.
This probability is, by Theorem (\ref{ThmYuleRankGivenT}), equivalent to counting all the possible rank functions on $\cT$ in
which $u$ has lower rank than $v$ and divide that number by all possible rank functions on $\cT$.
The algorithm {\sc Compare} will solve this problem in quadratic time.\\

\noindent
%\begin{algorithm} %[h!]
      %\begin{alg}
        {\bf Algorithm} {\sc Compare} ($\cT,u,v$) \index{algorithm C{\sc ompare}}\\
        {\bf Input}: A rooted phylogenetic tree $\cT$ and two distinct interior vertices $u$ and $v$.\\
        {\bf Output}: The probability $\bP_{u<v} :=\bP[r(u) < r(v) | \cT]$.
        \begin{algorithmic}[1]
          \STATE Denote the most recent common ancestor of $u$ and $v$ by $\rho_1$.
          \IF {$\rho_1 = v$}
            \STATE RETURN $\bP_{u<v} = 0$.
          \ENDIF
          \IF {$\rho_1 = u$}
            \STATE RETURN $\bP_{u<v} = 1$.
          \ENDIF
          \STATE Let $\cT_{\rho_1}$ be the subtree of $\cT$ which is induced by $\rho_1$.
          \STATE Delete the vertex $\rho_1$ from $\cT_{\rho_1}$. The two evolving subtrees are labeled $\cT_u$ and $\cT_v$ with $u \in \cT_u$ and $v \in \cT_v$.
          \STATE Run {\sc RankProb($\cT_u,u$)} and {\sc RankProb($\cT_v,v$)} to get $\bP[r(u)=i]$ on $\cT_u$ and $\bP[r(v)=i]$ on $\cT_v$ for all possible $i$.
          %\FOR{$i = |\Vr_{\cT_v}|+1, \ldots ,|\Vr_{\cT_u}|$}
          %      \STATE $\bP[r(v)=i]:=0$ 
           % \ENDFOR
            \FOR{$i = 1, \ldots ,|\Vr_{\cT_u}|$}
                \STATE $ucum(i) := \sum_{k=1}^i \bP[r(u)=i]$
            \ENDFOR
                \STATE $\bP_{u<v} := 0$
                \FOR{$i = 1, \ldots, |\Vr_{\cT_v}|$}
                \FOR{$j = 1, \ldots |\Vr_{\cT_u}|$}
                \STATE $p := \bP[r(v)=i] \cdot {i-1+j \choose j} \cdot {|\Vr_{\cT_v}|-i+|\Vr_{\cT_u}|-j \choose |\Vr_{\cT_u}|-j} \cdot ucum(j) \qquad (\ast)$
                    \STATE $\bP_{u<v} := \bP_{u<v}+p$
                \ENDFOR
            \ENDFOR
                \STATE $tot := {|\Vr_{\cT_u}|+|\Vr_{\cT_v}| \choose |\Vr_{\cT_v}|}$
                \STATE $\bP_{u<v} := \bP_{u<v}/tot$
                \STATE RETURN $\bP_{u<v}$
          \end{algorithmic}
        \bigskip
      %\end{alg}
%\end{algorithm}

\begin{thm}
The algorithm {\sc Compare} returns the value $$\bP_{u<v} =\bP[r(u) < r(v) | \cT].$$
\end{thm}
\begin{proof}
Note that the probability of $u$ having smaller rank than $v$ in tree $\cT_{\rho_1}$ equals the probability of $u$ having smaller
rank than $v$ in tree $\cT$, since for any rank function on $\cT_{\rho_1}$, there is the same number of linear extensions to get a rank function on the tree $\cT$.

So it is sufficient to calculate the probability $\bP_{u<v}$ in $\cT_{\rho_1}$.
If $\rho_1=u$, $u$ is before $v$ in $\cT$ and we return $\bP_{u<v}=1$.
If $\rho_1=v$, $v$ is before $u$ in $\cT$ and we return $\bP_{u<v}=0$.

In the following, let $\rho_1 \neq u, \rho_1 \neq v$.
The run of {\sc RankProb} gives us the probability $\bP[r(u)=i]$ in the tree $\cT_u$ and $\bP[r(v)=i]$ in $\cT_v$ for all $i$.
We want to combine these two linear orders. Assume that $r(v)=i$ and we insert $j$ vertices of $\cT_u$ before $v$.
Inserting $j$ vertices of $\cT_u$ into the linear order of $\cT_v$ before $v$ is possible in
${i-1+j \choose j}$ ways (see Corollary \ref{CorSequenceBin}).
Putting the remaining vertices in a linear order is possible in
${|\Vr_{\cT_v}|-i+|\Vr_{\cT_u}|-j \choose |\Vr_{\cT_u}|-j}$ ways.
The probability that the vertex $u$ is among the $j$ vertices which have smaller rank than $v$ is $\bP[r(u) \leq j] = ucum(j)$. There are $|r(\cT_u)|$ possible linear orders on $\cT_u$ and $|r(\cT_v)|$ possible linear orders on $\cT_v$.
The number of linear orders where vertex $v$ has rank $i$ in $\cT_v$, $v$ has rank $i+j$ in $\cT_{\rho_1}$ and $r(u)<i+j$ therefore equals
$$p'_{i,j} = \bP[r(v)=i] \cdot |r(\cT_v)| \cdot {i-1+j \choose j} \cdot {|\Vr_{\cT_v}|-i+|\Vr_{\cT_u}|-j \choose |\Vr_{\cT_u}|-j} \cdot ucum(j) \cdot |r(\cT_u)|$$
Adding up the $p'$ for each $i$ and $j$ gives us the number of linear orders where $u$ is earlier than $v$.

Combining a linear order on $\cT_v$ with a linear order on $\cT_u$ is possible in
$$tot := {|\Vr_{\cT_u}|+|\Vr_{\cT_v}| \choose |\Vr_{\cT_v}|}$$
different ways (see Corollary \ref{CorSequenceBin}). There are $|r(\cT_u)|$ linear orders on $\cT_u$ and $|r(\cT_v)|$ linear orders on $\cT_v$, so on $\cT_{\rho_1}$, we have
$$tot' := {|\Vr_{\cT_u}|+|\Vr_{\cT_v}| \choose |\Vr_{\cT_v}|} |r(\cT_v)| |r(\cT_v)|$$
linear orders.
Therefore we get
$$\bP_{u<v} = \frac{\sum_{i,j} p'_{i,j}}{tot'} = \frac{\sum_{i,j} p_{i,j}}{tot}$$
with $p_{i,j} = \bP[r(v)=i] \cdot {i-1+j \choose j} \cdot {|\Vr_{\cT_v}|-i+|\Vr_{\cT_u}|-j \choose |\Vr_{\cT_u}|-j} \cdot ucum(j)$. This shows that {\sc Compare} works correct.
\end{proof}
\begin{thm}
The runtime of {\sc Compare} is $O(|\Vr|^2)$.
\end{thm}
\begin{proof}
Again, note that the combinatorial factors ${n \choose k}$ for all $n,k \leq |\Vr|$ can be calculated in advance in quadratic time, see Remark (\ref{nchoosek}).
In the algorithm, those factors can then be obtained in constant time.

Contributions to the runtime from each line in {\sc Compare} (the runtime is always w.r.t. $|\Vr|$):\\
{\it Line 1:} linear time\\
{\it Line 2--7:} constant time\\
{\it Line 8:} linear time\\
{\it Line 9:} constant time\\
{\it Line 10:} quadratic time, since {\sc RankProb} has quadratic runtime\\
{\it Line 11--13:} linear time\\
{\it Line 14:} constant time\\
{\it Line 15--20:} quadratic time since $(\ast)$ has to be evaluated $|\Vr_{\cT_u}| \cdot |\Vr_{\cT_u}| \leq |\Vr_{\cT}|^2$ times\\
{\it Line 21--23:} constant time

Therefore, the overall runtime of {\sc Compare} is $O(|\Vr|^2)$.
\end{proof}

\begin{figure}
\begin{center}
\input{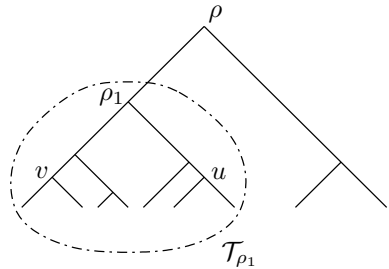}
\caption{Example for {\sc Compare}: Calculate the probability of $u<v$ in the displayed tree $\cT$.}
\label{FigCompare06}
\end{center}
\end{figure}

\begin{example}
Fig. \ref{FigCompare06} displays the tree $\cT$. We want to calculate the probability $\bP_{u<v}$, i.e. the probability
of vertex $u$ having a smaller rank than vertex $v$.\\

A run of the Python code attached in Appendix \ref{PythonCode} with input $(\cT, u, v)$ returns
$\bP_{u<v} = \frac {9}{20}$.
\end{example}

\section{Application of {\sc RankProb} - Estimating edge lengths in a Yule tree} \label{EstEdgeLength}

In \cite{Vos2006}, a primate supertree on 218 species was constructed with the MRP method (Matrix Representation using Parsimony analysis, see \cite{Baum1992,Ragan1992}). \index{supertree} \index{primates}
The resulting supertree is shown in Appendix \ref{Primates}. This tree has only 210 interior vertices. There are
six `soft' polytomies in the supertree, i.e. six vertices have more than two direct descendants because the exact resolution is unclear (i.e. the supertree is non-binary). \index{polytomy}

Since for most of the interior vertices, no molecular estimates were available, the edge lengths for the tree were estimated. 
Here, the length of an edge represents the
time between two speciation events.

A very common stochastic model for trees with edge lengths is the continuous-time Yule model. \index{Yule model!continuous-time} As in the discrete-time Yule model, at every point in time, each species is equally likely to split and give birth to two new species. The expected waiting time for the next speciation event in a tree with $n$ leaves is $1/n$. That is, each species at any given time has a constant speciation rate (normalized so that 1 is the expected time until it next speciates).

It was assumed that the primate tree $\cT_p$ evolved under the continuous-time Yule model.
%We have seen in example \ref{ExPrimatesYule} that in discrete time, the assumption of $\cT_p$ having evolved under the Yule model is reasonable.
%
In \cite{Vos2006}, $10^6$ rank functions on $\cT_p$ were drawn uniformly at random. For each of those rank functions, the expected time intervals, i.e. the edge lengths, between vertices were considered (the expected waiting time after the $(n-1)$th event until the $n$th event is $1/n$).

The authors of \cite{Vos2006} concluded their paper by asking for an analytical approach to the estimation of the edge length, and we provide this now.

\subsection{Analytical estimation of the edge length}

Let $(u,v)$ be an interior edge in $\cT$ with $u<_{\cT}v$.
Let $X$ be the random variable `length of the edge $(u,v)$' given that $\cT$ is generated according to the continuous-time Yule model.
%Let $\Omega = \{(i,j): i<j, i,j \in \{1, \ldots, n-1\} \}$ and let $A_{i,j}$ be the event that $r(u)=i, r(v)=j$.

\begin{figure}
\begin{center}
\input{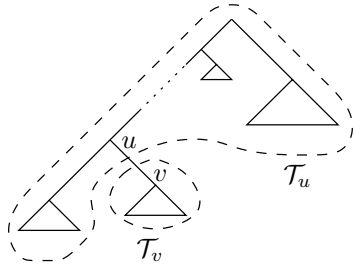}
\caption{Labeling the tree for estimating the edge lengths.}
\label{TreeEstEdge}
\end{center}
\end{figure}

The expected length $\bE[X]$ of the edge $(u,v)$ is given by
$$\bE[X] = \sum_{i,j} \bE[X|r(u)=i, r(v)=j] \bP[r(u)=i, r(v)=j].$$
Since under the continuous-time Yule model, the expected waiting time for the next event is $1/n$, we have
$$\bE[X|r(u)=i, r(v)=j] = \sum_{k=1}^{j-i} \frac{1}{i+k}.$$
It remains to calculate the probability $\bP[r(u)=i, r(v)=j]$.
We count all the possible rank functions where $r(u)=i$ and $r(v)=j$.
The subtree $\cT_v$ consists of $v$ and all its descendants.
The tree $\cT_u$ evolves from $\cT$ when we replace the subtree $\cT_v$ by a leaf, see Fig. \ref{TreeEstEdge}.

Note that $\bP[r(u)=i, r(v)=j]=0$ if $|\Vr_{\cT_u}|<j-1$.
Therefore, assume $|\Vr_{\cT_u}| \geq j-1$ in the following.

The number of rank functions in $\cT_u$ is denoted by $R_{\cT_u}$.
The probability $\bP[r(u)=i]$ can be calculated with {\sc RankProb}($\cT_u$, $u$).
So the number of rank functions in $\cT_u$ with $\bP[r(u)=i]$ is $\bP[r(u)=i]\cdot R_{\cT_u}$.

The number of rank functions in $\cT_v$ is denoted by $R_{\cT_v}$.
Let any linear order on the tree $\cT_u$ and $\cT_v$ be given.
Combining those two linear orders to an order on $\cT$, where $r(v)=j$ holds, means, that the vertices with rank $1,2,\ldots,j-1$ in $\cT_u$ keep their rank. Vertex $v$ gets rank $j$. The remaining $|\Vr_{\cT_u}|-(j-1)$ vertices
in $\cT_u$ and $|\Vr_{\cT_v}|-1$ vertices in $\cT_v$ have to be shuffled together.
According to Corollary (\ref{CorSequenceBin}), this can be done in
$$ {|\Vr_{\cT_u}|-(j-1) + |\Vr_{\cT_v}|-1 \choose |\Vr_{\cT_v}|-1} = {|\Vr_{\cT_u}| + |\Vr_{\cT_v}|-j \choose |\Vr_{\cT_v}|-1}$$
different ways. Overall, we have
$$\bP[r(u)=i]\cdot R_{\cT_u} \cdot R_{\cT_v} \cdot {|\Vr_{\cT_u}| + |\Vr_{\cT_v}|-j \choose |\Vr_{\cT_v}|-1}$$
different rank functions on $\cT$ with $r(u)=i$ and $r(v)=j$.
For the probability $\bP[r(u)=i, r(v)=j]$, we get
$$\bP[r(u)=i, r(v)=j] = \frac {\bP[r(u)=i]\cdot R_{\cT_u} \cdot R_{\cT_v} \cdot {|\Vr_{\cT_u}| + |\Vr_{\cT_v}|-j \choose |\Vr_{\cT_v}|-1}}{\sum_{i,j} \bP[r(u)=i]\cdot R_{\cT_u} \cdot R_{\cT_v} \cdot {|\Vr_{\cT_u}| + |\Vr_{\cT_v}|-j \choose |\Vr_{\cT_v}|-1}}$$
Since $R_{\cT_u}$ and $R_{\cT_v}$ are independent of $i$ and $j$, those factors cancel out, and we get
$$\bP[r(u)=i, r(v)=j] = \frac {\bP[r(u)=i]\cdot {|\Vr_{\cT_u}| + |\Vr_{\cT_v}|-j \choose |\Vr_{\cT_v}|-1}}{\sum_{i,j} \bP[r(u)=i] \cdot {|\Vr_{\cT_u}| + |\Vr_{\cT_v}|-j \choose |\Vr_{\cT_v}|-1}}$$
Further, we note that
$${|\Vr_{\cT_u}| + |\Vr_{\cT_v}|-j \choose |\Vr_{\cT_v}|-1} = \frac{(|\Vr_{\cT}|-j)!}{(|\Vr_{\cT_v}|-1)!(|\Vr_{\cT}|-j-(|\Vr_{\cT_v}|-1))!}$$
Again, since $(|\Vr_{\cT_v}|-1)!$ is independent of $i$ and $j$, this factor cancels out, and we are left with
$$\bP[r(u)=i, r(v)=j] = \frac {\bP[r(u)=i]\cdot \prod_{k=0}^{|\Vr_{\cT_v}|-2} (|\Vr_{\cT}|-j-k) }{\sum_{i,j} \bP[r(u)=i] \cdot \prod_{k=0}^{|\Vr_{\cT_v}|-2} (|\Vr_{\cT}|-j-k)}$$
Let $\Omega = \{(i,j): i<j, i,j \in \{1, \ldots, |\Vr| \}, |\Vr_{\cT_u}| \geq j-1\}$. With that notation, the expected edge length $\bE[X]$ is
\begin{eqnarray}
\bE[X] &=& \sum_{(i,j) \in \Omega}  \bE[X|r(u)=i, r(v)=j] \bP[r(u)=i, r(v)=j] \notag \\
&=& \sum_{(i,j) \in \Omega} \left[ \left( \sum_{k=1}^{j-i} \frac{1}{i+k} \right) \frac {\bP[r(u)=i]\cdot \prod_{k=0}^{|\Vr_{\cT_v}|-2} (|\Vr_{\cT}|-j-k) }{\sum_{(i,j) \in \Omega} \left[ \bP[r(u)=i] \cdot \prod_{k=0}^{|\Vr_{\cT_v}|-2} (|\Vr_{\cT}|-j-k) \right]} \right] \notag \\
&=& \frac {\sum_{(i,j) \in \Omega} \left[ \left( \sum_{k=1}^{j-i} \frac{1}{i+k} \right) \cdot \bP[r(u)=i]\cdot \prod_{k=0}^{|\Vr_{\cT_v}|-2} (|\Vr_{\cT}|-j-k)  \right]}{\sum_{(i,j) \in \Omega} \left[ \bP[r(u)=i] \cdot \prod_{k=0}^{|\Vr_{\cT_v}|-2} (|\Vr_{\cT}|-j-k) \right]}  \label{EqnExpEdgeLength}
\end{eqnarray}

\bigskip
\begin{rem}
With Equation (\ref{EqnExpEdgeLength}), we can estimate the length of all the interior edges. For the pendant edges, the approach above gives us no estimate though. All we know is that the time from the latest interior vertex, which has rank $n-1$, until the presence is expected to be at most $1/n$ where $n$ is the number of leaves.
\end{rem}
\begin{rem}
In a supertree, we can have interior vertices which are not fully resolved, i.e. an interior vertex can have more than two descendants, because the exact resolution is unclear. Our calculation for the expected edge length assumes a binary tree though.

However, we can calculate the expected edge length for each possible binary resolution of the supertree.
Assume the supertree $\cT$ has the possible binary resolutions $\cT_1, \ldots, \cT_m$.
For an edge $(u,v)$ in $\cT$ where $u<_{\cT}v$, the expected edge length is calculated in the trees $\cT_i$ for $i = 1, \ldots, m$. The expected edge length in $\cT_i$ is denoted by $e_i$ for $i = 1, \ldots, m$.

We calculate the expected edge length $\bE[X]$ of $(u,v)$ in the supertree $\cT$ by
\begin{equation}
\bE[X]=\frac{\sum_i e_i \bP[\cT_i]}{\sum_i \bP[\cT_i]} \label{EqnWeightSum}
\end{equation}
where the probability $\bP[\cT_i]$ is calculated according to Corollary (\ref{CorProbYule}).

Note that if $u$ is a vertex with more than two descendants in $\cT$, $v$ is in general not a direct descendant of $u$ in $\cT_i$. The value $e_i$ in resolution $\cT_i$ is then the sum of all expected edge lengths on the path from $u$ to $v$ in $\cT_i$.
\end{rem}

\begin{rem}
In the primate supertree in Appendix \ref{Primates}, there are six interior vertices with more than two descendants (vertex labels $48,63,148, 153,157$ and $200$). For the vertices labeled with $63$ and $200$, only one resolution is possible (up to the labeling).

The interior vertices with label $48$, $153$ and $157$ have three descendants each. So there are $3^3$ possible binary resolutions. The interior vertex $148$ has four leaf-descendants. There are two possible binary resolutions (up to the labeling). To calculate the expected edge lengths for the primate supertree, we therefore have to calculate the expected edge lengths on $3^3 \cdot 2$ binary trees and then calculate the weighted sum from Equation (\ref{EqnWeightSum}).
\end{rem}

\chapter{Speciation Rates} \label{SpeciRate}
This chapter was motivated by Craig Moritz and Andrew Hugall, biologists from Berkeley and Adelaide. They looked at a tree showing the relationships between a set of snails. Each of those snails lives either in rain forest or open forest. The tree has edge lengths assigned.
Moritz and Hugall asked if the rate of speciation \index{rate of speciation} is different for rain forest snails and open forest snails.
%The rate of speciation is the inverse of the average edge length for a class of species. So the rate of speciation is a measure of how fast a class of species produces splits in the evolutionary tree.

Mathematically, determining the rate of speciation is the following problem. The leaves are divided into two classes, $\alpha$ and $\beta$
(e.g. rain forest and open forest snails). Given the rate that a species belonging to class $\alpha$ changes to a species belonging to class $\beta$ (and vice versa),
we calculate the expected length of an edge between two species of group $\alpha$ (resp. $\beta$). This expected length is
an estimate for the inverse of the rate of speciation and is calculated in linear time.

\section{Some notation}
%schon in datei diplom1
\begin{defi} \index{character} \index{character!full} \index{character!binary} \label{DefChar}
Let $X'$ be a non-empty subset of $X$. Let $C$ be a non-empty set.
A $character$ on $X$ is a function $\chi: X' \rightarrow C$. $C$ is the {\it character state set} of $\chi$. \index{character state set}
If $X'=X$, we say $\chi$ is a $full~character$. If $|C|=2$, we say $\chi$ is a $binary~character$.
\end{defi}

\begin{figure}
\begin{center}
\input{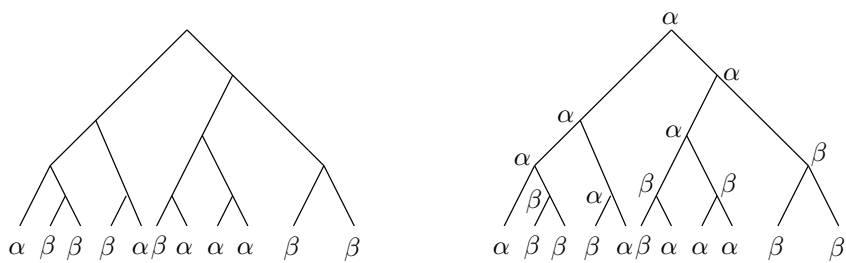}
\caption{A phylogenetic tree with a full character on the left and a phylogenetic state tree on the right (without the leaf labels).} %\label{FigTtilde}
\end{center}
\end{figure}

%Let X be a non-empty set.
\begin{defi} \index{phylogenetic state tree}
Let $\cT$ be a rooted phylogenetic $X$-tree with vertex set $V$
and leaf set $L \subset V$. Let $\chi$ be a full binary character
on $\cT$, $\chi: X \rightarrow \{\alpha, \beta\}$.
%Since in a phylogenetic $X$-tree, the map from $X$ to $L$ is a bijection, we can as well define $\chi$ as $\chi: L \rightarrow \{\alpha, \beta\}$.\\
Define $s: V \rightarrow \{\alpha, \beta\}$ with $s|_L = \chi
\circ \phi^{-1}$. $(\cT,s)$ is called a $phylogenetic~state~tree$, $s$ a $state$ $function$. \label{DefStateFunc} \index{state function}
\end{defi}
\label{DefEdgeLength} \index{edge!length} In the following, the phylogenetic state tree $(\cT,s)$ shall have
assigned a function $l: E \rightarrow \mathbb{R}^+$. $l$ shall
denote the edge lengths of $\cT$. Let $ \eta \in \{ \alpha, \beta
\}$ throughout this chapter. Let $v$ be any node in $(\cT,s)$ with
$s(v)=\eta$. We then say that the $state$ $of$ $v$ is $\eta$. \index{state of vertex} Let
$ \gamma \in \{ \alpha, \beta \} \times \{ \alpha, \beta \} $
throughout the chapter, i.e. $\gamma = (\gamma_1, \gamma_2)$ with
$\gamma_1, \gamma_2 \in  \{ \alpha, \beta \}$. An edge $e=(v_1,v_2)$ of
$(\cT,s)$ where $v_1 <_{\cT} v_2$ and
$s(v_1)=\gamma_1$, $s(v_2)=\gamma_2$ is called a $\gamma-edge$. \index{$\gamma$-edge}
\begin{figure}
\begin{center}
\input{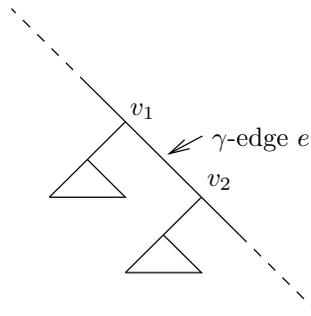}
\caption[Example of a $\gamma$-edge]{With $s(v_1)=\gamma_1$ and $s(v_2)=\gamma_2$, the edge $e=(v_1,v_2)$ is a $\gamma$-edge.} \label{FigGammaEdge}
\end{center}
\end{figure}
\section{Markov Chain Model} \index{Markov Chain model}

\begin{figure}
\begin{center}
\input{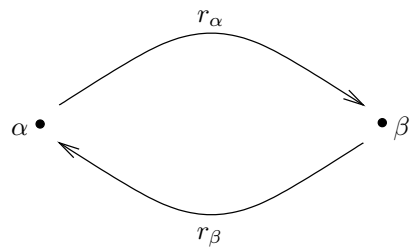}
\caption{Rate of the state change for a binary character}
\label{FigRate}
\end{center}
\end{figure}

Throughout evolution, assume that state $\alpha$ changes to state $\beta$ with rate $r_\alpha$ and state $\beta$ changes to state $\alpha$
with rate $r_\beta$, \label{DefRateState}
so the rates only depend upon the state of the last vertex (see Fig. \ref{FigRate}).
This means that the state change follows a Markov Chain model, and for that model, we want to calculate the transision matrix \index{transition matrix} \label{TransMat}
\[P(l(e))=
\begin{pmatrix}
    p_{\alpha \alpha}(l(e)) & p_{\alpha \beta}(l(e)) \\
    p_{\beta \alpha}(l(e)) & p_{\beta \beta}(l(e))
\end{pmatrix}
\]
where $p_{ \gamma_1 \gamma_2}(l(e)) = \bP \left[ (s(v_2)=
\gamma_2) | (s(v_1)= \gamma_1) \right]$ with $e=(v_1,v_2)$ and $v_1 <_{\cT} v_2$.

The rate matrix $R$ is defined as \label{RateMatrix} \index{rate matrix}
\[
R=\begin{pmatrix}
    -r_\alpha & r_\alpha \\
    r_\beta & -r_\beta
\end{pmatrix}
\]
Diagonalization of R yields
\[ R =
\begin{pmatrix}
    -r_\alpha & r_\alpha \\
    r_\beta & -r_\beta
\end{pmatrix} = S
\begin{pmatrix}
    0 & 0 \\
    0 & -(r_\alpha+r_\beta)
\end{pmatrix} S^{-1}
\]
with
\[S=
\begin{pmatrix}
    1 & r_\alpha \\
    1 & -r_\beta
\end{pmatrix}
\]
From stochastic processes, we know that the connection between the rate matrix and the transition matrix is
$$P'(l(e))=R P(l(e))$$
Solving this differential equation yields
$$P(l(e))=P(0) e^{R (l(e))}$$
with $P(0) = Id$ since $l(e)=0$ means staying in the vertex. Therefore $P(l(e))$ can be rewritten as
\begin{eqnarray}
P(l(e)) &=& e^{R (l(e))} \notag \\
&=& \exp \{ S
\begin{pmatrix}
    0 & 0 \\
    0 & -(r_\alpha+r_\beta)
\end{pmatrix} S^{-1} l(e)
\} \notag \\
&=& S \exp \{{ \begin{pmatrix}
    0 & 0 \\
    0 & -(r_\alpha+r_\beta)
\end{pmatrix} l(e) }\} S^{-1} \notag \\
&=& S \begin{pmatrix}
    1 & 0 \\
    0 & e^{-(r_\alpha+r_\beta)l(e)}
\end{pmatrix} S^{-1} \notag \\
&=&
\begin{pmatrix}
    \frac {1}{r_\alpha+r_\beta} \left( r_\beta+r_\alpha e^{-(r_\alpha+r_\beta)l(e)} \right) &
    \frac {r_\alpha}{r_\alpha+r_\beta} \left( 1- e^{-(r_\alpha+r_\beta)l(e)} \right) \\
    \frac {r_\beta}{r_\alpha+r_\beta} \left( 1- e^{-(r_\alpha+r_\beta)l(e)} \right) &
    \frac {1}{r_\alpha+r_\beta} \left( r_\alpha + r_\beta e^{-(r_\alpha+r_\beta)l(e)} \right) \notag
\end{pmatrix}
\end{eqnarray}
The initial probability of vertex $v$ being in state $\eta$ shall be $\pi_\eta$, $\eta \in \{ \alpha, \beta \}$. \label{InitProb} \index{initial probability distribution} It holds
$$\begin{pmatrix} \pi_\alpha & \pi_\beta \end{pmatrix} R = \begin{pmatrix} \pi_\alpha & \pi_\beta \end{pmatrix} \begin{pmatrix}
    -r_\alpha & r_\alpha \\
    r_\beta & -r_\beta
\end{pmatrix} = 0$$
so $$\pi = \begin{pmatrix} \pi_\alpha & \pi_\beta \end{pmatrix} =
\begin{pmatrix} \frac{r_\beta}{r_\alpha+r_\beta} & \frac{r_\alpha}{r_\alpha+r_\beta} \end{pmatrix}$$

Therefore, for any given phylogenetic tree $\cT$ with edge lengths
$l(e)$, the probability of its vertices being in states according
to a state function $s$ is
\begin{equation} \bP [s] = \pi_{s(\rho)} \prod_{\substack{e \in E \\ e=(v_1,v_2) \\ v_1 <_{\cT} v_2}} p_{s(v_1),s(v_2)} \label{EqnPTree} \end{equation}
Furthermore, it holds for any $e \in E$ with $e=(v_1,v_2)$
\begin{equation} p_{s(v_1),s(v_2)}(l(e)) = \frac {r_{s(v_1)}}{r_{s(v_2)}} p_{s(v_2),s(v_1)} (l(e)) \label{EqnPTrafo} \end{equation}

\section{Expected length of a $\gamma$-edge}
Given a phylogenetic tree $\cT$ with character $\chi$, edge length
$l(e)$ and rate matrix $R$, we want to calculate the expected
average length of a $\gamma$-edge over all $(\cT,s)$. The inverse
of this length is an estimate for the rate of speciation.

\begin{figure}
\begin{center}
\input{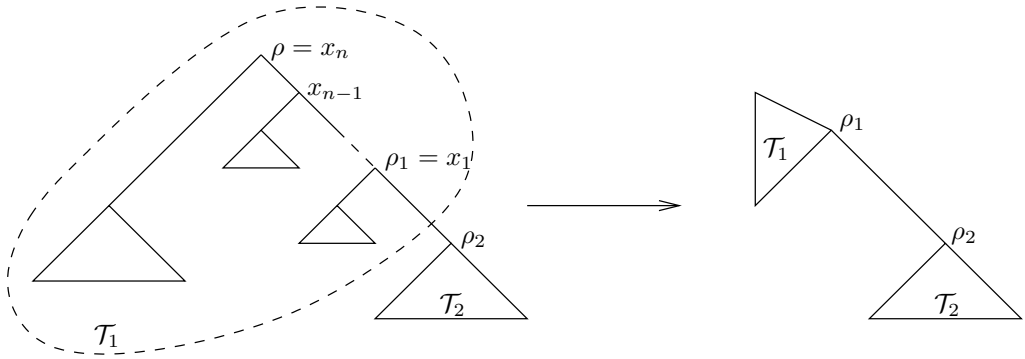}
\caption{Calculating the expected edge length: Defining $\cT_1$
and $\cT_2$} \label{FigNewRoot}
\end{center}
\end{figure}

Calculating the expected average length of a $\gamma$-edge over
all $(\cT,s)$ means calculating
$$ \bE_\chi \left[ \frac{\displaystyle \sum_{e\in E,~e~\gamma-edge}l(e)}{\#~of~\gamma-edges}\right]$$
where $ \bE_\chi$ denotes the expected value over all $s$ given
$s|_{L} = \chi$.
Trying to calculate this expected value turns out to give us very nasty recursion formulas.\\
So we change the problem slightly and try to calculate instead
$$\Psi_{\gamma}=\frac{\displaystyle \bE_\chi \left[ \sum_{e\in E,~e~\gamma-edge}l(e) \right]}{\bE_\chi \left[ \#~of~\gamma-edges \right]}$$ \label{CalcPsi}
Define the random variable
$$X_{\gamma}(e):=\left\{
\begin{array}{ll}
    1 & \hbox{if $e$ is $\gamma$-edge} \\
    0 & \hbox{else} \\
\end{array}
\right.$$
With that, we get
\begin{eqnarray}
\Psi_{\gamma}&=&\frac{\displaystyle \bE_\chi \left[ \sum_{e\in E,~e~\gamma-edge}l(e) \right]}{\bE_\chi \left[ \#~of~\gamma-edges \right]} \notag \\
&=& \frac{\displaystyle \bE_\chi \left[ \sum_{e \in E}l(e) X_\gamma (e) \right]}{\bE_\chi \left[\displaystyle  \sum_{e \in E} X_\gamma (e) \right]} \notag \\
&=& \frac{\displaystyle  \sum_{e \in E}l(e) \bP \left[  (X_\gamma(e) = 1) | \chi\right] }{\displaystyle  \sum_{e \in E}
\bP \left[  (X_\gamma(e) = 1) | \chi \right] } \label{EqnRho}
\end{eqnarray}
where $\bP \left[  (X_\gamma(e) = 1) | \chi \right]$ denotes the probability of
$e$ being a $\gamma$-edge given $s|_L = \chi$. So it is basically left to
calculate $\bP \left[ (X_\gamma(e) = 1) | \chi \right]$. To do so,
we first define two subtrees of $\cT$ (see also Fig.
\ref{FigNewRoot}). Denote the end vertices of $e$ by $\rho_1$ and
$\rho_2$ with $\rho_1 <_{\cT} \rho_2$. By deleting the
$\gamma$-edge $e$ in $\cT$, we get two new trees $\cT_1$ and
$\cT_2$, $\cT_1$ with $\rho_1 \in \cT_1$ and character $\chi_{1} =
\chi|_{\phi^{-1}(L_{\cT_1})}$, and $\cT_2$ with $\rho_2 \in \cT_2$
and character $\chi_{2} = \chi|_{\phi^{-1}( L_{\cT_2})}$ where $L_{\cT_i}$ denotes the set of leaves of $\cT_i$, $i \in \{1,
2\}$. The root in $\cT_i$ shall be $\rho_i$, so $\rho$ becomes
an ordinary vertex in $\cT_1$.

$\bP \left[\chi_i|(s(\rho_i)=\gamma_i) \right]$ shall denote the
probability of the character $\chi_i$ on the tree $\cT_i$ given $s(\rho_i)=\gamma_i$.
$\bP \left[\chi_{\cT \setminus \cT_2}|(s(\rho_1)=\gamma_1) \right]$ shall
denote the probability of the character $\chi_{\cT \setminus \cT_2}$ on the tree $\cT \setminus \cT_2$ given $s(\rho_1)=\gamma_1$. $\bP[\chi_{\cT},s]$ shall denote the probability of the character $\chi$ and the state function $s$ on the tree $\cT$. We denote the vertices on
the path from $\rho_1$ to $\rho$ by $\rho_1=x_1, x_2, \ldots
,x_{n-1}, x_n = \rho$. With (\ref{EqnPTree}) and
(\ref{EqnPTrafo}), it holds
\begin{eqnarray}
\bP \left[\chi_1,s \right] 
&=& \frac{\pi_{s(\rho_1)} \prod_{i=1}^{n-1} p_{s(x_i),s(x_{i+1})}}
{\pi_{s(\rho)} \prod_{i=1}^{n-1} p_{s(x_{i+1}),s(x_i)}} 
\bP \left[ \chi_{\cT \setminus \cT_2},s \right] \notag \\
&=& \frac{\pi_{s(\rho_1)} \prod_{i=1}^{n-1} p_{s(x_i),s(x_{i+1})}} {\pi_{s(\rho)} \prod_{i=1}^{n-1} \frac {r_{s(x_{i+1})}}
{r_{s(x_{i})}} p_{s(x_i),s(x_{i+1})}} \bP \left[ \chi_{\cT \setminus \cT_2},s \right] \notag \\
&=& \frac{\pi_{s(\rho_1)} r_{s(x_1)}} {\pi_{s(\rho)} r_{s(x_{n})}} \bP \left[ \chi_{\cT \setminus \cT_2},s \right] \notag \\
&=& \frac{ \frac{r_\alpha r_\beta} {r_\alpha + r_\beta} } { \frac{r_\alpha r_\beta} {r_\alpha + r_\beta}}
\bP \left[ \chi_{\cT \setminus \cT_2},s \right] \notag \\
&=& \bP \left[ \chi_{\cT \setminus \cT_2},s \right] \notag
\end{eqnarray}
This yields
\begin{eqnarray}
\bP \left[\chi_1|(s(\rho_1)=\gamma_1) \right] 
&=& \sum_{s: s(\rho_1)=\gamma_1} \bP \left[\chi_1,s \right] \notag \\
&=& \sum_{s: s(\rho_1)=\gamma_1} \bP \left[ \chi_{\cT \setminus \cT_2},s \right] \notag \\
&=& \bP \left[ \chi_{\cT \setminus \cT_2}|(s(\rho_1)=\gamma_1) \right] \notag
\end{eqnarray}
With that result, we get
\begin{align}
%& & \bP \left[ (X_\gamma(e) = 1) | \chi \right] \notag \\
%& & \notag \\
& \bP \left[ (X_\gamma(e) = 1) | \chi \right] = \notag \\
\notag \\
&= \frac{ \bP \left[(X_\gamma(e) = 1) \right]  \bP \left[\chi|(X_\gamma(e) = 1) \right] }
 { \bP \left[ \chi \right]  } \notag \\
&= \frac {  \pi_{\gamma_1} p_{\gamma_1 \gamma_2}(l(e))
 \bP \left[ \chi_{\cT \setminus
\cT_2}|(s(\rho_1)=\gamma_1) \right] \bP
\left[\chi_2|(s(\rho_2)=\gamma_2) \right] } {\displaystyle
\sum_{\gamma=(\gamma_1,\gamma_2)} \pi_{\gamma_1} p_{\gamma_1
\gamma_2}(l(e))
\bP \left[\chi_{\cT_1}|(s(\rho_1)=\gamma_1) \right] \bP \left[\chi_{\cT_2}|(s(\rho_2)=\gamma_2) \right]   } \notag \\
&= \frac{\pi_{\gamma_1} p_{\gamma_1 \gamma_2}(l(e)) \bP
\left[\chi_1|(s(\rho_1)=\gamma_1) \right] \bP
\left[\chi_2|(s(\rho_2)=\gamma_2) \right]} {\displaystyle
\sum_{\gamma=(\gamma_1,\gamma_2)} \pi_{\gamma_1} p_{\gamma_1
\gamma_2}(l(e)) \bP \left[\chi_{\cT_1}|(s(\rho_1)=\gamma_1)
\right] \bP \left[\chi_{\cT_2}|(s(\rho_2)=\gamma_2) \right]   }
\label{EqnPx}
\end{align}
$\bP \left[\chi_{i}|(s(\rho_i)=\gamma_i) \right]$ is calculated in a recursive way, starting from the bottom of the tree.
\begin{figure}
\begin{center}
\input{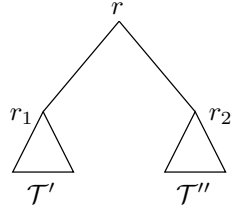}
\caption{Calculating the expected edge length: Defining
$\tilde{\cT}$} \label{FigTtilde2}
\end{center}
\end{figure}

Suppose we have the subtree $\tilde{\cT}$ as in Fig. \ref{FigTtilde2} and either $r_1$, $r_2$ are leaves or we know
$\bP \left[\chi_{\cT'}|(s(r_1)=\eta) \right]$ on tree $\cT'$,
$\bP \left[\chi_{\cT''}|(s(r_2)=\eta) \right]$ on tree $\cT''$,
for $\eta \in \{ \alpha, \beta \}$.
With that, we get the following $recursive$ $formulas$ for the probabilities on tree $\tilde{\cT}$.
\begin{itemize}
\item For $r_1$ and $r_2$ leaves:
$$\bP \left[\chi_{\tilde{\cT}}|(s(r)=\eta) \right] = \frac{p_{\eta \chi(r_1)} p_{\eta \chi(r_2)}}
{\displaystyle \sum_{\eta_1, \eta_2 \in \{\alpha, \beta\}} p_{\eta \eta_1} p_{\eta \eta_2}}$$
\item For $r_1$ leave, $r_2$ interior node:
$$\bP \left[\chi_{\tilde{\cT}}|(s(r)=\eta) \right] = \frac{\displaystyle \sum_{\eta_1 \in \{\alpha, \beta\}}
\bP \left[\chi_{\cT'}|(s(r_1)=\eta_1) \right] p_{\eta \chi(r_2)} p_{\eta \eta_1}}
{\displaystyle \sum_{\eta_1, \eta_2 \in \{\alpha, \beta\}} \bP \left[\chi_{\cT'}|(s(r_1)=\eta_1) \right]
p_{\eta \eta_2} p_{\eta \eta_1} } $$
\item For $r_1$ and $r_2$ interior nodes:
\begin{align}
&\bP \left[\chi_{\tilde{\cT}}|(s(r)=\eta) \right] = \notag \\
\notag \\
&{\displaystyle \sum_{\eta_1, \eta_2 \in \{\alpha, \beta\}}
\bP \left[\chi_{\cT'}|(s(r_1)=\eta_1) \right]
\bP \left[\chi_{\cT''}|(s(r_2)=\eta_2) \right]
p_{\eta \eta_1} p_{\eta \eta_2}} \notag
\end{align}
\end{itemize}
\index{algorithm E\textsc{dgeLength}}
\noindent
{\bf Algorithm} \textsc{EdgeLength} ($\cT, \chi$)\\
{\bf Input}: A rooted binary phylogenetic tree $\cT$ and a character $\chi$ on $\cT$ with state change rates $r_\alpha$ and $r_\beta$\\
{\bf Output}: The values $\Psi_\gamma$ for $\gamma \in \{\alpha, \beta\} \times \{\alpha, \beta\}$ ({\it cf.} Equation (\ref{EqnRho}))
\begin{itemize}
\item Define the subtrees $\cT_1$ and $\cT_2$ of $\cT$ as described above.
\item Calculate $\bP \left[\chi_{\cT_i}|(s(\rho_i)=\gamma_j) \right]$ for $i \in \{1,2 \}$, $j \in \{1,2 \}$, with the recursive formulas from above.
\item Evaluate $\bP \left[ (X_\gamma(e) = 1) | \chi \right]$ according to (\ref{EqnPx}) for all
$\gamma \in \{\alpha, \beta\} \times \{\alpha, \beta\}$.
\item Evaluate $\Psi_\gamma$ according to (\ref{EqnRho}) for all $\gamma \in \{\alpha, \beta\} \times \{\alpha, \beta\}$.
\end{itemize}
\begin{thm}
\textsc{EdgeLength} works correct, i.e. it returns $$\Psi_{\gamma}=\frac{\displaystyle \bE_\chi \left[ \sum_{e\in E,~e~\gamma-edge}(l(e)) \right]}
{\bE_\chi \left[ \#~of~\gamma-edges \right]}$$
The complexity is $O(|V|)$, so it is linear.
\end{thm}
\begin{proof}
The correctness of the algorithm follows from the construction above. It is left to verify the runtime.\\
Calculating the probabilities $\bP \left[(\chi_{\cT_i}|(s(\rho_i)=\gamma_j) \right]$ for $i \in \{1,2 \}$, $j \in \{1,2 \}$ with the recursive formulas requires
$O(|V|)$ calculations
since we have to evaluate one recursion formula for each vertex.
For each edge $e$, $\bP \left[ (X_\gamma(e) = 1) | \chi \right]$ can then be calculated according to (\ref{EqnPx})
with a constant number of calculations. So obtaining $\bP \left[ (X_\gamma(e) = 1) | \chi \right]$ for all $e$ requires $O(|E|)=O(|V|)$ calculations.
Calculating $\Psi_\gamma$ according to (\ref{EqnRho}) requires again $O(|E|)$ calculations. Therefore, the complexity is linear.
\end{proof}

%Define $\triangle_g^{\gamma} = av_{(e:~e~\gamma-edge)} (l(e))$, where l(e) is the length of $\gamma$-edge e in (\cT,g).\\
%The goal is to calculate $E_g[\triangle_g^{\gamma}]^{\gamma \in \{\alpha, \beta\}}$.\\
%The inverse of this expected value is the rate of speciation.\\
%\\
%Of course, $E_g[\triangle_f^{\gamma}]^{\gamma \in \{\alpha, \beta\}}$ = $\sum_g P(f|( g|_L = \chi))  \triangle_g^{\gamma}$.\\
%$P(g|( g|_L = \chi))$ can be calculated with the transition probabilities.\\
%But evaluating $P(g|( g|_L = \chi)) \triangle_g^{\gamma}$ for every possible g is very inefficient and can't be done for big trees.\\
%\\
%Our plan is to develop a polynomial algorithm to calculate $E_g[\triangle_g^{\gamma}]^{\gamma \in \{\alpha, \beta\}}$.
%The idea is to start in small subtrees and find $E_g$ with recursion formulas.
%

\chapter*{Outlook}
\addcontentsline{toc}{chapter}{Outlook}
There are several topics in the thesis which suggest further work.\\

In Chapter \ref{Martingale}, we conclude with the log-likelihood-ratio test for deciding if a tree evolved under Yule. The given bound for the power of the test, Equation (\ref{EqnPower}), depends on the bound for the Azuma inequality.
The bound $\ln n$ for the Azuma inequality was obtained in \ref{SecAzumaConst} by a lot of rough estimations. So we are very confident that there can be found a better bound $c \ln n$, with $c"<1$ being a constant. This would lead to an improved bound for the power of the log-likelihood-ratio test (i.e. one could show analytically that the log-likelihood-ratio test is very good even on trees with a small number of leaves).\\

The edge lengths estimation in Section \ref{EstEdgeLength} will be implemented by Rutger Vos in Perl for his library and in Java for Mesquite (Mesquite is a tree manipulation software suite). Once implemented, the algorithm can finally be applied to real data. One can then estimate the edge lengths of a constructed supertree.\\

%Further, it would be nice to code up the algorithm so that it can be applied to real data trees.
Section \ref{SpeciRate} provides an algorithm for calculating $\Psi_{\alpha,\alpha}$ and $\Psi_{\beta,\beta}$ which estimate the average edge lengths. Let $\psi_\alpha$ be the speciation rate for species of class $\alpha$ and let $\psi_\beta$ be the speciation rate for species of class $\beta$. One could test the hypothesis $\psi_{\alpha,\alpha}=\psi_{\beta,\beta}$ against
$\psi_{\alpha,\alpha} \neq \psi_{\beta,\beta}$ with the statistic $\frac{\Psi_{\alpha,\alpha}}{\Psi_{\beta,\beta}}$.
For evaluating this test, i.e. obtaining the Type I and Type II error, one can use simulations.\\

Further, in Section \ref{SpeciRate}, we assumed that the transition rates $r_\alpha$ and $r_\beta$ are given. An interesting open question is how to handle the problem without having these transition rates in advance.\\

\appendix
\chapter{List of Symbols}
%\begin{longtable}
%\caption{Table of Symbols}

\begin{longtable}{l l r}
%\allowdisplaybreaks
\emph{Symbol} & \emph{Meaning} & \emph{page}\\
\\
$\leq_T$ & partial order on the vertices of a tree $T$ & \pageref{DefPartOrder}\\
$\leq_{\cT}$ & partial order on the vertices of $\cT$ & \pageref{DefPartOrder}\\
$(2n-1)!!$ & $(2n-1) \times (2n-3) \ldots 3 \times 1$ & \pageref{Lem-n!!} \\
$(\cT,s)$ & phylogenetic state tree & \pageref{DefRank}\\
$(\cT,r)$ & ranked phylogenetic tree $\cT$ with rank function $r$ & \pageref{DefRank}\\
\\
$\alpha_{\cT,v}(i)$ & $|\{r: r(v)=i, r \in r(\cT) \}|$ & \pageref{DefiAlpha}\\
$\chi$ & character on a phylogenetic tree & \pageref{DefChar}\\
$\delta(v)$ & degree of vertex $v$ & \pageref{DefDegreeV}\\
$\lambda_v$ & number of elements of $\Vr$ that are descendants of $v$ & \pageref{LemNumbRank}\\
$\pi$ & initial probability distribution of Markov chain & \pageref{InitProb}\\
$\rho$ & root of a tree & \pageref{DefTree}\\
$\phi$ & labelling function of a phylogenetic tree $\cT$ & \pageref{DefPhyloTree}\\
$\Psi_\gamma$ & estimated length of a $\gamma$-edge & \pageref{CalcPsi}\\
\\
%$%\cP$ & poset & \pageref{DefPoset}\\
$\cT$ & phylogenetic $X$-tree & \pageref{DefPhyloTree} \\
$\cT_p$ & Primate supertree constructed in \cite{Vos2006} & \pageref{Primates}\\
$\cT_v$ & phylogenetic subtree of $\cT$ induced by vertex $v$ & \pageref{DefSubtreev}\\
$\cT_{X'}$ & phylogenetic subtree of $\cT$ with label set $X'$ & \pageref{DefSubtreeX'}\\
\\
$\bJ_p$ & Entropy of the probability distribution $p$ & \pageref{DefEntropy}\\
$\bP_{u<v}$ & Probability $\bP[r(u) < r(v) | \cT]$ & \pageref{ProbComp}\\
$\bP_U$ & Uniform distribution on $RB(X)$ & \pageref{ProbYU}\\
$\bP_U[\cT]$ & Probability of $\cT$ under the uniform model & \pageref{ProbYU}\\
$\bP_Y$ & Yule distribution on $RB(X)$ & \pageref{ProbYU}\\
$\bP_Y[\cT]$ & Probability of $\cT$ under the Yule model & \pageref{ProbYU}\\
\\
$c_n$ & Catalan number & \pageref{DefCatalan}\\
$C$ & set of character states & \pageref{DefChar}\\
$d(v)$ & number of direct descendants of vertex v & \pageref{DefNumbDes}\\
$d_{KL}(p,q)$ & Kullbach-Liebler distance between $p$ and $q$ & \pageref{DefKL}\\
$E, E_{\cT}$ & Edges of a phylogenetic tree $\cT$ & \pageref{DefGraph}\\
$l(e)$ & Length of edge $e$ in $\cT$ & \pageref{DefEdgeLength}\\
$L, L_{\cT}$ & Leaf set of a (phylogenetic) tree & \pageref{DefTree}\\
$p_{\gamma_1,\gamma_2}$ & probability of state change from $\gamma_1$ to $\gamma_2$ & \pageref{TransMat}\\
$P(l(e))$ & transition matrix of Markov chain,\\
 & dependent on edge length & \pageref{TransMat}\\
$r_\alpha (r_\beta)$ & rate of change from state $\alpha$ to $\beta$ ($\beta$ to $\alpha$) & \pageref{DefRateState}\\
$r, r_{\cT}$ & rank function of phylogenetic tree $\cT$ & \pageref{DefRank}\\
%$r(\cP)$ & Set of rank functions on $\cP$ & \pageref{RankPoset}\\
$r(\cT)$ & Set of rank functions on $\cT$ & \pageref{DefRank}\\
$rRB(n)$ & Set of ranked binary phylogenetic $X$-trees\\ & with $X=
\{1,2, \ldots n \}$ &
\pageref{TreeSetn} \\
$rRB(X)$ & Set of ranked binary phylogenetic $X$-trees & \pageref{RemTreeSet}\\
$R$ & rate matrix of a Markov chain & \pageref{RateMatrix}\\
$RB(n)$ & Set of binary phylogenetic $X$-trees\\ & with $X= \{1,2, \ldots n \}$ & \pageref{TreeSetn}\\
$RB(X)$ & Set of binary phylogenetic $X$-trees & \pageref{RemTreeSet}\\
$s$ & state function & \pageref{DefStateFunc}\\
$V, V_{\cT}$ & Set of vertices of a (phylogenetic) tree & \pageref{DefTree}\\
$\Vr, \Vr_{\cT}$ & Set of interior vertices of a (phylogenetic) tree & \pageref{DefTree}\\
\end{longtable}

\chapter{Algorithms coded in Python} \label{PythonCode}
\begin{verbatim}
# Rank functions
# Daniel Ford, Tanja Gernhard 2006
#
# Functions:
#
# rankprob(t,u) - returns the probability distribution 
#       of the rank of vertex "u" in tree "t"
# expectedrank(t,u) returns the expected rank
#       of vertex "u" and the variance
# compare(t,u,v) - returns the probability that "u" 
#       is below "v" in tree "t"



import random

# How we store the trees:
# The interior vertices of a tree with n leaves are 
#       labeled by 1...n-1
# Example input tree for all the algorithms below:
# The tree "t" below has n=9 leaves and the inner nodes have 
#       label 1...8
t1 = (((), (), {'leaves_below': 2, 'label': 4}), (), 
    {'leaves_below': 3, 'label': 3})
t2 = (((), (), {'leaves_below': 2, 'label': 7}), ((), (), 
    {'leaves_below': 2, 'label': 8}),
    {'leaves_below': 4, 'label': 6})
t3 = ((), (), {'leaves_below': 2, 'label': 5})
t4 = (t1,t3,{'leaves_below': 5, 'label': 2})
t = (t2,t4,{'leaves_below': 9, 'label': 1})


# Calculation of n choose j
# This version saves partial results for use later
nc_matrix = []	#stores the values of nchoose(n,j)
								# -- note: order of indices is reversed
def nchoose_static(n,j,nc_matrix):
    if j>n:
        return 0
    if len(nc_matrix)<j+1:
        for i in range(len(nc_matrix),j+1):
            nc_matrix += [[]]
    if len(nc_matrix[j])<n+1:
        for i in range(len(nc_matrix[j]),j):
           nc_matrix[j]+=[0]
        if len(nc_matrix[j])==j:
           nc_matrix[j]+=[1]    
        for i in range(len(nc_matrix[j]),n+1):
           nc_matrix[j]+=[nc_matrix[j][i-1]*i/(i-j)]
    return nc_matrix[j][n]

# dynamic programming verion
def nchoose(n,j):
    return nchoose_static(n,j,nc_matrix)	
    		#nc_matrix acts as a static variable


# get the number of descendants of u and of all vertices on the 
# path to the root (subroutine for rankprob(t,u))
def numDescendants(t,u):
    if t == ():
        return [False,False]
    if t[2]["label"]==u:
        return [True,[t[2]["leaves_below"]-1]]
    x = numDescendants(t[0],u)
    if x[0] == True:
        if t[1]==():
            n = 0
        else:
            n = t[1][2]["leaves_below"]-1
        return [True,x[1]+[n]]
    y = numDescendants(t[1],u)
    if y[0] == True:
        if t[0]==():
            n = 0
        else:
            n = t[0][2]["leaves_below"]-1
        return [True,y[1]+[n]]
    else:
        return [False,False]


# A version of rankprob which uses the function numDescendants
def rankprob(t,u):
    x = numDescendants(t,u)
    x = x[1]
    lhsm = x[0]
    k = len(x)
    start = 1
    end = 1
    rp = [0,1]
    step = 1
    while step < k:
        rhsm = x[step]
        newstart = start+1
        newend = end+rhsm+1
        rp2 = []
        for i in range(0,newend+1):
            rp2+=[0]
        for i in range(newstart,newend+1):
            q = max(0,i-1-end)
            for j in range(q,min(rhsm,i-2)+1):
                a = rp[i-j-1]*nchoose(lhsm + rhsm - (i-1),rhsm-j)
                 *nchoose(i-2,j)
                rp2[i]+=a
        rp = rp2
        start = newstart
        end = newend
        lhsm = lhsm+rhsm+1
        step+=1
    tot = float(sum(rp))
    for i in range(0,len(rp)):
        rp[i] = rp[i]/tot
    return rp


# For tree "t" and vertex "u" calculate the
# expected rank and variance
def expectedrank(t,u):
    rp = rankprob(t,u)
    mu = 0
    sigma = 0
    for i in range(0,len(rp)):
        mu += i*rp[i]
        sigma += i*i*rp[i]
    return (mu,sigma-mu*mu)


# GCD - assumes positive integers as input
# (subroutine for compare(t,u,v))
def gcd(n,m):
    if n==m:
        return n
    if m>n:
        [n,m]=[m,n]
    i = n/m
    n = n-m*i
    if n==0:
        return m
    return gcd(m,n)    


# Takes two large integers and attempts to divide them and give 
# the float answer without overflowing
# (subroutine for compare(t,u,v))
# does this by first taking out the gcd
def gcd_divide(n,m):
    x = gcd(n,m)
    n = n/x
    m = m/x
    return n/float(m)


# returns the subtree rooted at the common ancestor of u and v 
# (subroutine for compare(t,u,v))
# return
# True/False - have we found u yet
# True/False - have we found v yet
# the subtree - if we have found u and v
# the u half of the subtree
# the v half of the subtree
def subtree(t,u,v):
    if t == ():
        return [False,False,False,False,False]
    [a,b,c,x1,x2]=subtree(t[0],u,v)
    [d,e,f,y1,y2]=subtree(t[1],u,v)
    if (a and b):
        return [a,b,c,x1,x2]
    if (d and e):
        return [d,e,f,y1,y2]
    #
    x = (a or d or t[2]["label"]==u)
    y = (b or e or t[2]["label"]==v)
    #
    t1 = False
    t2 = False
    #
    if a:
	      t1 = x1
    if b:
        t2 = x2
    if d:
        t1 = y1
    if e:
        t2 = y2
    #
    if x and (not y):
	      t1 = t
    elif y and (not x):
        t2 = t
    #
    if t[2]["label"]==u:
        t1 = t
    if t[2]["label"]==v:
        t2 = t
    return [x,y,t,t1,t2]


# Gives the probability that vertex labeled v is 
# below vertex labeled u
def compare(t,u,v):
    [a,b,c,d,e] = subtree(t,u,v)
    if not (a and b):
        print "This tree does not have those vertices!"
        return 0
    if (c[2]["label"]==u):
        return 1.0
    if (c[2]["label"]==v):
        return 0.0
    tu = d
    tv = e    
    usize = d[2]["leaves_below"]-1
    vsize = e[2]["leaves_below"]-1
    x = rankprob(tu,u)      
    y = rankprob(tv,v)      
    for i in range(len(x),usize+2):
         x+=[0]
    xcumulative = [0]
    for i in range(1,len(x)):
        xcumulative+=[xcumulative[i-1]+x[i]]
    rp = [0]
    for i in range(1,len(y)):
        rp+=[0]
        for j in range(1,usize+1):
            a = y[i]*nchoose(i-1+j,j)*nchoose(vsize-i+usize-j,
              usize-j)*xcumulative[j]
    	    rp[i]+=a
    tot = nchoose(usize+vsize,vsize)
    return sum(rp)/float(tot)
\end{verbatim}

\chapter{Primate Supertree} \label{Primates}

\bigskip
%test
%
%\begin{figure}
\begin{center}
\includegraphics[scale=0.7]{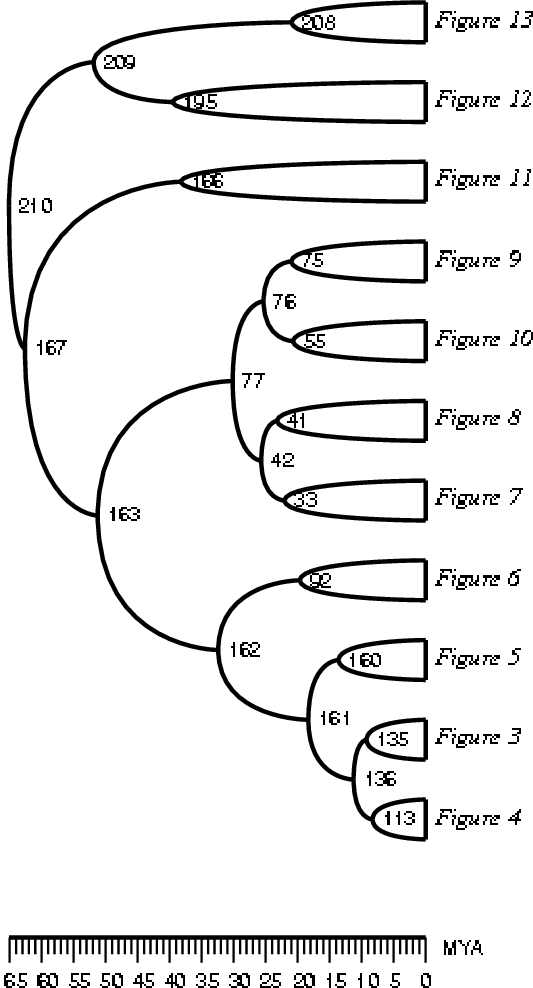}
%\caption{Primates - Figure 1}
%\label{Fig1}
\end{center}
%\end{figure}

\begin{figure}
\begin{center}
\includegraphics[scale=0.5]{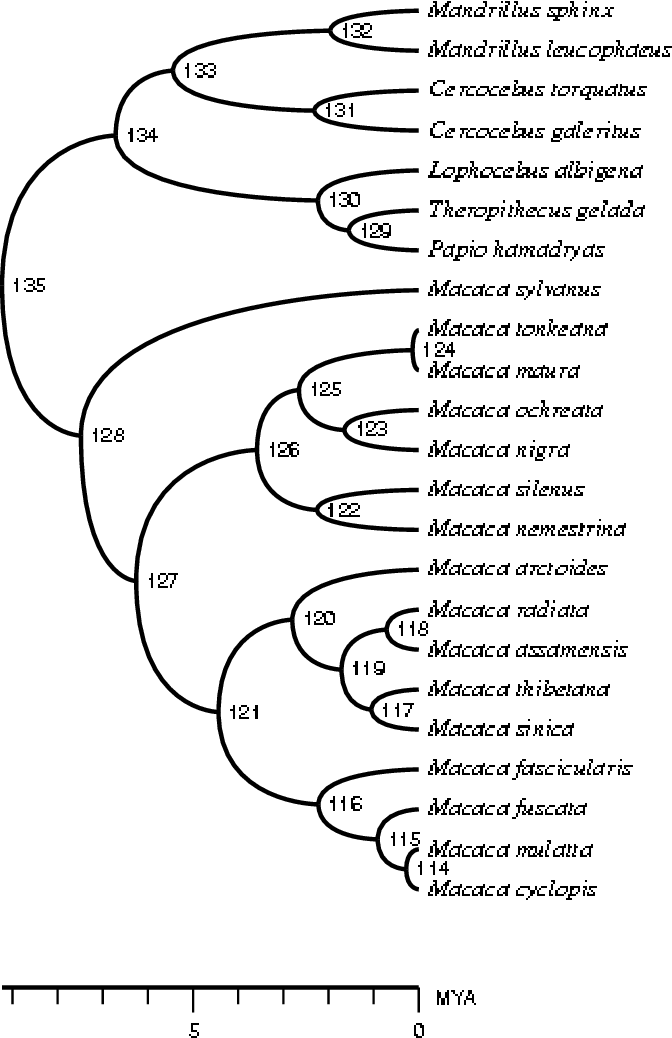}
\caption{Primate Supertree - Figure 3}
\label{Fig3}
\end{center}
\end{figure}

\begin{figure}
\begin{center}
\includegraphics[scale=0.5]{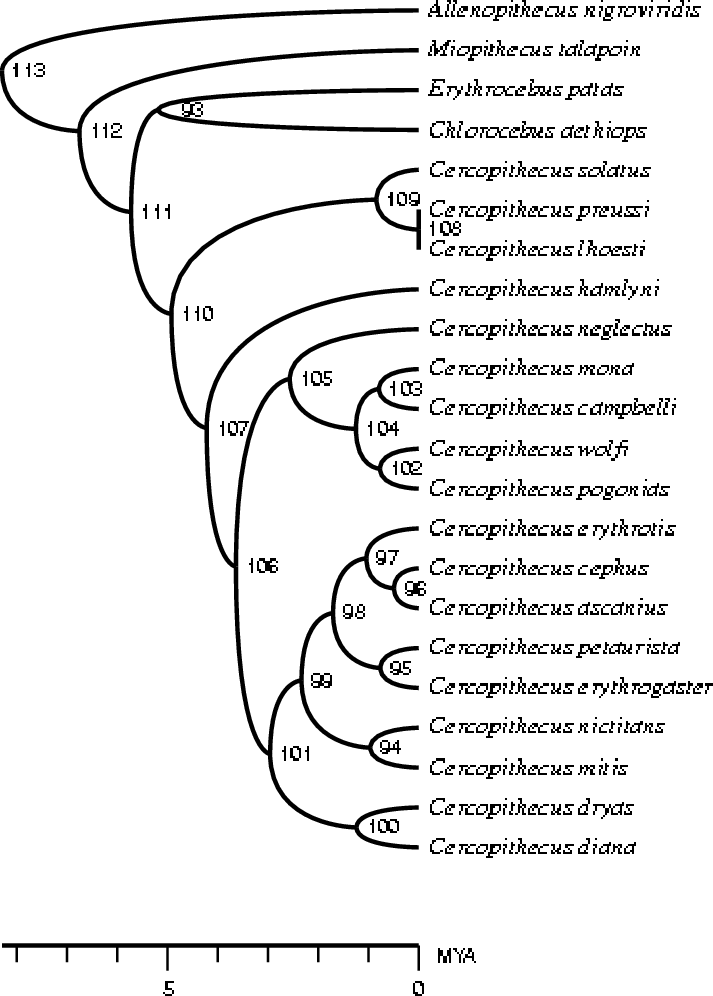}
\caption{Primate Supertree - Figure 4}
\label{Fig4}
\end{center}
\end{figure}

\begin{figure}
\begin{center}
\includegraphics[scale=0.5]{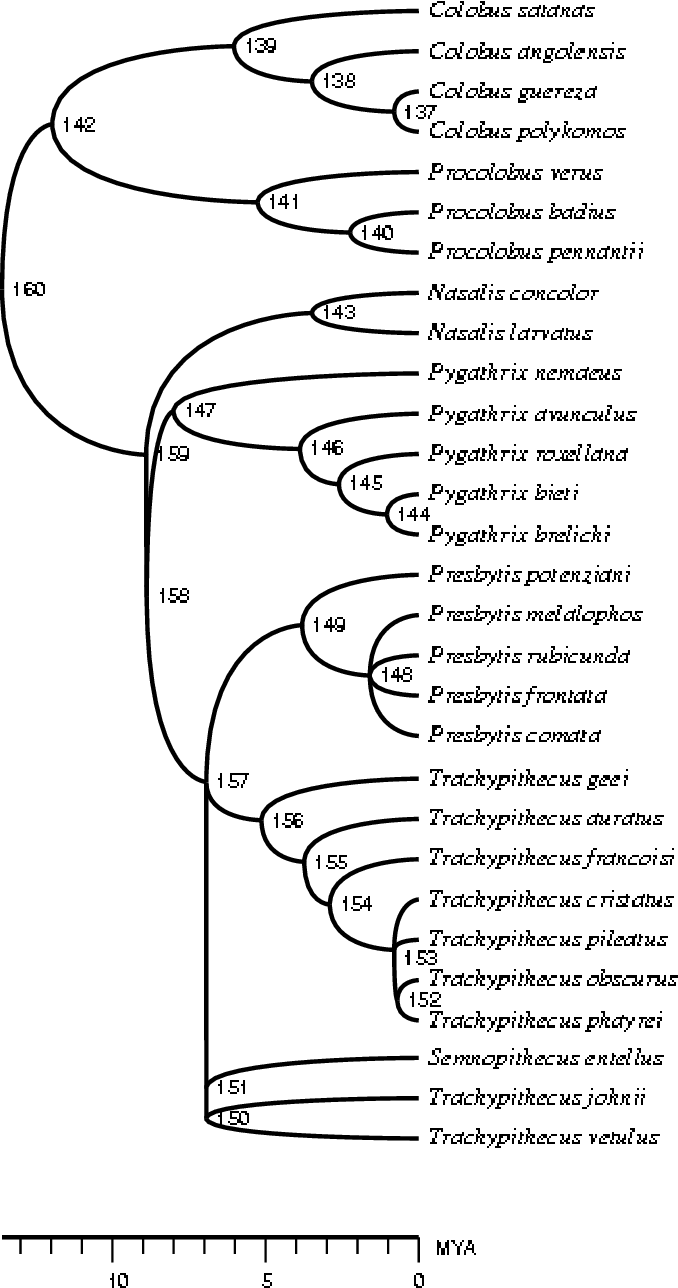}
\caption{Primate Supertree - Figure 5}
\label{Fig5}
\end{center}
\end{figure}

\begin{figure}
\begin{center}
\includegraphics[scale=0.5]{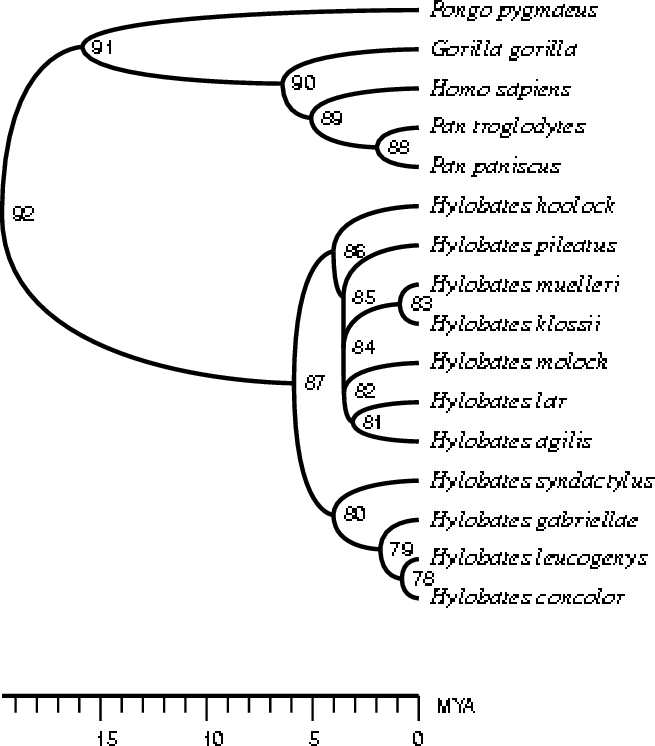}
\caption{Primate Supertree - Figure 6}
\label{Fig6}
\end{center}
\end{figure}

\begin{figure}
\begin{center}
\includegraphics[scale=0.5]{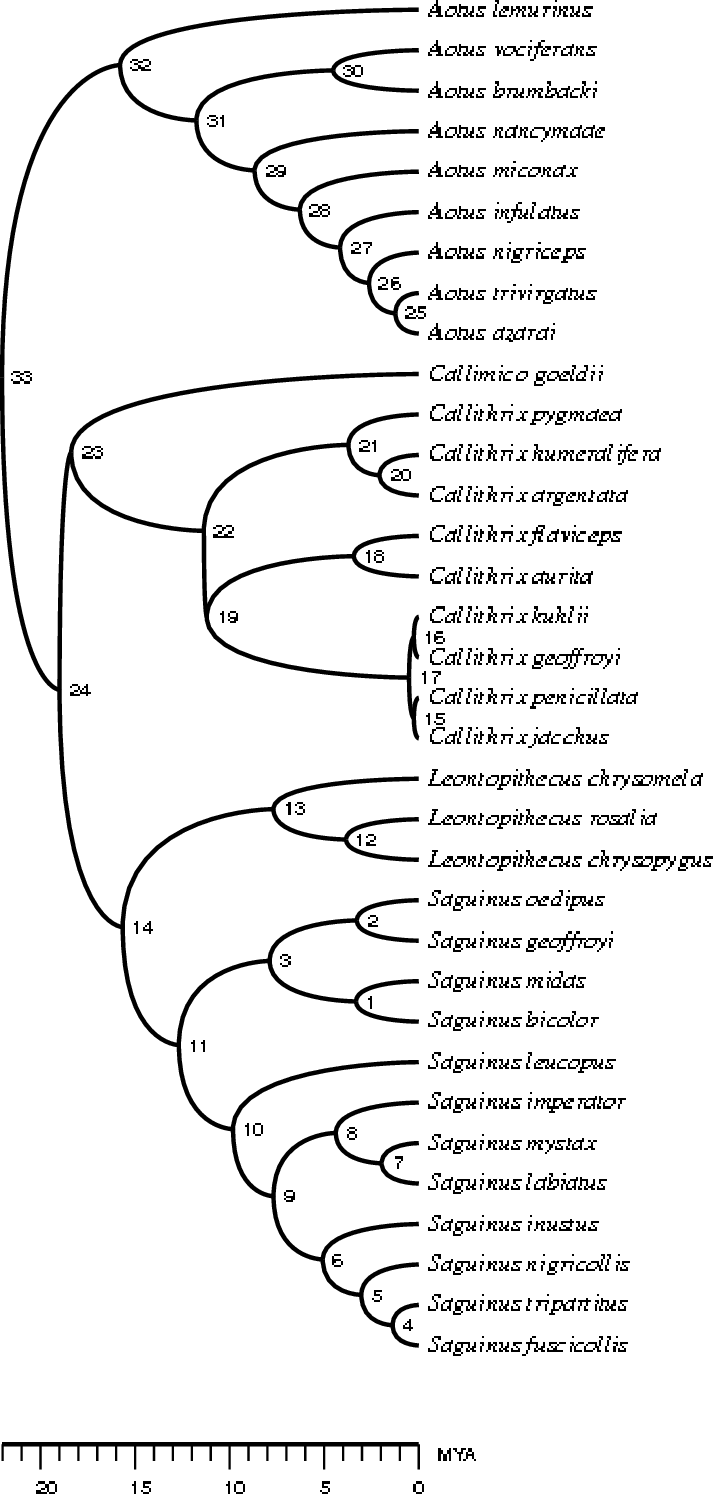}
\caption{Primate Supertree - Figure 7}
\label{Fig7}
\end{center}
\end{figure}

\begin{figure}
\begin{center}
\includegraphics[scale=0.5]{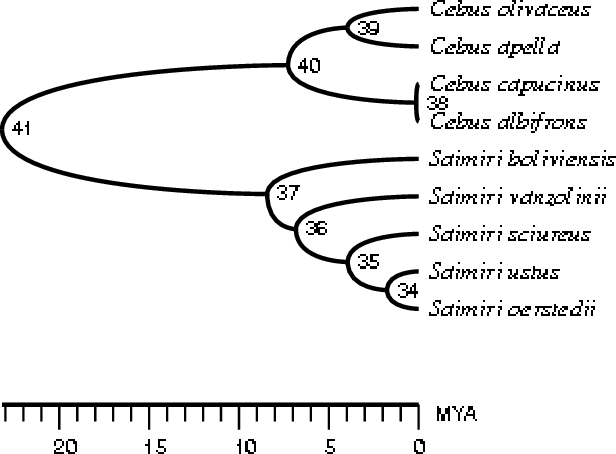}
\caption{Primate Supertree - Figure 8}
\label{Fig8}
\end{center}
\end{figure}

\begin{figure}
\begin{center}
\includegraphics[scale=0.5]{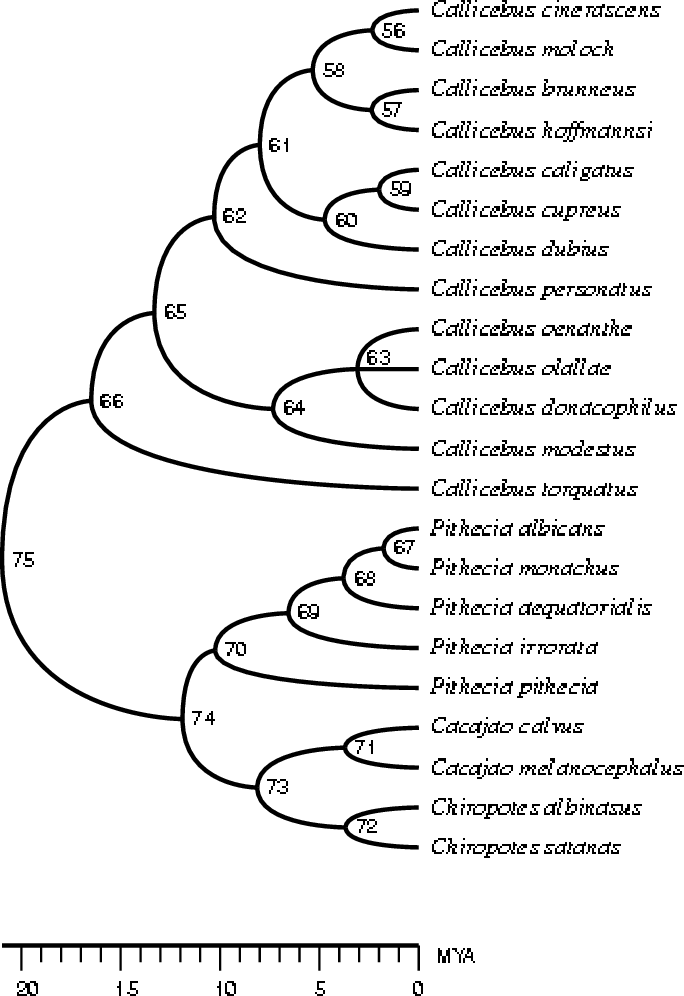}
\caption{Primate Supertree - Figure 9}
\label{Fig9}
\end{center}
\end{figure}

\begin{figure}
\begin{center}
\includegraphics[scale=0.5]{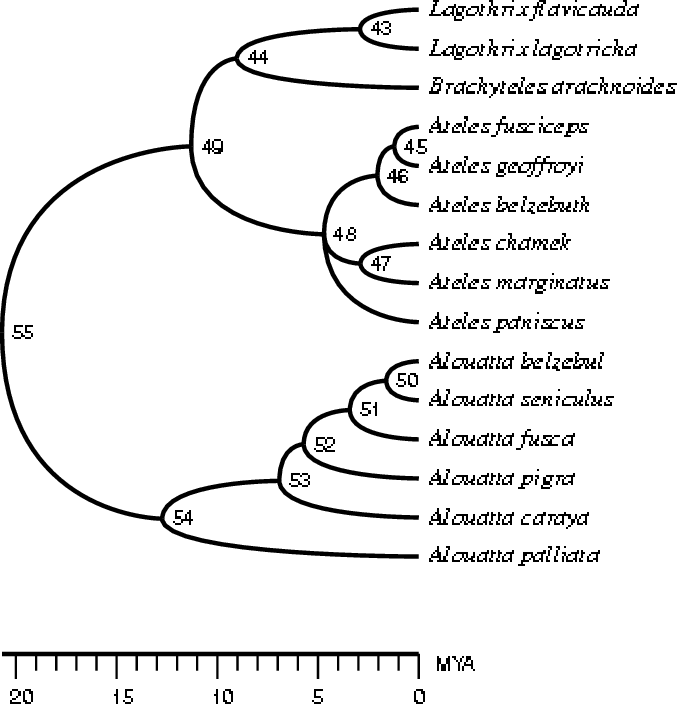}
\caption{Primate Supertree - Figure 10}
\label{Fig10}
\end{center}
\end{figure}

\begin{figure}
\begin{center}
\includegraphics[scale=0.5]{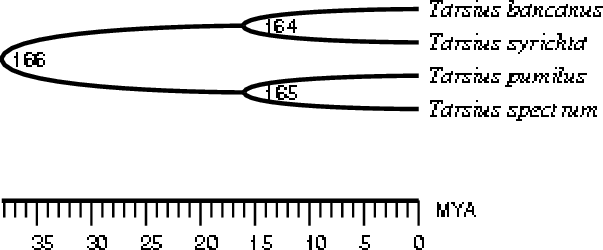}
\caption{Primate Supertree - Figure 11}
\label{Fig11}
\end{center}
\end{figure}

\begin{figure}
\begin{center}
\includegraphics[scale=0.5]{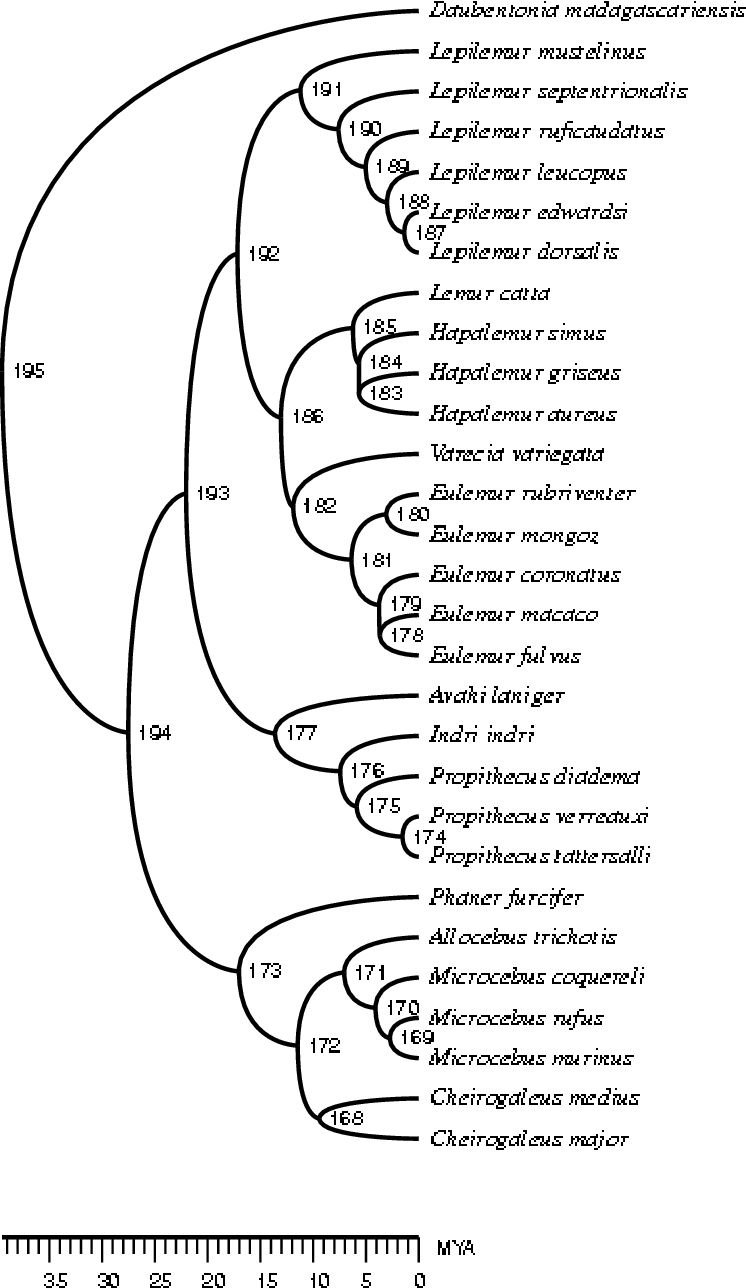}
\caption{Primate Supertree - Figure 12}
\label{Fig12}
\end{center}
\end{figure}

\begin{figure}
\begin{center}
\includegraphics[scale=0.5]{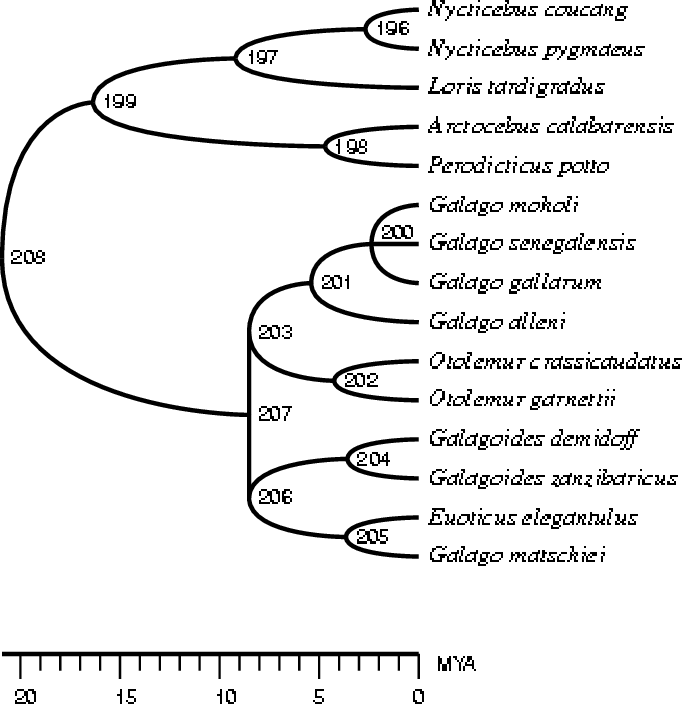}
\caption{Primate Supertree - Figure 13}
\label{Fig13}
\end{center}
\end{figure}

%\includegraphics[scale=1.0]{fig2.eps} 
%\includegraphics[scale=1.0]{test1.eps}  % oder.ps
%
%test
%\includegraphics{fig2_2.pdf} 
%test
%\includegraphics[bb= 0 0 100 100]{fig2_3.pdf} 
%test
%\includegraphics{fig2_4.pdf} 
%test
%\includegraphics{fig2_5.pdf} 
%test
%\includegraphics[scale=0.5]{test4.eps}
%\includegraphics[scale=0.7]{fig3.eps}

%\include{bibliography}
\bibliographystyle{abbrv} 
\bibliography{bibliography1}

\printindex

\end{document}